\documentclass[a4paper,11pt]{article}

\usepackage{amsmath,amssymb}
\usepackage[T1]{fontenc}
\usepackage{latexsym}
\usepackage[amsmath,thmmarks]{ntheorem}
\usepackage{enumerate}
\usepackage{fancyhdr}
\usepackage{theorem}
\usepackage{scrpage}
\usepackage{ucs}
\usepackage{latexsym}
\usepackage{dsfont}				
\usepackage{mathrsfs}			
\usepackage[active]{srcltx}		
\usepackage{color}
\usepackage{enumitem}

\DeclareMathOperator{\supp}{supp}
\DeclareMathOperator{\dist}{dist}
\DeclareMathOperator{\diam}{diam}
\DeclareMathOperator{\loc}{loc}
\DeclareMathOperator{\Lip}{Lip}

\numberwithin{equation}{section}

\newtheorem{Def}{Definition}[section]
\newtheorem{Lemma}[Def]{Lemma}
\newtheorem{Cor}[Def]{Corollary}
\newtheorem{Theorem}[Def]{Theorem}

\newtheorem{Prop}[Def]{Proposition}

{\theorembodyfont{\rmfamily} 
\newtheorem{Notation}[Def]{Notation}}
{\theorembodyfont{\rmfamily} \newtheorem{Remark}[Def]{Remark}}
{\theoremstyle{nonumberplain} }
{\theoremstyle{nonumberplain}
\theorembodyfont{\upshape}
\theoremheaderfont{\normalfont \bfseries}
\theoremsymbol{\ensuremath{_\square}}
\theoremseparator{:}
\newtheorem{Proof}{Proof}}

\let\Re=\relax
\DeclareMathOperator{\Re}{Re}

\newcommand{\R}{\ensuremath{\mathbb{R}}}
\newcommand{\C}{\ensuremath{\mathbb{C}}}
\newcommand{\N}{\ensuremath{\mathbb{N}}}
\newcommand{\Z}{\ensuremath{\mathbb{Z}}}

\newcommand{\bH}{\ensuremath{\mathbb{H}}}
\newcommand{\Eins}{\ensuremath{\mathds{1}}}

\newcommand{\calD}{\ensuremath{\mathcal{D}}}

\newcommand{\calE}{\ensuremath{\mathcal{E}}}
\newcommand{\calM}{\ensuremath{\mathcal{M}}}

\newcommand{\calC}{\ensuremath{\mathcal{C}}}
\newcommand{\calR}{\ensuremath{\mathcal{R}}}
\newcommand{\calQ}{\ensuremath{\mathcal{Q}}}

\newcommand{\calO}{\ensuremath{\mathcal{O}}}
\newcommand{\scrC}{\ensuremath{\mathscr{C}}}
\newcommand{\scrA}{\ensuremath{\mathscr{A}}}
\newcommand{\norm}[1]{\left\|#1\right\|}
\newcommand{\abs}[1]{\left|#1\right|}
\newcommand{\abbs}[1]{|#1|}

\newcommand{\skp}[1]{\langle #1 \rangle}

\newcommand{\eps}{\ensuremath{\varepsilon}}
\newcommand{\dis}{\displaystyle}
\newcommand{\molMeps}{\calM_0^{1,2,M,\eps}(L)}

\newcommand{\molMdeps}{\calM_0^{1,2,M,\eps}(L^{\ast})}

\title{A $T(1)$-Theorem for non-integral operators}
\author{Dorothee Frey and Peer Christian Kunstmann}

\begin{document}

\maketitle

\begin{abstract}
Let $X$ be a space of homogeneous type and let $L$ be a sectorial operator with bounded holomorphic functional calculus on $L^2(X)$. We assume that the semigroup $\{e^{-tL}\}_{t>0}$ satisfies Davies-Gaffney estimates. 
Associated to $L$ are certain approximations of the identity. We call an operator $T$ a non-integral operator if compositions involving $T$ and these approximations satisfy certain weighted norm estimates. 
 The Davies-Gaffney and the weighted norm estimates are together a substitute for the usual kernel estimates on $T$ in Calder\'on-Zygmund theory.\\
In this paper, we show, under the additional assumption that a vertical Littlewood-Paley-Stein square function associated to $L$ is bounded on $L^2(X)$, that a non-integral operator $T$ is bounded on $L^2(X)$ if and only if
$T(1) \in BMO_L(X)$ and $T^{\ast}(1) \in BMO_{L^{\ast}}(X)$.
Here, $BMO_L(X)$ and $BMO_{L^{\ast}}(X)$ denote the recently defined $BMO(X)$ spaces associated to $L$ that generalize the space $BMO(X)$ of John and Nirenberg. \\
Generalizing a recent result due to F. Bernicot, we show a second version of a $T(1)$-Theorem under weaker off-diagonal estimates, which gives a positive answer to a question raised by him.
As an application, we prove $L^2(X)$-boundedness of a paraproduct operator associated to $L$. We moreover study criterions for a $T(b)$-Theorem to be valid. 
\end{abstract}

\tableofcontents

\section{Introduction}

The term \emph{$T(1)$-Theorem} originally denotes a famous result of David and Journ\'{e} \cite{DavidJourne}, which characterizes the boundedness of Calder\'{o}n-Zygmund operators on $L^2(\R^n)$. 
In short, they prove that a Calder\'{o}n-Zygmund operator $T$ is bounded on $L^2(\R^n)$ if and only if it is weakly bounded (in some appropriate sense) and $T(1), T^{\ast}(1) \in BMO(\R^n)$.
 What is fascinating about this theorem is that it is both - a deep result of crucial importance and a theorem that can be formulated in only one sentence. \\
Many examples of operators, such as the Calder\'{o}n commutators and pseudo-differential operators, can be covered by this result. For others, such as the Cauchy integral operator along Lipschitz curves, the $T(1)$-Theorem is not directly applicable. This led to the development of a $T(b)$-Theorem, see \cite{McIntoshMeyer}, \cite{DavidJourneSemmes}, where the function $1$ is replaced by a para-accretive function $b$.
There exist numerous variants and generalizations, among them local $T(b)$-Theorems, see e.g. \cite{Christ2}, \cite{Hofmann}, generalizations to  non-homogeneous spaces, see e.g. \cite{NTV}, and operator-valued versions, see e.g. \cite{Hytoenen1}.\\

Even though in practice many operators fall under the scope of the Calder\'{o}n-Zygmund theory, there are still numerous operators of interest that do not.
It is well known, that an $L^2$-bounded Calder\'on-Zygmund operator is automatically also bounded on $L^p$ for all $p \in (1,\infty)$. This makes the Calder\'on-Zygmund theory not applicable to operators which are bounded on $L^p$ only for a range of $p$ strictly smaller than $(1,\infty)$.\\
Examples of such operators include operators that are related to a sectorial operator $L$ in $L^2$, with domain $\calD(L)$ and range $\calR(L)$, where the corresponding semigroup $\{e^{-tL}\}_{t>0}$ is bounded on $L^p$ only for a range of $p$ strictly smaller than $(1,\infty)$. In this case, one cannot work, as it has frequently been done in the last two decades, with pointwise Gaussian estimates for the semigroup, but has to work with generalized Gaussian estimates, Davies-Gaffney estimates or other off-diagonal estimates instead. \\
Aiming at a unified treatment of some of these operators, an $L^p$ theory was developed for operators that lie \emph{beyond} Calder\'{o}n-Zygmund theory, still - or even more - being ``singular'' in some sense and generalizing the concept of Calder\'{o}n-Zygmund operators. See e.g. \cite{DuongMcIntosh}, \cite{CoulhonDuong}, \cite{BlunckKunstmann}, \cite{ACDH}, \cite{AuscherD} and \cite{Auscher}. Actually, many ideas used in the  study of such operators are generalizations of methods developed in Calder\'{o}n-Zygmund theory.
  Those operators have also been called \emph{non-integral operators}, reflecting the property that the operators under consideration can no longer be represented as an integral operator with a Calder\'{o}n-Zygmund kernel, sometimes even not with any other kernel in a suitable sense (besides the Schwartz kernel). 
The main idea in this concept (already present in \cite{DuongRobinson}, \cite{DuongMcIntosh}) is to use approximation operators that are constructed via $H^{\infty}$-functional calculus as introduced in \cite{McIntosh}, e.g. the semigroup $\{e^{-tL}\}_{t>0}$ as an approximation of the identity and  the derivative $\{t\partial_te^{-tL}\}_{t>0}$ for the construction of a resolution of the identity. The H\"ormander condition for a Calder\'on-Zygmund operator is then replaced by weighted norm estimates, also called off-diagonal estimates, on compositions involving $T$ and these approximations.\\
Closely related to this theory are results on generalizations of operators and function spaces, that were originally constructed via the Laplacian and Littlewood-Paley theory. This includes versions of Hardy spaces $H^p_L$ and corresponding spaces $BMO_L$ that are associated to $L$, see e.g. \cite{AuscherDuongMcIntosh}, \cite{DuongYan2}, \cite{DuongYan3}, \cite{AuscherMcIntoshRuss}, \cite{BernicotZhao}, \cite{HofmannMayboroda}, \cite{HofmannMayborodaMcIntosh}, \cite{HLMMY}, \cite{DuongLi} and the study of Riesz transforms, e.g. in \cite{KatoSquare}, \cite{HofmannMartell}, \cite{BlunckKunstmann2}.\\

The present paper is devoted to a corresponding $L^2$ theory for such non-integral operators.\\
We assume $X$ to be a space of homogeneous type and let $L$ be a sectorial operator with bounded holomorphic functional calculus on $L^2(X)$. We assume that the semigroup $\{e^{-tL}\}_{t>0}$ satisfies Davies-Gaffney estimates, an $L^p-L^2$ estimate for some $p<2$ and an $L^2-L^q$ estimate for some $q>2$. 
Standard examples of operators that satisfy our assumptions are elliptic operators in divergence form with bounded measurable complex coefficients, see e.g. \cite{Auscher}, Schr\"odinger operators with singular potentials, see e.g. \cite{LSV}, and Laplace-Beltrami operators on complete Riemannian manifolds with non-negative Ricci curvature, see e.g. \cite{Davies1}, \cite{Grigoryan}.\\

Let us be a bit more precise on the term ``non-integral operator'': We consider operators $T: \calD(L) \cap \calR(L) \to L^2_{\loc}(X)$ with $T^{\ast}: \calD(L^{\ast}) \cap \calR(L^{\ast}) \to L^2_{\loc}(X)$ such that for functions $\psi_1,\psi_2 \in \Psi$ (where $\Psi$ denotes the   set consisting of bounded holomorphic functions on a sector with decay at zero and infinity) with suitable decay at zero the following off-diagonal estimates are valid:
\begin{align} \label{intro-eqT1} 
		\norm{T\psi_1(tL)f}_{L^2(B_2)} + \norm{T^{\ast}\psi_2(tL^{\ast})f}_{L^2(B_2)}
					&\leq C \left(1+\frac{\dist(B_1,B_2)^{2m}}{t} \right)^{-\gamma} \norm{f}_{L^2(B_1)}
\end{align}
for some $\gamma>0$, for all $t>0$, all balls $B_1,B_2$ with radius $r=t^{1/2m}$ and all $f \in L^2(X)$ supported in $B_1$.\\
These off-diagonal estimates replace H\"{o}lder or H\"{o}rmander conditions on the kernel of Calder\'{o}n-Zygmund operators. 
 Similar estimates were already used in \cite{HofmannMayboroda}, Theorem 3.2, to show boundedness of some operator $T: H_L^1(X) \to L^1(X)$ under the assumption that $T$ is bounded on $L^2(X)$.
 The relation of our assumptions on $T$ and those used in Theorem \ref{H1-bddness} is given by Lemma \ref{self-improving-est} below.
Moreover, observe that the estimates in  \eqref{intro-eqT1}  are not only ``off-diagonal'' assumptions, but also include the ``on-diagonal'' case for $\dist(B_1,B_2)=0$. In contrast to the standard $T(1)$-Theorem of \cite{DavidJourne}, we therefore do not require a weak boundedness property in Theorem \ref{T1-Theorem-intro} below. \\

On the Euclidean space $\R^n$ let us denote by $G_L$ the vertical Littlewood-Paley-Stein square function associated to $L$, i.e. let $\dis G_L(f)(x):=\left(\int_0^{\infty} \abs{t\nabla e^{-t^{2m}L} f(x)}^2 \,\frac{dt}{t} \right)^{1/2}$ for all $x \in \R^n$ and all $f \in L^2(\R^n)$. Then the main result, Theorem \ref{T1-Theorem} below, reads as follows:

\begin{Theorem} \label{T1-Theorem-intro}
		Let $L$ be the sectorial operator of order $2m$ as specified above such that $G_L$ and $G_{L^{\ast}}$ are bounded on $L^2(\R^n)$. Let $T$ be a non-integral operator satisfying \eqref{intro-eqT1} for sufficiently large $\gamma>0$. Then $T$ is bounded on $L^2(\R^n)$ if and only if
\begin{align*}
		T(1) \in BMO_L(\R^n) \qquad \text{and} \qquad T^{\ast}(1) \in BMO_{L^{\ast}}(\R^n).
\end{align*} 
\end{Theorem}

Here, $T(1)$ and $T^{\ast}(1)$ are appropriately defined linear functionals on a subspace of $H^1_L(\R^n)$ and $H^1_{L^{\ast}}(\R^n)$, respectively.
If the space $\R^n$ is replaced by an arbitrary space $X$ of homogeneous type, we require in addition the validity of some Poincar\'{e} inequality and have to reformulate the boundedness of the Littlewood-Paley-Stein square functions.\\

The most important tool in the proof of Theorem \ref{T1-Theorem-intro} is a paraproduct associated to $L$. With the help of a Fefferman-Stein criterion for Carleson measures and elements of $BMO_L(X)$, it is shown in \cite{Frey} (see Theorem \ref{paraproduct} below) that for every $b \in BMO_L(X)$ the operator
 \begin{equation} \label{intro-para}
 	\Pi_b: f \mapsto  \int_0^{\infty} \tilde{\psi}(t^{2m}L) [\psi(t^{2m}L)b \cdot A_t(e^{-t^{2m}L}f)] \, \frac{dt}{t}
 \end{equation}
is bounded on $L^2(X)$, where $\psi$, $\tilde{\psi} \in \Psi$ with sufficient decay at zero, and $A_t$ denotes some averaging operator. After subtracting $\Pi_{T(1)}$ from the operator $T$, the remaining term can be dealt with by Poincar\'e inequalities, quadratic estimates and almost orthogonality arguments. In absence of pointwise kernel estimates, off-diagonal estimates become a crucial technical ingredient, cf. Section \ref{sect-off-diag} below. \\

Let us also mention the following extension property of non-integral operators, that is shown in Corollary \ref{Lp-bddness} below. If $T$ satisfies \eqref{intro-eqT1} and is bounded on $L^2(X)$, then it extends to a bounded operator $T: H^p_L(X) \to L^p(X)$ for $p \in [1,2)$, $T: L^{p}(X) \to H^p_{L}(X)$ for $p \in (2,\infty)$ and $T: L^{\infty}(X) \to BMO_L(X)$.
Such a property is similar to the behaviour of Calder\'{o}n-Zygmund operators, in respect of the fact that every Calder\'{o}n-Zygmund operator, that is bounded on $L^2(X)$, is automatically also bounded on $L^p(X)$ for all $p \in (1,\infty)$.\\
For a second order elliptic operator $L$ in divergence form, we denote by $(p_{-}(L),p_{+}(L))$ and $(q_{-}(L),q_{+}(L))$ the interior of the interval of $L^p(X)$ boundedness of $\{e^{-tL}\}_{t>0}$ and $\{\sqrt{t}\nabla e^{-tL}\}_{t>0}$, respectively. In \cite{Auscher} it is shown that $p_{-}(L)= q_{-}(L)$ and $p_{+}(L) \geq q_{+}(L)$. Then for $p \in  (p_{-}(L),p_{+}(L))$, as shown in \cite{HofmannMayborodaMcIntosh}, there holds $H^p_L(X)=L^p(X)$, and therefore $T$ is bounded on $L^p(X)$ for all $p \in (p_{-}(L),p_{+}(L))$. For other types of operators $L$, one can obtain similar results via generalized Gaussian estimates, cf. Proposition \ref{Hp-equiv} below. 
However, these results on $L^p(X)$ boundedness also show that Theorem \ref{T1-Theorem-intro} is not applicable to operators, such as the Riesz transform $\nabla L^{-1/2}$, which are only bounded on $L^p(X)$ for $p \in (q_{-}(L),q_{+}(L))$, in the case that $p_{+}(L) > q_{+}(L)$.\\

While the work was in preparation, we learned that a similar $T(1)$-Theorem has also been proved by Bernicot, cf. \cite{Bernicot}. The main difference to our result is, that a crucial assumption in \cite{Bernicot} are pointwise bounds on the kernels of the semigroup $\{e^{-tL}\}_{t>0}$. Moreover, it is assumed that the conservation properties $e^{-tL}(1)=1$ and $e^{-tL^{\ast}}(1)=1$ hold. On the other hand, the assumed off-diagonal estimates on the operator $T$ are slightly weaker than \eqref{intro-eqT1}. \\
In Theorem \ref{weak-T1-Theorem} below, we show, with the same methods used in the proof of Theorem \ref{T1-Theorem-intro}, a second version of a $T(1)$-Theorem under weaker off-diagonal estimates. 
This generalizes the result of \cite{Bernicot} and answers the question of Bernicot, raised in \cite{Bernicot}, whether such a result could be obtained assuming only off-diagonal estimates instead of pointwise bounds on the kernel of the semigroup.\\

Under the additional assumption that $e^{-tL}$ is bounded on $L^{\infty}(X)$ uniformly in $t>0$, we then apply Theorem \ref{weak-T1-Theorem} to prove the boundedness of the paraproduct operator $\tilde{\Pi}_f$ on $L^2(X)$, where $\tilde{\Pi}_f$ is defined by
\begin{align*} 
		\tilde{\Pi}_f(g):= \int_0^{\infty} \psi(t^{2m}L) [e^{-t^{2m}L}g \cdot e^{-t^{2m}L}f ] \,\frac{dt}{t}
\end{align*}
for $f \in L^{\infty}(X)$, $g \in L^2(X)$ and $\psi \in \Psi$ with sufficient decay at zero and infinity.\\

Moreover, we study conditions for a $T(b)$-Theorem to be valid for an accretive function $b \in L^{\infty}(X)$. The conditions are given in terms of certain Schur conditions as they are used in a continuous version of the Cotlar-Knapp-Stein lemma. That is, we assume that for large enough $\delta>0$
\[
 	\norm{\tilde{\psi}(sL) M_b \psi(tL)}_{L^2(X)\to L^2(X)} 
		\leq C \min\left( \frac{s}{t}, \frac{t}{s} \right)^{\delta} \norm{b}_{L^{\infty}(X)}
\]
uniformly for all $s,t>0$, where $\psi,\tilde{\psi} \in \Psi$ and $M_b$ denotes the multiplication operator with $b$. \\

The article is organized as follows: In Section 2, we set some notation and summarize the most important definitions and results for spaces of homogeneous type and holomorphic functional calculus.
In Section 3, we define three different notions of off-diagonal estimates and collect essential properties of those. Moreover, we fix our assumptions on the operator $L$ and give a short review on Hardy and BMO spaces associated to operators. Section 4 contains the main results of this article. We first show how to define $T(1)$ and $T^{\ast}(1)$ (a problem which has not been addressed in \cite{Bernicot}) and introduce the notion of Poincar\'e inqualities on spaces of homogeneous type. We continue with statement and proof of a $T(1)$-Theorem for non-integral operators, Theorem \ref{T1-Theorem-intro}, and of a second version and explain the corresponding $L^p$ theory. Finally, we study criterions for a $T(b)$-Theorem.
In Section 5, we give the proofs of some auxiliary results concerning off-diagonal estimates for certain operators.\\

Throughout the article, the letter ``$C$'' will denote (possibly different) positive constants that are independent of the essential variables. We will frequently write $a \lesssim b$ for non-negative quantities $a,b$, if $a \leq Cb$ for some $C$.

\section{Preliminaries}

\subsection{Spaces of homogeneous type}
\label{sect-homogSpace}

In the following we will always assume $X$ to be a space of homogeneous type. 
More precisely, we assume that $(X,d)$ is a metric space and $\mu$ is a nonnegative Borel measure on $X$ with $\mu(X)=\infty$ which satisfies the \emph{doubling condition}:\\
 There exists a constant $A_1 \geq 1$ such that for all $x \in X$ and all $r>0$
 \begin{align*} 
 		V(x,2r) \leq A_1 V(x,r) < \infty,
 \end{align*}
 where we set $B(x,r):=\{y \in X\,:\, d(x,y)<r\}$, $V(U):=\mu(U)$ for an open set $U \subseteq X$ and $V(x,r):=\mu(B(x,r))$. 
 Note that the doubling property implies the following strong homogeneity property: There exists a constant $A_2>0$  and some $n>0$ such that for all $\lambda \geq 1$, for all $x \in X$ and all $r>0$
 \begin{align} \label{doublingProperty2}
 		V(x,\lambda r) \leq A_2 \lambda^n V(x,r).
 \end{align}
 In a Euclidean space with the Lebesgue measure, the parameter $n$ corresponds to the dimension of the space.
 There also exist constants $C$ and $D$, $0 \leq D \leq n$, so that
\begin{align} \label{DoublingProperty3}
 	V(y,r) \leq C \left(1+\frac{d(x,y)}{r}\right)^D V(x,r)
\end{align}
uniformly for all $x,y \in X$ and $r>0$.
For $D=n$, this is a direct consequence of \eqref{doublingProperty2} and the triangle inequality. If $X=\R^n$, then $D$ can be chosen to be $0$.
 For more details on spaces of homogeneous type, we refer to \cite{CoifmanWeiss2}.\\
We fix some element $x_0 \in X$ that is henceforth denoted by $0$. The ball $B_0:=B(0,1)$ is then referred to as \emph{unit ball}.\\
For a ball $B \subseteq X$ we denote by $r_B$ the radius of $B$ and set
\begin{equation} \label{annuli}
 	S_0(B):=B \qquad \text{and} \qquad S_j(B):= 2^jB \setminus 2^{j-1} B \quad \text{for} \; j=1,2,\ldots,
\end{equation}
where $2^jB$ is the ball with the same center as $B$ and radius $2^jr_B$. \\

We recall from \cite{David3}, \cite{Christ2} the following construction of an analogue of a dyadic grid on Euclidean spaces for spaces of homogeneous type.

 \begin{Lemma} \label{ChristCubes}
 Let $(X,d,\mu)$ be a space of homogeneous type. Then there exists a collection $\calQ := \{Q_{\alpha}^k  \,:\, k \in \Z, \, \alpha \in I_k\}$ of open subsets of $X$, where $I_k$ is some index set, a constant $\delta \in (0,1)$ and constants $C_1,C_2>0$ such that
 \begin{enumerate}[label=(\roman*)]
 	\item $\mu(X \setminus \bigcup_{\alpha} Q_{\alpha}^k) =0$ for each fixed $k$ and $Q_{\alpha}^k \cap Q_{\beta}^k = \emptyset$ if $\alpha \neq \beta$;
 	\item for any $\alpha,\beta,k,l$ with $l \geq k$, either $Q_{\beta}^l \subseteq Q_{\alpha}^k$ or $Q_{\beta}^l \cap Q_{\alpha}^k=\emptyset$; 
 	\item for each $l<k$ and each $\alpha \in I_k$ there is a unique $\beta$ such that $Q_{\alpha}^k \subseteq Q_{\beta}^l$;
 	\item $\diam(Q_{\alpha}^k) \leq C_1 \delta^k$;
 	\item each $Q_{\alpha}^k$ contains some ball $B(z^k_{\alpha},C_2\delta^k)$, where $z^k_{\alpha} \in X$.
 \end{enumerate}
  \end{Lemma}
 
 One can think of $Q^k_{\alpha} \in \calQ$ as being a dyadic cube with sidelength $\delta^k$ centered at $z^k_{\alpha}$.
By abuse of notation we will sometimes call the elements of the collection $\calQ$ ``cubes''.\\

We fix the following notation for further reference. It describes the covering of a dilated ball $2^jB$ with elements of $\calQ$ whose diameters are related to the radius of the ball $B$.

\begin{Notation} \label{Cube-Notation}
 Let $B=B(x_B,r_B)$ be an arbitrary ball in $X$. With the notation as in Lemma \ref{ChristCubes}, we define $k_0 \in \Z$ to be the uniquely determined integer satisfying
 \begin{align}	\label{Def-k0}
 			C_1 \delta^{k_0} \leq r_B < C_1 \delta^{k_0-1}
 \end{align} and for each $j \in \N$ we define $k_j \in \Z$ to be the integer satisfying
 \begin{align} \label{Def-kj}
 		 \delta^{-k_j} \leq 2^j < \delta^{-k_j-1}.
 \end{align}
 We further define for each $j \in \N$ the index set $M_j$ related to the ball $B=B(x_B,r_B)$ by 
\begin{align} \label{Def-Mj}
 	M_j:= \{ \beta \in I_{k_0} \,:\, Q_{\beta}^{k_0} \cap B(x_B,C_1 \delta^{k_0-k_j-2}) \neq \emptyset\},
\end{align}
representing all ``cubes'' out of $\calQ$ with ``sidelength'' approximately equal to $r_B$ that have non-empty intersection with the dilated ball $2^jB$. More precisely, we observe that Lemma \ref{ChristCubes} yields -- modulo null sets of $\mu$ -- for every $j \in \N$ the following inclusions:
\begin{align}  \label{inclusions-cubes}
	2^jB \subseteq B(x_B,C_1\delta^{k_0-k_j-2})
		\subseteq \bigcup_{\beta \in M_j} Q_{\beta}^{k_0}
		\subseteq B(x_B,2C_1\delta^{k_0-k_j-2})
		\subseteq \delta^{-2} 2^{j+1} B.
\end{align}
The first and the fourth inclusions are simple consequences of the definition of $k_0$ and $k_j$, whereas the second one follows from Lemma \ref{ChristCubes} (i) and the third one uses Lemma \ref{ChristCubes} (iv). 
Further, Lemma \ref{ChristCubes} yields that the sets $Q_{\beta}^{k_0}$, $\beta \in M_j$, are disjoint and for each $\beta \in M_j$ there exists some $z^{k_0}_{\beta} \in X$ such that
\begin{align}
		B(z^{k_0}_{\beta},c_1r_B) \subseteq Q_{\beta}^{k_0} \subseteq B(z_{\beta}^{k_0},r_B)
\end{align}
for some $c_1 \in (0,1)$ independent of $j$ and $\beta$ due to Lemma \ref{ChristCubes} (v) and (iv).
\end{Notation}

\begin{Remark} \label{Covering} 
The cardinality of the set $M_j$ defined in \eqref{Def-Mj} is bounded from above by a constant times $2^{jn}$.
This fact is in analogy to  the case of Euclidean spaces, meaning that for an arbitrary ball $B=B(x_B,r_B)$ in $X$, one can cover the dilated ball $2^jB=B(x_B,2^jr_B)$ by approximately $2^{jn}$ disjoint ``cubes'' out of $\calQ$ of diameter approximately equal to $r_B$. 
The argument is a simple modification of the one  given in \cite{CoifmanWeiss2}, Chapitre III. 
\end{Remark}

\subsection{Notation}

Let $f \in L^1_{\loc}(X)$. 
We denote the average of $f$ over an open set $U \in X$ by
\[
 	\skp{f}_U:= \frac{1}{V(U)} \int_U f(x) \,d\mu(x).
\]

\paragraph{Averaging operator}
Let $t>0$. With the notation as in Lemma \ref{ChristCubes}, we denote by $k_0 \in \Z$ the unique integer satisfying  \begin{align} \label{Def-k_0}
		C_1 \delta^{k_0} \leq t < C_1 \delta^{k_0-1}.
\end{align}
Then for almost every $x \in X$ there exists a unique $\alpha \in I_{k_0}$ such that $x \in Q_{\alpha}^{k_0}$.
We will therefore define the uncentered averaging operator $A_t$ with respect to ``dyadic cubes'' by
\begin{align} \label{avOp}
	A_tf(x) := \frac{1}{V(Q_{\alpha}^{k_0})} \int_{Q_{\alpha}^{k_0}} f(y) \,d\mu(y), \qquad \text{for almost all} \ x \in X,
\end{align}
for every $f \in L^1_{\loc}(X)$, where $Q_{\alpha}^{k_0}$ is the uniquely determined open set out of the collection $\{Q_{\beta}^{k_0}\}_{\beta \in I_{k_0}}$ with $x \in Q_{\alpha}^{k_0}$.
Observe that the operator $A_t$ is constant on each open set $Q_{\alpha}^{k_0}$ and that
$
		A_t f = \sum_{\alpha \in I_{k_0}} \skp{f}_{Q_{\alpha}^{k_0}} \Eins_{Q_{\alpha}^{k_0}},
$ 
where $k_0$ is determined by \eqref{Def-k_0}.
Moreover, there exists a constant $C>0$ such that for almost every $x \in X$ and every $f \in L^1_{\loc}(X)$
\begin{align} \label{avOp-est}
  \abs{A_tf(x)} \leq C \frac{1}{V(x,t)} \int_{B(x,t)} \abs{f(y)} d\mu(y).
\end{align}
\paragraph{Maximal operators}
We denote by $\calM$ the uncentered Hardy-Littlewood maximal operator, i.e. for a measurable function $f: X \to \C$ and a point $x \in X$ we set
\[
 	\calM f (x) = \sup_{\substack{r>0 \\ y \in B(x,r)}} \frac{1}{V(y,r)} \int_{B(y,r)} \abs{f(z)} \,d\mu(z).
\]
 Further, for $p \in [1,\infty)$, we denote by $\calM_p$ the \emph{p-maximal operator}, i.e. for a measurable function $f: X \to \C$ we set
$
 	\calM_pf = [\calM(\abs{f}^p)]^{1/p}. 
$
Recall that $\calM_p$ is bounded on $L^q(X)$ for every $q \in (p,\infty]$, but not on $L^p(X)$.

\paragraph{Tent spaces and Carleson measures} 
For any $x \in X$ , we denote by $\Gamma(x)$ the \emph{cone} of aperture $1$ with vertex $x$, namely
$
		\Gamma(x):=\{(y,t) \in X \times (0,\infty) \,:\, d(y,x) < t\}.
$
If $O$ is an open subset of $X$, then the \emph{tent} over $O$, denoted by $\hat{O}$, is defined as
$
		\hat{O}:= \{(x,t) \in X \times (0,\infty) \,:\, \dist(x,O^c) \geq t\}.
$
For any measurable function $F$ on $X \times (0,\infty)$, the conical square function $\scrA F$ is defined by 
\[
 	\scrA F(x) := \left(\iint_{\Gamma(x)} \abs{F(y,t)}^2 \; \frac{d\mu(y)}{V(x,t)}\frac{dt}{t} \right)^{1/2}, \qquad x \in X,
\]
and the Carleson function $\scrC F$ by
\[
 	\scrC F(x):= \sup_{B \,:\, x \in B} \left(\frac{1}{V(B)} \iint_{\hat{B}} \abs{F(y,t)}^2 \frac{d\mu(y)dt}{t}\right)^{1/2}, \qquad x \in X,
\]
where the supremum is taken over  all balls $B$ in $X$ that contain $x$. 
One then defines on $X \times (0,\infty)$ the tent spaces
\begin{align*}
 	T^1(X) &:= \{F: X \times (0,\infty) \to \C  \ \text{measurable} \,;\, \norm{F}_{T^1(X)}:=\norm{\scrA F}_{L^1(X)} < \infty \}, \\
	T^{\infty}(X) & := \{F: X \times (0,\infty) \to \C \ \text{measurable} \,;\, \norm{F}_{T^{\infty}(X)} := \norm{\scrC F}_{L^{\infty}(X)} < \infty\}.
\end{align*}
A \emph{Carleson measure} is a Borel measure $\nu$ on $X \times (0,\infty)$ such that 
\[
 	\norm{\nu}_{\calC} \, := \sup_B  \frac{1}{V(B)} \iint_{\hat{B}} \,\abs{d\nu} < \infty,
\]
where the supremum  is taken over all balls $B$ in $X$.\\
For more details on tent spaces and Carleson measures, we refer to \cite{CoifmanMeyerStein}.

\subsection{Holomorphic functional calculus}

We only state the most important definitions and results. For more details on holomorphic functional calculi we refer to \cite{McIntosh}, \cite{AlbrechtDuongMcIntosh}, \cite{LevicoNotes} and \cite{Haase}. \\
For $0 \leq \omega < \sigma < \pi$ we define the closed and open sectors in the complex plane $\C$ by
\begin{align*}
 	S_{\omega+} 	:= \{\zeta \in \C \setminus \{0\}\,:\, \abs{\arg \zeta} \leq \omega\} \cup \{0\}, \qquad 
	\Sigma^0_{\sigma} 		:= \{\zeta \in \C \,:\, \zeta \neq 0, \abs{\arg \zeta} < \sigma\}.
\end{align*}
We denote by $H(\Sigma^0_{\sigma})$ the space of all holomorphic functions on $\Sigma^0_{\sigma}$.
We further define 
\begin{align*}
 	H^{\infty}(\Sigma_{\sigma}^0)		&:= \{\psi \in H(\Sigma^0_{\sigma}) \,:\, \norm{\psi}_{L^{\infty}(\Sigma_{\sigma}^0) }< \infty \}, \\
	\Psi_{\alpha,\beta} (\Sigma^0_{\sigma}) 
				& := \{ \psi \in H(\Sigma^0_{\sigma}) \,:\, \exists C : \abs{\psi(\zeta)} \leq C \abs{\zeta}^\alpha(1+\abs{\zeta}^{\alpha+\beta})^{-1} \text{ for every } \zeta \in \Sigma_\sigma^0\}
\end{align*}
for every $\alpha,\beta>0$ and $\Psi(\Sigma_{\sigma}^0):= \bigcup_{\alpha,\beta>0} \Psi_{\alpha,\beta}(\Sigma_{\sigma}^0)$.

\begin{Def}
Let $\omega \in [0, \pi)$. A closed operator $L$ in a Hilbert space $H$ is said to be  \emph{sectorial of angle $\omega$}
if $\sigma(L) \subseteq S_{\omega+}$ and, for each $\sigma > \omega$, there exists a constant $C_{\sigma}>0$ such that
\[
 	\norm{(\zeta I-L)^{-1}} \leq C_{\sigma} \abs{\zeta}^{-1}, \qquad \zeta \notin S_{\sigma+}.
\]
\end{Def}

\begin{Remark}
Let  $\omega \in [0, \pi)$ and let $L$ be a sectorial operator of angle $\omega$ in a Hilbert space $H$. Then $L$ has dense domain in $H$. If $L$ is assumed to be injective, then $L$ also has dense range in $H$. See e.g. \cite{CDMcY}, Theorem 2.3 and Theorem 3.8.
\end{Remark}

Let $\omega < \theta < \sigma<\pi$ and let $L$ be a sectorial operator of angle $\omega \in [0,\pi)$ in a Hilbert space $H$. 
Then for every $\psi \in \Psi(\Sigma_{\sigma}^0)$
\begin{align} \label{Def-functcalc}
 	\psi(L):=\frac{1}{2\pi i} \int_{\partial \Sigma_{\theta}^0} \psi(\lambda) (\lambda I-L)^{-1} \,d\lambda
\end{align}
defines a bounded operator on $H$. 
 By sectoriality of $L$ the integral in \eqref{Def-functcalc} is well-defined, and an extension of Cauchy's theorem shows that the definition is independent of the choice of $\theta \in (\omega,\sigma)$.\\
Let $L$ be in addition injective and set $\psi(z):=z(1+z)^{-2}$. Then $\psi(L)$ is injective and has dense range in $H$. For $f \in H^{\infty}(\Sigma_{\sigma}^0)$
 one can define by
 \begin{align*}
 			f(L):= [\psi(L)]^{-1}(f \cdot \psi)(L) 
 \end{align*}
a closed operator in $H$. We say that $L$ has a \emph{bounded} $H^{\infty}(\Sigma_{\sigma}^0)$ \emph{functional calculus} if there exists a constant $c_{\sigma}>0$ such that for all $f \in H^{\infty}(\Sigma^0_{\sigma})$, there holds $f(L) \in B(H)$ with
\[
 	\norm{f(L)} \leq c_{\sigma} \norm{f}_{L^\infty(\Sigma_{\sigma}^0)}.
\]

One can show that $L$ has a bounded holomorphic functional calculus on $H$ if and only if the following quadratic estimates are satisfied:\\
For some (all) $\sigma \in (\omega,\pi)$ and some $\psi \in \Psi(\Sigma_{\sigma}^0)\setminus \{0\}$ there exists some $C>0$ such that for all $x \in H$ 
\begin{align} \label{square-functions}
		C^{-1} \norm{x}^2 \leq \int_0^{\infty} \norm{\psi(tL)x}^2 \,\frac{dt}{t} \leq C \norm{x}^2.
\end{align}

Moreover, if $\psi,\tilde{\psi} \in \Psi(\Sigma_{\sigma}^0)\setminus\{0\}$ are chosen to satisfy $\int_0^{\infty} \psi(t) \tilde{\psi}(t) \, \frac{dt}{t}=1$, then the functional calculus of $L$ on $H$ yields the following \emph{Calder\'{o}n reproducing formula}: For every $f \in H$ 
\begin{align*}
		\int_0^{\infty} \psi(t^{2m}L) \tilde{\psi}(t^{2m}L) f \,\frac{dt}{t} = f	\qquad \text{in}\ H.
\end{align*}
Observe that for given $\psi \in \Psi(\Sigma_{\sigma}^0)\setminus\{0\}$ and given $\alpha,\beta>0$, one can always find a function $\tilde{\psi}  \in \Psi_{\alpha,\beta}(\Sigma_{\sigma}^0)\setminus\{0\}$ such that $\int_0^{\infty} \psi(t) \tilde{\psi}(t) \, \frac{dt}{t}=1$.

\section{Off-diagonal estimates and definition of $H^p_L(X)$ and $BMO_L(X)$}

In the following, $m\geq 1$ will be a fixed constant, and $2m$ represents the order of the sectorial operator $L$.

\subsection{Davies-Gaffney and other off-diagonal estimates}
\label{sect-off-diag}

We introduce the following three different notions of off-diagonal estimates (compare also \cite{AuscherMartell2}).

\paragraph{Davies-Gaffney estimates} 
We say that the family of operators $\{S_t\}_{t>0}$ satisfies \emph{Davies-Gaffney estimates ($L^2$ off-diagonal estimates)} if there exist constants $C,c>0$ such that for arbitrary open sets $E,F \subseteq X$
\begin{equation} \label{L2offdiag}
 	\norm{S_tf}_{L^2(F)} \leq C e^{-\left(\frac{\dist(E,F)^{2m}}{ct}\right)^{\frac{1}{2m-1}}} \norm{f}_{L^2(E)},
\end{equation}
for every $t>0$ and every $f \in L^2(X)$ supported in $E$. 

\paragraph{Off-diagonal estimates}
We say that a family of operators $\{S_t\}_{t>0}$ satisfies \emph{$L^2$ off-diagonal estimates of order $\gamma$}, $\gamma>0$, if there exists a constant $C>0$ such that for arbitrary open sets $E, F \subseteq X$
\[
 	\norm{S_tf}_{L^2(F)} \leq C \left(1+\frac{\dist(E,F)^{2m}}{t} \right)^{-\gamma} \norm{f}_{L^2(E)},
\]
for every $t>0$ and every $f \in L^2(X)$ supported in $E$.

\paragraph{Weak off-diagonal estimates}
We say that a family of linear operators $\{S_t\}_{t>0}$ satisfies \emph{weak $L^2$ off-diagonal estimates of order $\gamma$}, $\gamma>0$, if there exists a constant $C>0$ such that for every $t>0$, arbitrary balls $B_1,B_2 \in X$ with radius $r=t^{1/2m}$ and every $f \in L^2(X)$ supported in $B_1$
\begin{equation} \label{weak-off-diag-est}
 	\norm{S_tf}_{L^2(B_2)} \leq C \left(1+\frac{\dist(B_1,B_2)^{2m}}{t} \right)^{-\gamma} \norm{f}_{L^2(B_1)}.
\end{equation}

Unless otherwise specified, we always mean by (weak) off-diagonal estimates the definition of (weak) $L^2$ off-diagonal estimates.\\

We collect some important properties of the different concepts of off-diagonal estimates. \\

Operator families that satisfy off-diagonal or Davies-Gaffney estimates are uniformly bounded on $L^2(X)$. This is a direct consequence of the definition by taking $E=F=X$. 
For operator families that satisfy weak off-diagonal estimates, we obtain the following.

\begin{Lemma} \label{uniformL2-bound}
 Assume that the family of operators $\{S_t\}_{t>0}$ satisfies weak $L^2$ off-diagonal estimates of order $\gamma>\frac{n}{2m}$.
Then $S_t$ is bounded on $L^2(X)$ uniformly in $t>0$, i.e. there exists a constant $C>0$ such that for all $f \in L^2(X)$ and every $t>0$
\[
 	\norm{S_tf}_{L^2(X)} \leq C \norm{f}_{L^2(X)}.
\]
 \end{Lemma}
 
 \begin{Proof}
Let $t>0$ and $f,g \in L^2(X)$. We split $X$ with the help of Lemma \ref{ChristCubes} into ``cubes'' out of $\calQ$ with diameter approximately equal to $t^{1/2m}$ and then order them into annuli around one fixed ``cube'' to get an estimate for the distance of the ``cubes''.
With the notation as in Lemma \ref{ChristCubes}, let $k_0 \in \Z$ be the integer satisfying $C_1 \delta^{k_0} \leq t^{1/2m} < C_1 \delta^{k_0-1}$.
In addition, for every $\alpha \in I_{k_0}$ we denote by $B_{\alpha}$ the ball $B(z_{\alpha}^{k_0},t^{1/2m})$ and observe that Lemma \ref{ChristCubes} (iv) and (v) yield the inclusion $Q_{\alpha}^{k_0} \subseteq B_{\alpha}$.
Then by assumptions 
\begin{align} \label{uniformL2-bound-eq1} \nonumber
 	\abs{\skp{S_tf,g}}
		& \leq \sum_{\alpha \in I_{k_0}} \sum_{\beta \in I_{k_0}} \abs{\skp{S_t \Eins_{Q_{\alpha}^{k_0}} f, \Eins_{Q_{\beta}^{k_0}}g}} \\ \nonumber
		& \lesssim  \sum_{\alpha \in I_{k_0}} \sum_{\beta \in I_{k_0}}
				\left(1+\frac{\dist(B_{\alpha},B_{\beta})^{2m}}{t}\right)^{-\gamma}
				\norm{f}_{L^2(Q_{\alpha}^{k_0})} \norm{g}_{L^2(Q_{\beta}^{k_0})} \\ \nonumber
		& \leq \left(\sum_{\alpha \in I_{k_0}} \sum_{\beta \in I_{k_0}} \left(1+\frac{\dist(B_{\alpha},B_{\beta})^{2m}}{t}\right)^{-\gamma} \norm{f}^2_{L^2(Q_{\alpha}^{k_0})} \right)^{1/2} \\
		& \qquad \qquad \times \left(\sum_{\alpha \in I_{k_0}} \sum_{\beta \in I_{k_0}} \left(1+\frac{\dist(B_{\alpha},B_{\beta})^{2m}}{t}\right)^{-\gamma} \norm{g}^2_{L^2(Q_{\beta}^{k_0})} \right)^{1/2}, 
\end{align}
using the Cauchy-Schwarz inequality in the last step.
Let $\alpha \in I_{k_0}$ be fixed and let $j \in \N$. As in Notation \ref{Cube-Notation} we define the index set $M_j$ related to the ball $B_{\alpha}$ by
\begin{align*}
 	M_j:= \{ \beta \in I_{k_0} \,:\, Q_{\beta}^{k_0} \cap B(z_{\alpha}^{k_0},C_1 \delta^{k_0-k_j-2}) \neq \emptyset\}.
\end{align*}
The inclusions \eqref{inclusions-cubes} from Notation \ref{Cube-Notation} yield that if $z_{\beta}^{k_0} \in S_j(B_{\alpha})$, then $\beta \in M_j$ and, by definition of the annulus, $\dist(B_{\alpha},B_{\beta}) \gtrsim 2^{j} t^{1/2m}$ for every $j \geq 3$.
We therefore get for fixed $\alpha \in I_{k_0}$
\begin{align} \label{uniformL2-bound-eq2} \nonumber
 	  \sum_{\beta \in I_{k_0}} \left(1+\frac{\dist(B_{\alpha},B_{\beta})^{2m}}{t}\right)^{-\gamma} 
	&	 \leq \sum_{j=0}^{\infty} \sum_{\substack{\beta \in I_{k_0} \\ z^{k_0}_{\beta} \in S_j(B_{\alpha})}}
				 \left(1+\frac{\dist(B_{\alpha},B_{\beta})^{2m}}{t}\right)^{-\gamma} \\
		&  \lesssim \sum_{j=0}^{\infty} \sum_{\beta \in M_j} (1+2^j)^{-2m\gamma} 
	\lesssim \sum_{j=0}^{\infty} 2^{-2m\gamma j} 2^{nj},
\end{align}
where we used the result of Remark \ref{Covering} in the last step, saying that the cardinality of $M_j$ is less than a constant times $2^{jn}$.
On the other hand, the disjointness of the cubes $\{Q_{\alpha}^{k_0}\}_{\alpha \in I_{k_0}}$ implies that
\begin{align*}
 	\sum_{\alpha \in I_{k_0}} \norm{f}^2_{L^2(Q_{\alpha}^{k_0})} \leq \norm{f}_{L^2(X)}^2.
\end{align*}
Hence, the expression in the first bracket of \eqref{uniformL2-bound-eq1} is bounded by a constant times $\norm{f}^2_{L^2(X)}$. Repeating the same procedure for the second bracket with the roles of $\alpha$ and $\beta$ interchanged and $f$ replaced by $g$ finally shows that $\abs{\skp{S_t f,g}} \lesssim \norm{f}_{L^2(X)} \norm{g}_{L^2(X)}$. 
\end{Proof}

\begin{Remark} \label{annuli-est}
  Let $\{S_t\}_{t>0}$ be a family of linear operators on $L^2(X)$ that satisfies weak off-diagonal estimates of order $\gamma>0$. 
 Then there exists a constant $C>0$ such that for an arbitrary ball $B \in X$ with radius $r_B=t^{1/2m}$, for all $j \in \N_0$ and all $f,g \in L^2(X)$ with $\supp f \subseteq B$ and $\supp g \subseteq S_j(B)$
\begin{align} \label{annuli-est2}
 	\abs{\skp{S_t f,g}} 
		\lesssim 2^{jn/2} \left(1+\frac{\dist(B,S_j(B))^{2m}}{t}\right)^{-\gamma} \norm{f}_{L^2(B)} \norm{g}_{L^2(S_j(B))}.
\end{align}
The proof is similar to the one of Lemma \ref{uniformL2-bound}.
\end{Remark}

\begin{Remark} \label{Def-on-Linfty}
Let $\{S_t\}_{t>0}$ be a family of linear operators on $L^2(X)$ that satisfies weak off-diagonal estimates of order $\gamma> \frac{n}{2m}$. Then, for every $t>0$ and every ball $B$ in $X$, the operator $S^{\ast}_t$ also acts from $L^2(B)$ to $L^1(X)$ and one can thus define $S_t$ as an operator from $L^{\infty}(X)$ to $L^2_{\loc}(X)$ via duality.
 This works as follows:\\
Let $f \in L^{\infty}(X)$ and $t>0$. Further, let $B=B(x_B,t^{1/2m})$ be some ball in $X$ and $\varphi \in L^2(X)$ with $\supp \varphi \subseteq B$.
Splitting $X$ into annuli around $B$, we obtain
\begin{align*}
 	\abs{\skp{f,S_t^{\ast}\varphi}}
		&\leq \sum_{j=0}^{\infty} \abs{\skp{f,\Eins_{S_j(B)}S_t^{\ast}\varphi}} \\
		&\lesssim \sum_{j=0}^{\infty} 2^{jn/2} \left(1+\frac{\dist(B,S_j(B))^{2m}}{t}\right)^{-\gamma} \norm{f}_{L^2(S_j(B))} \norm{\varphi}_{L^2(B)} \\
		& \lesssim  \norm{f}_{L^{\infty}(X)} \norm{\varphi}_{L^2(B)} V(B)^{1/2} \sum_{j=0}^{\infty} 2^{jn}(1+2^j)^{-2m\gamma} 
		\lesssim \norm{f}_{L^{\infty}(X)} \norm{\varphi}_{L^2(B)} V(B)^{1/2} ,
\end{align*}
since $\gamma>\frac{n}{2m}$.
Thus, for every $t>0$ we can define $S_t f$ for $f \in L^{\infty}(X)$ via duality as
\[
 	\skp{S_t f, \varphi} := \skp{f, S^{\ast}_t \varphi},
\]
where $\varphi \in L^2(X)$ is supported in some ball in $X$. 
\end{Remark}

We continue with another important observation:
All notions of off-diagonal estimates are stable under composition.
\begin{Lemma}  \label{DaviesGaffneyComp} 
Let $\{S_t\}_{t>0}$ and $\{T_t\}_{t>0}$ be two families of bounded linear operators on $L^2(X)$.\\
 (i) If $\{S_t\}_{t>0}$ and $\{T_t\}_{t>0}$ satisfy Davies-Gaffney estimates, then $\{S_sT_t\}_{s,t>0}$ satisfies Davies-Gaffney estimates in $\max(s,t)$.\\
(ii)  If  $\{S_t\}_{t>0}$ and $\{T_t\}_{t>0}$ satisfy off-diagonal estimates of order $\gamma$ and $\delta$, respectively, then $\{S_sT_t\}_{s,t>0}$ satisfies off-diagonal estimates of order $\min(\gamma,\delta)$ in $\max(s,t)$.
\end{Lemma}

Part (i) for $m=1$ is proven in \cite{HofmannMartell}, Lemma 2.3. The proof for arbitrary $m$ and part (ii) follow along the same lines. \\

For the corresponding result for families of operators that satisfy weak off-diagonal estimates (Proposition \ref{weak-composition} below), we first state some auxiliary results.

\begin{Remark} \label{CubeSummation2}
 Let $s,t>0$ with $t\leq s$ and let $B$ be an arbitrary ball in $X$ with radius $t$. As in Notation \ref{Cube-Notation}, let $k_0$ be the uniquely determined integer satisfying $C_1 \delta^{k_0} \leq t < C_1\delta^{k_0-1}$ and for each $\beta \in I_{k_0}$ let $B_{\beta}:=B(z_{\beta}^{k_0},t)$, where $z_{\beta}^{k_0}$ is given by Lemma \ref{ChristCubes}. Further, suppose that $\gamma>n$. Then  for every $\eps>0$ with $\gamma\geq n+\eps$
\begin{align*}
 	\sum_{\beta \in I_{k_0}} \left(1+\frac{\dist(B,B_{\beta})}{s}\right)^{-\gamma}
		 = \sum_{j=0}^{\infty} \sum_{\substack{\beta \in I_{k_0} \\ z_{\beta}^{k_0} \in S_j(B)}} 
					\left(1+\frac{\dist(B,B_{\beta})}{s}\right)^{-\gamma} 
		 \lesssim \sum_{j=0}^{\infty} \sum_{\beta \in M_j} \left(1+ \frac{2^jt}{s}\right)^{-(n+\eps)},
\end{align*}
using the fact that for every $j \geq3$ and all $\beta$ with $z_{\beta}^{k_0} \in S_j(B)$ there holds $\dist(B,B_{\beta}) \gtrsim 2^jt$ and $\beta \in M_j$, where $M_j$ was defined in \eqref{Def-Mj}. Moreover, Remark \ref{Covering} shows that $\#M_j \lesssim 2^{jn}$, therefore the above is bounded by a constant times
\begin{align*}
		\sum_{j=0}^{\infty} 2^{jn} \left(1+\frac{2^jt}{s}\right)^{-(n+\eps)}
				\lesssim \left(\frac{s}{t}\right)^{n+\eps} \sum_{j=0}^{\infty} 2^{jn} 2^{-j(n+\eps)}
					\lesssim 	\left(\frac{s}{t}\right)^{n+\eps},
\end{align*}
since we assumed $t \leq s$.\\
Thus, we finally obtain the following: For every $\eps>0$ there exists a constant $C>0$ such that for all $t \leq s$ and every $\gamma\geq n+\eps$
\begin{align} \label{CubeSummation-eq}
 \sum_{\beta \in I_{k_0}} \left(1+\frac{\dist(B,B_{\beta})}{s}\right)^{-\gamma}
			\leq C \left(\frac{s}{t}\right)^{n+\eps},
\end{align}
where $B$ is an arbitrary ball in $X$ with radius $t$ and the balls $B_{\beta}=B(z_{\beta}^{k_0},t)$ are specified above.
In view of the assumption $t \leq s$, one obviously aims at an application for sufficiently small chosen $\eps>0$.
\end{Remark}

Fundamental for the proof of Proposition \ref{weak-composition} is the following lemma.
It can be considered as an analogue of certain estimates for compositions of integral operators, see e.g. \cite{Grafakos}, Appendix K.1.
\begin{Lemma} \label{composition-lemma}
 Let $s,t>0$ with $t\leq s$ and let $B_1, B_2$ be two arbitrary balls in $X$ with radius $t$. If $\gamma, \delta > n$, then for every $\eps>0$ there exists some constant $C>0$ such that
\begin{align} \label{appK1}\nonumber
 	& \sum_{\beta \in I_{k_0}} \left(1+\frac{\dist(B_1,B_{\beta})}{s}\right)^{-\gamma}
				\left(1+\frac{\dist(B_{\beta},B_2)}{s}\right)^{-\delta} \\ 
			& \qquad \qquad \qquad \qquad \qquad \qquad \leq C \left(\frac{s}{t}\right)^{n+\eps} \left(1+\frac{\dist(B_1,B_2)}{s}\right)^{-\min(\gamma,\delta)},
\end{align}
where $B_{\beta}=B(z_{\beta}^{k_0},t)$, $k_0 \in \Z$ is uniquely determined by $C_1 \delta^{k_0} \leq t < C_1\delta^{k_0-1}$ and the index set $I_{k_0}$ and $z_{\beta}^{k_0}$ are given in Lemma \ref{ChristCubes}.
\end{Lemma}

\begin{Proof}
Let $\eps>0$. We denote by $\Sigma$ the left-hand side of \eqref{appK1}. If $\frac{\dist(B_1,B_2)}{s} \leq 3$, then we get, according to \eqref{CubeSummation-eq},
\[
 	\Sigma \leq \sum_{\beta \in I_{k_0}} \left(1+\frac{\dist(B_1,B_{\beta})}{s}\right)^{-\gamma} 
			\lesssim \left(\frac{s}{t}\right)^{n+\eps} 
			\lesssim \left(\frac{s}{t}\right)^{n+\eps} \left(1+\frac{\dist(B_1,B_2)}{s}\right)^{-\min(\gamma,\delta)}.
\]
If otherwise $\frac{\dist(B_1,B_2)}{s} \geq 3$, we split the space $X$ into two parts. For this purpose, we set $\rho:=\dist(B_1,B_2)$ and define $G:=\{x \in X\,:\, \dist(x,B_2) < \frac{\rho}{2}\}$. 
Then, for every $\beta \in I_{k_0}$ with $z_{\beta}^{k_0} \in G$ we have the estimate
\begin{align*}
	\dist(B_1,B_{\beta}) \geq \dist(B_1,G)-t \geq \frac{1}{2}\dist(B_1,B_2) - \frac{1}{3}\dist(B_1,B_2) = \frac{1}{6}\dist(B_1,B_2).
\end{align*}
Using \eqref{CubeSummation-eq}, this yields
\begin{align} \label{comp-lemma-eq1} \nonumber
	 &	\sum_{\substack{\beta \in I_{k_0} \\ z_{\beta}^{k_0} \in G}} 
									\left(1+\frac{\dist(B_1,B_{\beta})}{s}\right)^{-\gamma}
				\left(1+\frac{\dist(B_{\beta},B_2)}{s}\right)^{-\delta} \\ 
			&   \lesssim 	\left(1+\frac{\dist(B_1,B_2)}{s}\right)^{-\gamma}
					\sum_{\substack{\beta \in I_{k_0} \\ z_{\beta}^{k_0} \in G}} 
					\left(1+\frac{\dist(B_{\beta},B_2)}{s}\right)^{-\delta} 
			 \lesssim \left(\frac{s}{t}\right)^{n+\eps} 	\left(1+\frac{\dist(B_1,B_2)}{s}\right)^{-\min(\gamma,\delta)}.
\end{align}
Similarly, if $\beta \in I_{k_0}$ with $z_{\beta}^{k_0} \in X \setminus G$, we obtain the estimate $\dist(B_2,B_{\beta}) \gtrsim \dist(B_1,B_2)$. Hence, we can argue as before and end up with the same bound as in \eqref{comp-lemma-eq1} for the sum over all $\beta \in I_{k_0}$ with $z_{\beta}^{k_0} \in X \setminus G$. 
\end{Proof}

Now, we can show the following.
\begin{Prop} \label{weak-composition}
Let $\{S_t\}_{t>0}$ and $\{T_t\}_{t>0}$ be two families of linear operators on $L^2(X)$ that satisfy weak off-diagonal estimates of order $\gamma> \frac{n}{2m}$ and $\delta > \frac{n}{2m}$, respectively.
Then there exists some constant $C>0$  such that for every $t>0$ and arbitrary balls $B_1,B_2 \in X$ with radius $r=t^{1/2m}$
\[
 	\norm{S_t T_tf}_{L^2(B_2)} 
		\leq C\left(1+\frac{\dist(B_1,B_2)^{2m}}{t} \right)^{-\min(\gamma,\delta)} \norm{f}_{L^2(B_1)},
\]
for all $f \in L^2(X)$ supported in $B_1$.
\end{Prop}

\begin{Proof}
 Let $k_0 \in \Z$ be defined by \eqref{Def-k0}, so that $C_1\delta^{k_0} \leq t^{1/2m} < C_1\delta^{k_0-1}$. Moreover, let $I_{k_0}$ be the index set defined in Lemma \ref{ChristCubes} and denote for every $\beta \in I_{k_0}$ by $B_{\beta}$ the ball $B(z_{\beta}^{k_0},t^{1/2m})$. Lemma \ref{ChristCubes} then yields in particular that $X = \bigcup_{\beta \in I_{k_0}} B_{\beta}$ up to a nullset. \\
 Since we assumed $\gamma,\delta > \frac{n}{2m}$, we can apply Lemma \ref{composition-lemma} (now with $t^{1/2m}$ instead of $t$) and get for every $f \in L^2(X)$ with $\supp f \subseteq B_1$ by assumption on the operators
\begin{align*}
 	\norm{S_t T_tf}_{L^2(B_2)}
			& \leq \sum_{\beta \in I_{k_0}} \norm{S_t \Eins_{B_{\beta}} T_tf}_{L^2(B_2)} \\
			& \lesssim \sum_{\beta \in I_{k_0}} \left(1+\frac{\dist(B_2,B_{\beta})^{2m}}{t}\right)^{-\gamma}
					\left(1+\frac{\dist(B_{\beta},B_1)^{2m}}{t}\right)^{-\delta} \norm{f}_{L^2(B_1)} \\
			& \lesssim \left(1+\frac{\dist(B_1,B_2)^{2m}}{t}\right)^{-\min(\gamma,\delta)} \norm{f}_{L^2(B_1)}.
\end{align*}
\end{Proof}

One can apply weak off-diagonal estimates for balls with some radius different from the scale of the operator family in the following way.

\begin{Remark}	\label{weak-estimates-st}
 Let $\{S_t\}_{t>0}$ be a family of linear operators on $L^2(X)$ that satisfies weak off-diagonal estimates of order $\gamma>0$. 
Let $s,t>0$ and let $B_1,B_2$ be two arbitrary balls in $X$ with radius $s$. Then for every $f \in L^2(X)$ with $\supp f \subseteq B_1$
\[
 	\norm{S_t f}_{L^2(B_2)} \lesssim \max\left\{1, \left(\frac{s}{t^{1/2m}}\right)^n \right\} \left(1+\frac{\dist(B_1,B_2)^{2m}}{t}\right)^{-\gamma} \norm{f}_{L^2(B_1)}.
\]
For the proof, one splits $X$ into ``cubes'' out of $\calQ$ with diameter approximately equal to $t^{1/2m}$ and argues similar to the proof of Lemma \ref{uniformL2-bound}.
\end{Remark}

\subsection{Assumptions on $L$}
\label{sect-assumpOperator}

We fix our assumptions on the operator $L$. 
Unless otherwise specified, we will assume the following.

\begin{enumerate}[label= \textbf{(H\arabic*)},ref=H\theenumi]
\item \label{H1} 
The operator $L$ is an injective, sectorial operator in $L^2(X)$ of angle $\omega$, where $0 \leq \omega < \pi/2$. Further, $L$ has a bounded $H^{\infty}(\Sigma_\sigma^0)$ functional calculus for some (all) $\omega < \sigma < \pi$.
\item \label{H2}
	The operator $L$ generates an analytic semigroup $\{e^{-tL}\}_{t>0}$ satisfying Davies-Gaffney estimates, i.e. there exist constants $C,c>0$ such that for arbitrary open subsets $E,F \subseteq X$ 
\begin{equation} \label{DaviesGaffney}
		\norm{e^{-tL} f}_{L^2(F)}
			\leq C \exp\left[ - \left(\frac{\dist(E,F)^{2m}}{ct}\right)^{\frac{1}{2m-1}}\right] \norm{f}_{L^2(E)} 
\end{equation}
for every $t>0$ and every $f \in L^2(X)$ with $\supp f \subseteq E$.
\item \label{H3} 
	The semigroup $\{e^{-tL}\}_{t>0}$ satisfies an $L^{\tilde{p}}-L^2$ off-diagonal estimate for some $\tilde{p} \in (1,2)$ and an $L^2-L^{\tilde{q}}$ off-diagonal estimate for some $\tilde{q} \in (2,\infty)$, i.e. there exists a constant $C>0$ and some $\eps>0$ such that for every $t>0$, every $j \in \N_0$ and for an arbitrary ball $B$ in $X$ with radius $r=t^{1/2m}$
\begin{align} \label{Lp-L2-estimate}
 	\norm{e^{-tL} \Eins_{S_j(B)} f}_{L^2(B)} \leq C 2^{-j(\frac{n}{\tilde{p}}+\eps)} V(B)^{\frac{1}{2}-\frac{1}{\tilde{p}}} \norm{f}_{L^{\tilde{p}}(S_j(B))}
\end{align}
and
\begin{align} \label{L2-Lq-estimate}
 	\norm{e^{-tL} \Eins_B g}_{L^{\tilde{q}}(S_j(B))} \leq C 2^{-j(\frac{n}{\tilde{q}'}+\eps)} V(B)^{\frac{1}{\tilde{q}}-\frac{1}{2}} \norm{g}_{L^2(B)}	
\end{align}
for all $f \in L^{\tilde{p}}(X)$ and all $g \in L^2(X)$.
Here, $\tilde{q}'$ is the conjugate exponent of $\tilde{q}$ defined by $\frac{1}{\tilde{q}}+\frac{1}{\tilde{q}'}=1$.
\end{enumerate}
Observe that \eqref{L2-Lq-estimate} is just the dual estimate of \eqref{Lp-L2-estimate}. That is, if $L$ satisfies \eqref{L2-Lq-estimate} with exponent $\tilde{q}$, then $L^{\ast}$ satisfies \eqref{Lp-L2-estimate} with exponent $\tilde{q}'$ and vice versa.

\begin{Remark}
(i) One can show the following self-improving property of Davies-Gaffney estimates:\\
\emph{
Assume that \eqref{H1} is satisfied. If condition \eqref{DaviesGaffney} holds for all balls $B_1, B_2$ in $X$, then the assertion is also true for arbitrary open sets $E,F$ of $X$ (in general with different constants $C,c>0$).}\\
The proof is similar to the proof of Lemma \ref{uniformL2-bound} (cf. also \cite{AuscherMartell2}, Proposition 3.2(b)).\\
(ii) In the special case of non-negative self-adjoint operators $L$ and $m=1$, Coulhon and Sikora show in \cite{CoulhonSikora} that condition \eqref{DaviesGaffney} is equivalent to the following:\\
\emph{
There exist some constants $C,c>0$ such that for arbitrary open sets $E,F$ in $X$ with $\mu(E)<\infty$ and $\mu(F)<\infty$ and all $t>0$}
\[
 	\abs{\skp{e^{-tL} \Eins_{E},\Eins_{F}}} \leq C \exp\left(-\frac{\dist(E,F)^{2}}{ct}\right) \mu(E)^{1/2} \mu(F)^{1/2}.
\]
This is the form of Davies-Gaffney conditions as they were considered in \cite{Davies1}, for instance.
\end{Remark}

\begin{Remark}
If there exists a constant $C>0$ such that $V(x,r) \geq C r^n$ for all $x \in X$ and all $r>0$, then \eqref{H3} is a consequence of the following estimates:\\

Let $\tilde{p} \in (1,2)$ and $\tilde{q} \in (2,\infty)$. There exist constants $C,c>0$ such that for arbitrary open sets $E,F \subseteq X$
\begin{align} \label{Lp-L2-estimate-2}
 	\norm{e^{-tL} f}_{L^2(F)} \leq C t^{-\frac{n}{2m}(\frac{1}{\tilde{p}}-\frac{1}{2})} \exp\left[-\left(\frac{\dist(E,F)^{2m}}{ct}\right)^{\frac{1}{2m-1}}\right] \norm{f}_{L^{\tilde{p}}(E)}	
\end{align}
and
\begin{align} \label{L2-Lq-estimate-2}
 	\norm{e^{-tL} g}_{L^{\tilde{q}}(F)} \leq C t^{-\frac{n}{2m}(\frac{1}{2}-\frac{1}{\tilde{q}})} \exp\left[-\left(\frac{\dist(E,F)^{2m}}{ct}\right)^{\frac{1}{2m-1}}\right] \norm{g}_{L^2(E)}	
\end{align}
for every $t>0$ and every $f \in L^{\tilde{p}}(X)$ and $g \in L^2(X)$ supported in $E$.\\

The proof is obvious. If \eqref{Lp-L2-estimate-2} is satisfied, then, in particular, $e^{-tL}: L^{\tilde{p}}(X) \to L^2(X)$ is bounded for every $t>0$. Analogously, if \eqref{L2-Lq-estimate-2} is satisfied, then $e^{-tL}: L^2(X) \to L^{\tilde{q}}(X)$ is bounded for every $t>0$. For sufficient conditions for \eqref{Lp-L2-estimate-2} in terms of off-diagonal estimates of annular type, we refer to \cite{AuscherMartell2}, Proposition 3.2.
We refer to \cite{BlunckKunstmann4} and \cite{AuscherMartell2} in general for further comparison of these types of off-diagonal estimates. 
\end{Remark}


One can show that the Davies-Gaffney estimates imply $L^2$ off-diagonal estimates for more general operator families associated to $L$. The proof of \cite{HofmannMayborodaMcIntosh}, Lemma 2.28, carries over with only minor changes to our more general setting.

\begin{Prop}  \label{H-inf-offdiag}
 Let $L$ satisfy \eqref{H1} and \eqref{H2}. Let $\sigma \in (\omega, \frac{\pi}{2}), \ \psi \in \Psi_{\alpha,\beta}(\Sigma_{\sigma}^0)$ for some $\alpha, \beta >0$. Then, for any $\varphi \in H^{\infty}(\Sigma_{\sigma}^0)$, the family of operators $\{\psi(tL)\varphi(L)\}_{t>0}$ satisfies $L^2$ off-diagonal estimates of order $\alpha$, with the constant controlled by $\norm{\varphi}_{L^{\infty}(\Sigma_{\sigma}^0)}$.
\end{Prop}

The following almost orthogonality lemma is a slight generalization of \cite{HofmannMayborodaMcIntosh}, Lemma 4.6. 

 \begin{Lemma} \label{off-diag-st}
Let $\sigma \in (\omega,\frac{\pi}{2})$, $\alpha,\beta>0$ and $\psi \in \Psi_{\alpha,\beta}(\Sigma_{\sigma}^0)$.
Let further $\delta>0$ and $\varphi \in H^{\infty}(\Sigma_{\sigma}^0)$ with $\varphi(z) =\calO(\abs{z}^{\delta})$ for $\abs{z} \to 0$. Then for any $a \geq 0$ with $a \leq \delta$ and $a<\beta$, there is a family of operators $\{T_{s,t}\}_{s,t>0}$ such that
 \[
  		\varphi(tL) \psi(sL) = \left(\frac{t}{s}\right)^a T_{s,t}, \qquad s,t>0,
 \]
 where $\{T_{s,t}\}_{s,t>0}$ satisfies $L^2$ off-diagonal estimates in $s$ of order $\alpha + a$ uniformly in $t>0$.
 \end{Lemma}

 \begin{Proof}
 Let $\psi$, $\varphi$ as given in the assumptions and let $s,t>0$. For every $a>0$ with $a \leq \delta$ and $a<\beta$ we write 
 \begin{align*}
  	& \varphi(tL) \psi(sL) 
 		 = \left(\frac{t}{s}\right)^a (tL)^{-a} \varphi(tL) (sL)^a \psi(sL) 
 		  = \left(\frac{t}{s}\right)^a T_{s,t}		
 \end{align*}
  with $T_{s,t}:=(tL)^{-a}\varphi(tL) (sL)^a \psi(sL)$. 
 Since we assumed $\delta \geq a$ and $\varphi \in H^{\infty}(\Sigma_{\sigma}^0)$, there exists a constant $C>0$ such that for every $z \in \Sigma_{\sigma}^0$ with $\abs{z} \leq 1$ there holds
 $
  	\abs{z^{-a}\varphi(z)} \leq c\abs{z}^{-a} \abs{z}^\delta \leq C
 $
 and, obviously, for every $z \in \Sigma_{\sigma}^0$ with $\abs{z} \geq 1$ also
 $
  	\abs{z^{-a}\varphi(z)} \leq C.
 $
 Hence, the function $z \mapsto z^{-a}\varphi(z)$, $z \in \Sigma_{\sigma}^0$, belongs to $H^{\infty}(\Sigma_{\sigma}^0)$ with
 \[
  	\sup_{t>0} \norm{(t\,\cdot\,)^{-a}\varphi(t\,\cdot\,)}_{L^{\infty}(\Sigma_{\sigma}^0)} \leq C.
 \]
 As the function $z \mapsto z^a \psi(z)$ is in $\Psi_{\alpha+a,\beta-a}(\Sigma_{\sigma}^0)$,
 Proposition \ref{H-inf-offdiag} yields that
 $\{T_{s,t}\}_{s,t>0}$ satisfies $L^2$ off-diagonal estimates in $s$ of order $\alpha+a$ uniformly in $t >0 $.
For $a=0$ the claim follows directly from Proposition \ref{H-inf-offdiag}.
 \end{Proof}

We end with an observation on conservation properties of the semigroup. For a proof, we refer to \cite{Frey}, Lemma 2.9.

\begin{Lemma} \label{psi-remark}
Let $L$ satisfy \eqref{H1}, \eqref{H2} and let $\sigma \in (\omega, \frac{\pi}{2})$.\\
(i) Let $\gamma>\frac{n}{4m}$. For every ball $B \subseteq X$ there exists some constant $C_B>0$ such that for all $t>0$ 
\begin{align*}
		\norm{e^{-tL^{\ast}}}_{L^2(B)\to L^1(X\setminus 4B)} \leq C_B t^{\gamma}.
\end{align*}
In particular, one can define $e^{-tL}$ via duality as an operator from $L^{\infty}(X)$ to $L^2_{\loc}(X)$.\\
(ii) Let $\alpha>0$, $\beta>\frac{n}{4m}$ and $\psi \in \Psi_{\beta,\alpha}(\Sigma_{\sigma}^0)$. Moreover, let $b \in L^{\infty}(X)$. 
If for every $t>0$ 
\[
 	e^{-tL}(b)=b  \qquad \text{in} \ L^2_{\loc}(X),
\]
then for every $t>0$
\[
 	\psi(tL)(b)=0 \qquad \text{in} \ L^2_{\loc}(X).
\]
\end{Lemma}

 \subsection{Hardy and BMO spaces associated to operators}
 \label{sect-HardySpaces}

In the following, we will always assume that the operator $L$ satisfies the assumptions \eqref{H1} and \eqref{H2} and that $\sigma \in (\omega,\frac{\pi}{2})$.\\
We summarize the most important facts about Hardy and BMO spaces associated to $L$. For more details and proofs of the results, we refer to  \cite{HofmannMayboroda}, \cite{HofmannMayborodaMcIntosh}, \cite{HLMMY} and \cite{DuongLi}. 
The proofs given there carry over with only minor changes to our more general setting. \\


For  $\psi \in \Psi(\Sigma_{\sigma}^0)\setminus \{0\}$ and $(x,t) \in X \times (0,\infty)$ we set $Q_{\psi,L}f(x,t):=\psi(t^{2m}L)f(x)$.

\begin{Def} \label{Def-Hp}
(i) Let $1 \leq p \leq 2$ and let $\psi_0 \in \Psi(\Sigma_{\sigma}^0)$ be defined by $\psi_0(z):=ze^{-z}$. Define $H^p_L(X)$  to be the completion of the space
\begin{equation} \label{Def-bH1}
		\bH^p_L(X):= \{f \in L^2(X) \,:\, \scrA Q_{\psi_0,L}f \in L^p(X)\},
\end{equation}
 with respect to the norm
$
 			\norm{f}_{H^p_{\psi_0,L}(X)} := \norm{\scrA Q_{\psi_0,L}f}_{L^p(X)} =\norm{Q_{\psi,L}f}_{T^p(X)}. 
$\\
(ii) Let $2<p<\infty$. Define
$
			H^p_L(X) :=(H^{p'}_{L^{\ast}}(X))',
$
where $\frac{1}{p}+\frac{1}{p'}=1$ and $L^{\ast}$ is the adjoint operator of $L$.
\end{Def}

Observe that $H^2_L(X)=L^2(X)$ by \eqref{H1} and known square function estimates. \\

In both cases, for $p \leq 2$ and for $p>2$, there is a characterization of $H^p_L(X)$ by general square functions constructed via functions $\psi \in \Psi(\Sigma_{\sigma}^0) \setminus \{0\}$ with a certain decay at infinity and at zero, respectively. 
For a proof, we refer to Corollary 4.21 of \cite{HofmannMayborodaMcIntosh}.

\begin{Theorem} \label{Charact-Hp}
	Let $\alpha>0$ and $\beta>\frac{n}{4m}$. Further, let either $1\leq p \leq 2$ and $\psi \in \Psi_{\alpha,\beta}(\Sigma_{\sigma}^0) \setminus \{0\}$ or $2 \leq p < \infty$ and $\psi \in \Psi_{\beta,\alpha}(\Sigma_{\sigma}^0) \setminus \{0\}$.
	Define $H^p_{\psi,L}(X)$ to be the completion of the space
	\[
				\bH^p_{\psi,L}(X) :=\{f \in L^2(X) \,:\,  \scrA Q_{\psi,L}f \in L^p(X)\},
	\]
	 with respect to the norm
$
 			\norm{f}_{H^p_{\psi,L}(X)} := \norm{\scrA Q_{\psi,L}f}_{L^p(X)}. 
$
Then $H^p_L(X)=H^p_{\psi,L}(X)$, with equivalence of norms.
\end{Theorem}

Next, we recall the definition of the space $BMO_L(X)$.
One first defines a space $\calE_M(L)$ in such a way that for every $f \in \calE_M(L)$ there holds $(I-e^{r_B^{2m}L})^M f \in L^2_{\loc}(X)$, and therefore the expression in \eqref{defBMO-1} is well-defined.

\begin{Def}
Let $\eps>0$, $M \in \N$ and let $\phi \in \calR(L^M) \subseteq L^2(X)$ with $\phi = L^M \nu$ for some $\nu \in \calD(L^M)$. Introduce the norm
\[
 	\norm{\phi}_{\molMeps} := \sup_{j \geq 0} \left[2^{j\eps} V(2^jB_0)^{1/2} \sum_{k=0}^M \norm{L^k \nu}_{L^2(S_j(B_0))}\right],
\]
where $B_0$ is the unit ball centered at $0$ with radius 1 (cf. Section \ref{sect-homogSpace}), and set
\begin{align} \label{def-molMeps}
 	\molMeps:=\{\phi \in \calR(L^M) \,:\, \norm{\phi}_{\molMeps} < \infty\}.
\end{align}
One denotes by $(\molMeps)'$ the dual of $\molMeps$. 
For any $M \in \N$, let $\calE_M(L)$ be defined by
\[
 	\calE_M(L):= \bigcap_{\eps>0} (\molMdeps)'.
\]
\end{Def}

\begin{Remark}
Let $M \in \N$ and $\eps>0$. 
Then for every $f \in (\molMdeps)'$ and every $t>0$, one can define $(I-e^{-t^{2m}L})^M f$ and $(I-(I+t^{2m}L)^{-1})^Mf$ via duality as elements of $L^2_{\loc}(X)$.
\end{Remark}

\begin{Def}
 Let  $M \in \N$. An element $f \in \calE_M(L)$ is said to belong to $BMO_{L,M}(X)$ if
\begin{equation} \label{defBMO-1}
 	\norm{f}_{BMO_{L,M}(X)} := 
		\sup_{B \subseteq X} \left(\frac{1}{V(B)} \int_B \abs{(I-e^{-r_B^{2m}L})^M f(x)}^2 d\mu(x) \right)^{1/2} < \infty,
\end{equation}
where the supremum is taken over all balls $B$ in $X$.
\end{Def}

One can then show the following duality result. For a proof, we refer to \cite{DuongLi}, Theorem 3.23 and 3.24. 

 \begin{Theorem} \label{Duality}
Let $M>\frac{n}{4m}$. Then
$
 	(H^1_L(X))' = BMO_{L^{\ast},M}(X).
$
 \end{Theorem}

 In particular, the theorem yields that the definition of $BMO_{L,M}(X)$ is independent of the choice of $M>\frac{n}{4m}$. This leads to the following definition.
\begin{Def} \label{DefBMO-uniform}
The space $BMO_L(X)$ is defined by
$
		BMO_L(X):=BMO_{L,M}(X),
$
where $M \in \N$ with $M>\frac{n}{4m}$.
\end{Def}

The relation of elements of $BMO_L(X)$ and Carleson measures can be described as follows. 

  \begin{Prop} \label{LemmaCarlesonMeasure}
Let $M \in \N, \; M > \frac{n}{4m}$. Further, let $\alpha>0$, $\beta> \frac{n}{4m}$ and $\psi \in \Psi_{\beta,\alpha}(\Sigma_{\sigma}^0)\setminus\{0\}$. Then the operator
 \[
  	f \mapsto \psi(t^{2m}L)f
 \]
 maps $BMO_{L}(X) \to T^{\infty}(X)$, i.e. for every $f \in BMO_{L}(X)$ is
 \begin{equation} \label{CarlesonMeasure}
  	\nu_{\psi,f} := \abs{\psi(t^{2m}L)f(y)}^2 \frac{d\mu(y)\, dt}{t}
 \end{equation}
 a Carleson measure and there exists a constant $C_{\psi}>0$ such that for all $f \in BMO_{L}(X)$
 \[
  	\norm{\nu_{\psi,f}}_{\calC} \leq C_{\psi} \norm{f}^2_{BMO_{L}(X)}.
 \]
Conversely, if $f \in \calE_M(L)$ satisfies the controlled growth bound 
 \begin{equation} \label{growthBound}
  	\int_{X} \frac{\abs{(I-(I+L)^{-1})^M f(x)}^2}{(1+d(x,0))^{\eps_1}V(0,1+d(x,0))} \, d\mu(x) < \infty
 \end{equation}
for some $\eps_1>0$, and if $\nu_{\psi,f}$ defined in \eqref{CarlesonMeasure} is a Carleson measure, then $f \in BMO_L(X)$ and
 \[
  	\norm{f}_{BMO_L(X)}^2 \leq C \norm{\nu_{\psi,f}}_{\calC}.
 \]
 \end{Prop}

For a special choice of $\psi$, namely $\psi(z)=z^Me^{-z}$, the result is due to \cite{HofmannMayboroda}, Theorem 9.1. In the generality as stated above, the first part of the result is due to \cite{HofmannMayborodaMcIntosh}, Proposition 4.13. The second part is new and can be shown by combining the proof of \cite{HofmannMayboroda}, Theorem 9.1 with Lemma 3.17 of \cite{Frey}.\\


The spaces $H^p_L(X)$ form a complex interpolation scale.
For a proof, we refer to \cite{HofmannMayborodaMcIntosh}, Lemma 4.24, where the authors reduce the problem to complex interpolation of tent spaces.

\begin{Prop} \label{Interpolation}
 Let $L$ be an operator satisfying \eqref{H1} and \eqref{H2}. Let $1 \leq p_0<p_1<\infty$ and $0<\theta<1$.
 Then 
 \begin{align*}
 		[H^{p_0}_L(X),H^{p_1}_L(X)]_{\theta} & = H^p_L(X) \qquad \text{where} \;\, 1/p=(1-\theta)/p_0 + \theta/p_1, \\
 				[H^{p_0}_L(X),BMO_L(X)]_{\theta} &= H^p_L(X) \qquad \text{where} \;\, 1/p=(1-\theta)/p_0.
 \end{align*}
\end{Prop}

The next result is a slight generalization of \cite{HofmannMayboroda}, Theorem 3.2 and complements \cite{BlunckKunstmann}, Theorem 1.1.  

\begin{Prop} \label{H1-bddness}
Let $M \in \N$, $M>\frac{n}{4m}$. Assume that $T$ is a linear or a non-negative sublinear operator defined on $L^2(X)$ such that
$
 	T: L^2(X) \to L^2(X)
$
is bounded and $T$ satisfies the following weak off-diagonal estimates:\\
There exists some $\gamma > \frac{n}{2m}$ and a constant $C>0$ such that for every $t>0$, arbitrary balls $B_1,B_2 \in X$ with radius $r=t^{1/2m}$ and every $f \in L^2(X)$ supported in $B_1$
\begin{align} \label{H1-eq1}
 	\norm{T(I-e^{-tL})^M(f)}_{L^2(B_2)} &\leq C_T \left(1+\frac{\dist(B_1,B_2)^{2m}}{t}\right)^{-\gamma} \norm{f}_{L^2(B_1)},  \\ \label{H1-eq2}
	\norm{T(tLe^{-tL})^M(f)}_{L^2(B_2)} &\leq C_T \left(1+\frac{\dist(B_1,B_2)^{2m}}{t}\right)^{-\gamma} \norm{f}_{L^2(B_1)}.
\end{align}
Then
$
 	T: H^1_L(X) \to L^1(X)
$
is bounded and there exists some $C>0$, independent of $C_T$, such that for all  $f \in H^1_L(X)$ 
\begin{align*}
		\norm{Tf}_{L^1(X)} \leq C C_T \norm{f}_{H^1_L(X)}.
\end{align*}
\end{Prop}

\begin{Remark}
 If \eqref{H1-eq1} and \eqref{H1-eq2} are satisfied for arbitrary open sets $E,F \subseteq X$, one only requires a decay of order $\gamma>\frac{n}{4m}$.
\end{Remark}

A sufficient condition and a detailed proof for the equivalence of $H^p_L(X)$ and $L^p(X)$ is given in \cite{Uhl}, Theorem 4.19.
We refer the reader to a comparison with assumption \eqref{H3}.

\begin{Prop} \label{Hp-equiv}
Let $L$ satisfy \eqref{H1}, \eqref{H2}.	If for some $p_0 \in [1,2)$, there exist constants $C,c>0$ such that for all $x,y \in X$ and all $t>0$ 
\begin{align*}
		& \norm{\Eins_{B(x,t^{1/2m})} e^{-tL} \Eins_{B(y,t^{1/2m})}}_{L^{p_0}(X) \to L^{p'_0}(X)} \\
				& \qquad \qquad \qquad \qquad \leq C V(x,t^{1/2m})^{-(\frac{1}{p_0} - \frac{1}{p'_0})} \exp\left(-\left(\frac{d(x,y)^{2m}}{ct}\right)^{\frac{1}{2m-1}}\right),
\end{align*}
then 
\begin{align*}
		H^p_L(X) = L^p(X), \qquad p_0<p < p'_0.
\end{align*}	
\end{Prop}

For further relationships between $H^p_L(X)$ and $L^p(X)$ in the case of second order elliptic operators in divergence form, we refer to \cite{HofmannMayborodaMcIntosh}, Proposition 9.1.

 \section{A $T(1)$-Theorem for non-integral operators}
 \label{sect-T1Theorem}

In this section, we state and prove our first main result, Theorem \ref{T1-Theorem}. We first fix our assumptions on the operator $T$, clarify how, under these assumptions, the expressions $T(1)$ and $T^{\ast}(1)$ can be defined and discuss Poincar\'{e} estimates on metric spaces, that will be used in the proof of the main result.\\
After statement and proof of Theorem \ref{T1-Theorem}, which gives sufficient conditions for $L^2$ boundedness of $T$, we explain how to extend $T$ to Hardy spaces $H^p_L(X)$ for $p \neq 2$ and give necessary conditions for the boundedness of $T$ on $L^2(X)$.
We then add a second version of the $T(1)$-Theorem under weaker off-diagonal estimates and apply this version to prove the boundedness of a paraproduct operator on $L^2(X)$. We finally present a possible approach towards a $T(b)$-Theorem.\\
Throughout this section, we will always assume $L$ to be an operator satisfying  \eqref{H1}, \eqref{H2} and \eqref{H3}. \\


Let us fix our main assumptions on the operator $T$.

 \paragraph{Assumption}
 Let $\sigma \in (\omega, \frac{\pi}{2})$, and let $\alpha \geq 1$ and $\beta>\frac{n}{4m} +[\frac{n}{4m}]+1$.\\
 Let $T: \calD(L) \cap \calR(L) \to L^2_{\loc}(X)$ be a linear operator with $T^{\ast}: \calD(L^{\ast}) \cap \calR(L^{\ast}) \to L^2_{\loc}(X)$, which satisfies the following off-diagonal estimates:

\begin{enumerate}[label= \textbf{(OD\arabic*)\boldmath$_{\gamma}$},align=left,leftmargin=*,ref=(OD\arabic*)$_{\gamma}$]
\item \label{OD1} 
There exists a function $\psi_1 \in \Psi_{\beta,\alpha}(\Sigma_{\sigma}^0)\setminus\{0\}$, some $\gamma > 0$ and a constant $C>0$  such that $\psi_1(L)$ is injective and for every $t>0$, arbitrary balls $B_1,B_2 \in X$ with radius $r=t^{1/2m}$ and every $f \in L^2(X)$ supported in $B_1$ 
\begin{equation} \label{off-diag-est1}
 	\norm{T\psi_1(tL)f}_{L^2(B_2)} \leq C \left(1+\frac{\dist(B_1,B_2)^{2m}}{t} \right)^{-\gamma} \norm{f}_{L^2(B_1)}.
\end{equation}
\item \label{OD2}
There exists a function $\psi_2 \in \Psi_{\beta,\alpha}(\Sigma_{\sigma}^0)\setminus\{0\}$, some $\gamma > 0$ and a constant $C>0$  such that $\psi_2(L^{\ast})$ is injective and for every $t>0$, arbitrary balls $B_1,B_2 \in X$ with radius $r=t^{1/2m}$ and every $f \in L^2(X)$ supported in $B_1$ 
\begin{equation} \label{off-diag-est2}
 	\norm{T^{\ast}\psi_2(tL^{\ast})f}_{L^2(B_2)} \leq C \left(1+\frac{\dist(B_1,B_2)^{2m}}{t} \right)^{-\gamma} \norm{f}_{L^2(B_1)}.
\end{equation}
\end{enumerate}

Whenever we say that a linear operator $T$ satisfies \ref{OD1} or \ref{OD2}, we mean that $T$  satisfies \eqref{off-diag-est1} or \eqref{off-diag-est2}, respectively, for $\sigma,\alpha,\beta,\psi_1,\psi_2,C$ as specified above. The parameter $\gamma>0$ will be specified in each situation separately.\\

The assumptions that $\psi_1(L)$ and $\psi_2(L^{\ast})$ are injective are only used to define $T(1)$ and $T^{\ast}(1)$ in an appropriate way. If in applications it is clear how to do this, then the assumptions can be omitted. In that case, one can also relax the assumptions $\alpha \geq 1$ and $\beta>\frac{n}{4m} +[\frac{n}{4m}]+1$ to $\alpha>0$ and $\beta>\frac{n}{4m}$.

 \subsection{Definition of $T(1)$ and $T^{\ast}(1)$}
 \label{sect-DefT1}

Before we can state our $T(1)$-Theorem, we first have to clarify how to understand the expressions $T(1)$ and $T^{\ast}(1)$ for a linear operator $T: \calD(L) \cap \calR(L) \to L^2_{\loc}(X)$ with $T^{\ast}: \calD(L^{\ast}) \cap \calR(L^{\ast}) \to L^2_{\loc}(X)$, that satisfies \ref{OD2} and \ref{OD1}, respectively, for some $\gamma>\frac{n}{2m}$.
We confine ourselves to the definition of $T^{\ast}(1)$. How to define $T(1)$ will then be obvious.
We emphasize that this problem has not been addressed in \cite{Bernicot}.\\

The first observation is a simple consequence of Remark \ref{Def-on-Linfty}. If $T: \calD(L) \cap \calR(L) \to L^2_{\loc}(X)$ is a linear operator that satisfies \ref{OD1} for some $\gamma>\frac{n}{2m}$, then $\psi_1(tL^{\ast})T^{\ast}(1)$ can be defined via duality as an element of $L^2_{\loc}(X)$, i.e.
$
		\skp{\psi_1(tL^{\ast})T^{\ast}(1),\varphi} := \skp{1,T \psi_1(tL) \varphi}
$
for all $\varphi \in L^2(X)$ that are supported in some ball $B \subseteq X$.\\
We then define a space $Y^{\psi,\eps}(L)$ that will replace the space $\molMeps$ defined in \eqref{def-molMeps}.

\begin{Def} \label{Def-molSpace}
  Let $\eps>0$ and let $\alpha\geq 1$ and $\beta>\frac{n}{4m} +[\frac{n}{4m}]+1$. Let $\psi \in \Psi_{\beta,\alpha}(\Sigma_{\sigma}^0)\setminus\{0\}$ such that $\psi(L)$ is injective. 
  We define 
\begin{align*}
		 Y^{\psi,\eps}(L):=\{\phi =\psi(L)b \,:\, b \in L^2(X), \; \lim_{j \to \infty} 2^{j\eps} V(2^jB_0)^{1/2} \norm{b}_{L^2(S_j(B_0))} =0\},
\end{align*}
with the norm given by
\begin{align*}
		\norm{\phi}_{Y^{\psi,\eps}(L)}:=\sup_{j \geq 0} \left[2^{j\eps} V(2^jB_0)^{1/2} \norm{b}_{L^2(S_j(B_0))}\right].
\end{align*}
In addition, we define
\begin{align*}
		 Y^{\psi,\eps}_c(L):=\{\phi =\psi(L)b \in  Y^{\psi,\eps}(L) \,:\, \supp b \subseteq B \ \text{for some ball} \ B \subseteq X\}
\end{align*}
and 
\begin{align*}
		\calE_{\psi}(L):= \bigcap_{\eps>0} (Y^{\psi,\eps}(L^{\ast}))'.
\end{align*}
\end{Def}

\begin{Remark}
For every $\psi \in  \Psi(\Sigma_{\sigma}^0)$ as specified in Definition \ref{Def-molSpace} and every $\eps>0$, the space $Y^{\psi,\eps}(L)$ is a Banach space and $Y^{\psi,\eps}_c(L)$ is a dense subset of $Y^{\psi,\eps}(L)$. Moreover, the following inclusion holds:\\
Let $M \in \N$ with $M>\frac{n}{4m}$. Let further $\alpha\geq 1$, $\beta>\frac{n}{4m}+M$ and $\psi \in  \Psi_{\beta,\alpha}(\Sigma_{\sigma}^0)\setminus \{0\}$ such that $\psi(L)$ is injective. Then for every $\eps>0$ with $\frac{\eps}{2m}\leq \beta-(M+\frac{n}{4m})$
\begin{align*}
		Y^{\psi,\eps}(L) \subseteq \molMeps. 
\end{align*}
The result is also true for functions $\psi \in \Psi(\Sigma_{\sigma}^0)$ with $z \mapsto z^{-M} \psi(z) \in H^{\infty}(\Sigma_{\sigma}^0)$ such that the family of operators $\{(tL)^{-M}\psi(tL)\}_{t>0}$ satisfies Davies-Gaffney estimates. In this case, the inclusion is valid for all $\eps>0$. \\

For the proof, let $\phi \in  	Y^{\psi,\eps}(L)$, where $\phi=\psi(L)b$ for some $b \in L^2(X)$. Since $\beta>\frac{n}{4m}+M$, there obviously holds $\phi \in \calR(L^M)$. In addition, we have to show that $\norm{\phi}_{\molMeps} < \infty$ (see \eqref{newMol-eq3} below for a definition of the norm). First, observe that 
\begin{align} \label{newMol-bound}
		\norm{b}_{L^2(X)} \leq \sum_{j=0}^{\infty} \norm{b}_{L^2(S_j(B_0))} \leq C_{\eps} V(B_0)^{-1/2} \norm{\phi}_{Y^{\psi,\eps}(L)}
\end{align}
for some constant $C_{\eps}>0$ only depending on $\eps>0$.
Moreover, observe that for every $k=0,1,\ldots,M$, the function $z \mapsto z^{-(M-k)} \psi(z)$ is an element of  $\Psi_{\beta-M,\alpha}(\Sigma_{\sigma}^0)$. Thus, Proposition \ref{H-inf-offdiag} yields that the operator family $\{(tL)^{-(M-k)}\psi(tL)\}_{t>0}$ satisfies off-diagonal estimates of order $\beta-M$. 
Let us now write $b =\Eins_{R_j}b + \Eins_{(R_j)^c}b$ with
 \begin{align*}
  	R_j 	& = 2^{j+2} B_0, \qquad \text{if} \ j=0,1,2, \\
 	R_j		& = 2^{j+2} B_0 \setminus 2^{j-2}B_0, \qquad \text{if} \ j=3,4,\ldots.
 \end{align*}
For every $k=0,1,\ldots,M$ and all $j \in \N_0$, we obtain due to the boundedness of $L^{-(M-k)} \psi(L)$ on $L^2(X)$
\begin{align} \label{newMol-eq1}
		\norm{L^{-(M-k)} \psi(L) \Eins_{R_j} b}_{L^2(S_j(B_0))}
				\lesssim \norm{b}_{L^2(R_j)} 
				\lesssim 2^{-j\eps} V(2^jB_0)^{-1/2} \norm{\phi}_{Y^{\psi,\eps}(L)},
\end{align}
where in the last step $R_j$ is splitted into four annuli. On the other hand, the off-diagonal estimates for $\{(tL)^{-(M-k)}\psi(tL)\}_{t>0}$, \eqref{newMol-bound} and the doubling property \eqref{doublingProperty2} yield
\begin{align} \label{newMol-eq2} \nonumber
		\norm{L^{-(M-k)} \psi(L) \Eins_{(R_j)^c} b}_{L^2(S_j(B_0))}
			&\lesssim \left(1+\dist(S_j(B_0), (R_j)^c)^{2m}\right)^{-(\beta-M)} \norm{b}_{L^2(R_j)^c} \\ \nonumber
			&\lesssim 2^{-2m(\beta-M)j} \norm{b}_{L^2(X)} \\
			& \lesssim 2^{-2m(\beta-M)j} 2^{jn/2} V(2^jB_0)^{-1/2} \norm{\phi}_{Y^{\psi,\eps}(L)}.
\end{align}
We therefore obtain from \eqref{newMol-eq1}, \eqref{newMol-eq2} and the assumption $\frac{\eps}{2m}\leq \beta-(M+\frac{n}{4m})$
\begin{align} \label{newMol-eq3} \nonumber
		\norm{\phi}_{\molMeps}
				&=\sup_{j \geq 0} \left[2^{j\eps} V(2^jB_0)^{1/2} \sum_{k=0}^M \norm{L^{-(M-k)}\psi(L)b}_{L^2(S_j(B_0))}\right] \\ \nonumber
				& \lesssim   \sup_{j \geq 0} 2^{j\eps} \left(2^{-j\eps} + 2^{-2m(\beta-(M+\frac{n}{4m}))j} \right)\norm{\phi}_{Y^{\psi,\eps}(L)}
				 \lesssim  \norm{\phi}_{Y^{\psi,\eps}(L)}.
\end{align}
Since Davies-Gaffney estimates imply off-diagonal estimates of any order, the second case is then obvious.
\end{Remark}

Let us now define $T^{\ast}(1)$ as an element of $\calE_{\psi_1}(L^{\ast})$ in the following way.

\begin{Lemma} \label{Def-T1-Lemma}
Let $T: \calD(L) \cap \calR(L) \to L^2_{\loc}(X)$ be a linear operator that satisfies \emph{\ref{OD1}} for some $\gamma>\frac{n}{2m}$. 
Then $T^{\ast}(1)$ can be defined as an element of $\calE_{\psi_1}(L^{\ast})$ by 
\begin{align} \label{Def-T1-Eq}
		\skp{T^{\ast}(1),\phi}:=\lim_{R \to \infty} \skp{\psi_1(L^{\ast}) T^{\ast}(\Eins_{B(0,R)}),b}
\end{align}
for every $\phi \in Y^{\psi_1,\eps}_c(L)$ with $\phi=\psi_1(L)b$ and every $\eps>0$.
\end{Lemma}

In the same way, one can then also define $T(1)$ as an element of $\calE_{\psi_2}(L)$ under the assumption that $T^{\ast}: \calD(L^{\ast}) \cap \calR(L^{\ast}) \to L^2_{\loc}(X)$ satisfies \ref{OD2} for some $\gamma >\frac{n}{2m}$.

\begin{Proof}
Let $\gamma>\frac{n}{2m}$. The assumption \ref{OD1} yields, according to Remark \ref{Def-on-Linfty}, the following estimate: There exists some constant $C_T>0$ such that for every ball $B$ in $X$ with radius $r_B>0$ and every $f \in L^{\infty}(X)$
\begin{equation} \label{DefT1-BMO-eq1}
 			V(B)^{-1/2} \norm{\psi_1(r_B^{2m}L^{\ast}) T^{\ast}(f)}_{L^2(B)} \leq C_T \norm{f}_{L^{\infty}(X)}.
\end{equation}
As mentioned before, the left hand side of \eqref{DefT1-BMO-eq1} is well-defined via duality. 
With the help of the above estimate, we can now define $T^{\ast}(1)$ as an element of $\calE_{\psi_1}(L^{\ast})=\bigcap_{\eps>0} (Y^{\psi_1,\eps}(L))'$ as follows.\\
Let $\eps>0$. We define for every $R>0$ a linear functional $\ell_R$ on $Y^{\psi_1,\eps}_c(L)$ by setting
\begin{align*}
		\ell_R(\phi):= \skp{\psi_1(L^{\ast}) T^{\ast}(\Eins_{B(0,R)}),b}
\end{align*}
for every $\phi =\psi_1(L) b \in Y^{\psi_1,\eps}_c(L)$. Observe that $\ell_R(\phi)$ is well-defined, since $b$ is supported in some ball of $X$ and $\psi_1(L^{\ast}) T^{\ast}(\Eins_{B(0,R)})$ is via duality defined  as an element of $L^2_{\loc}(X)$. 
Using the definition of $\norm{\,.\,}_{Y^{\psi_1,\eps}(L)}$ and \eqref{DefT1-BMO-eq1}, we obtain
\begin{align} \label{DefT1-BMO-eq2} \nonumber
		\abs{\ell_R(\phi)} 
				&= \abs{\skp{\psi_1(L^{\ast}) T^{\ast}(\Eins_{B(0,R)}),b}} 
				 \leq \sum_{j=0}^{\infty} \norm{\psi_1(L^{\ast}) T^{\ast}(\Eins_{B(0,R)})}_{L^2(S_j(B_0))} \norm{b}_{L^2(S_j(B_0))} \\ \nonumber
				& \lesssim  \sum_{j=0}^{\infty}  2^{-j\eps} V(2^jB_0)^{-1/2}  \norm{\psi_1(L^{\ast}) T^{\ast}(\Eins_{B(0,R)})}_{L^2(S_j(B_0))} \norm{\phi}_{Y^{\psi_1,\eps}(L)} 
				 \lesssim C_T \norm{\phi}_{Y^{\psi_1,\eps}(L)},
\end{align}
where the implicit constants are independent of $R>0$. Thus, 
$
		\sup_{R>0}	\abs{\ell_R(\phi)} \lesssim  C_T \norm{\phi}_{Y^{\psi_1,\eps}(L)}.
$
Following the estimates in \eqref{DefT1-BMO-eq1} and \eqref{DefT1-BMO-eq2}, we can moreover show that $(\ell_R(\phi))_R$ is a Cauchy sequence for every $\phi \in Y^{\psi_1,\eps}_c(L)$. Hence, $\lim_{R \to \infty} \ell_R(\phi)$ exists.
Since $Y^{\psi_1,\eps}_c(L)$ is dense in $Y^{\psi_1,\eps}(L)$ and $\eps>0$ was arbitrary, we can now define $T^{\ast}(1) \in \calE_{\psi_1}(L^{\ast})$ by \eqref{Def-T1-Eq}.
\end{Proof}

 \subsection{Poincar\'{e} inequalities}
 \label{sect-Poincare}
 
For the proof of our $T(1)$-Theorem, we require some kind of Poincar\'{e} inequality. 
We follow the approach of Haj{\l}asz and Koskela in \cite{HajlaszKoskela1} and \cite{HajlaszKoskela2}, who give generalizations of Poincar\'{e} inequalities and Sobolev spaces on metric spaces.
Our basic tool will be the following definition, which is taken from Chapter 2 of \cite{HajlaszKoskela2}.

\begin{Def} \label{Def-Poincare}
Assume that $u \in L^1_{\loc}(X)$ and a measurable function $g \geq 0$ satisfy the inequality
\begin{align} \label{Poincare-pair}
		\frac{1}{V(B)} \int_B \abs{u(x) - \skp{u}_B} \,d\mu(x) 
				\leq C_P \, r_B \left(\frac{1}{V(\lambda B)} \int_{\lambda B} g(x)^p \,d\mu(x)\right)^{1/p},
\end{align}
on each ball $B$ in $X$, where $r_B$ is the radius of $B$ and $p>0$, $\lambda \geq 1$, $C_P>0$ are fixed constants.
We then say that the pair $(u,g)$ satisfies a $p$-Poincar\'{e} inquality.
\end{Def}

\begin{Remark}
	If $u \in \Lip(\R^n)$, $g=\abs{\nabla u}$ and $p \geq 1$, then \eqref{Poincare-pair} is a corollary of the classical Poincar\'{e} inequality
\begin{align} \label{classical-Poincare}
	\left(\frac{1}{V(B)} \int_{B} \abs{u(x)-\skp{u}_B}^p \,dx\right)^{1/p}
				\leq C(n,p) \, r_B \left(\frac{1}{V(B)} \int_B \abs{\nabla u(x)}^p \,dx\right)^{1/p}.
\end{align}
\end{Remark}

It is therefore natural to consider a pair $(u,g)$ that satisfies a $p$-Poincar\'{e} inequality as a Sobolev function and its ``gradient''. We refer to \cite{HajlaszKoskela2} for a survey on the topic and examples of pairs $(u,g)$ on certain metric spaces that satisfy a $p$-Poincar\'{e} inequality.\\

Let us formulate the required assumption.

\paragraph{Assumption}
Let $L$ satisfy \eqref{H1}, \eqref{H2} and \eqref{H3}.
\begin{enumerate}[label= \textbf{(P)},ref=P]
\item \label{P}
Assume that for every $f \in L^2(X)$ there exists a measurable function $g: X \times (0,\infty) \to \C$ such that
\renewcommand{\labelenumii}{(\roman{enumii})}
\begin{enumerate}
\item for all $t>0$ there holds $g_t:=g(\,.\,,t) \geq 0$, and the pair $(e^{-t^{2m}L}f,g_t)$ satisfies a $p$-Poincar\'{e} inequality of the form \eqref{Poincare-pair} for some $p<2$ and with constants  $\lambda \geq 1, \; C_P>0$ independent of $t$ and $f$;
\item for all $t>0$ there holds $g_t \in L^2(X)$, and there exists a constant $C>0$ independent of $f$ with 
\begin{align*}
		\int_0^{\infty} t^2 \norm{g_t}_{L^2(X)}^2 \,\frac{dt}{t} \leq C \norm{f}_{L^2(X)}^2.
\end{align*}
\end{enumerate}
\end{enumerate}
\begin{enumerate}[label= \textbf{(P$^{\ast}$)},ref=P$^{\ast}$]
\item \label{Pd} Assume that \eqref{P} holds with $L$ replaced by $L^{\ast}$.
\end{enumerate}

\begin{Remark}
If $X$ is the Euclidean space $\R^n$, then the Poincar\'{e} inequality is automatically satisfied for the pairs $(e^{-t^{2m}L}f, \abbs{\nabla e^{-t^{2m}L}f})$  and $(e^{-t^{2m}L^{\ast}}f, \abbs{\nabla e^{-t^{2m}L^{\ast}}f})$, see e.g. \cite{GilbargTrudinger}, (7.45). 
In this case, (ii) is just the assumption that the vertical Littlewood-Paley-Stein square function is bounded on $L^2(\R^n)$. For elliptic second order operators in divergence form, this can easily be shown with the help of the ellipticity condition, see e.g. \cite{Auscher}, Section 6.1. 
In general, (ii) is fulfilled whenever the Riesz transforms $\nabla L^{-1/2m}$, $\nabla (L^{\ast})^{-1/2m}$
are bounded on $L^2(\R^n)$, since then
\begin{align*}
	\int_0^{\infty} \norm{t\nabla e^{-t^{2m}L}f}_{L^2(\R^n)}^2 \,\frac{dt}{t}
			= \int_0^{\infty} \norm{\nabla L^{-1/2m} (t^{2m}L)^{1/2m} e^{-t^{2m}L}f}_{L^2(\R^n)}^2 \,\frac{dt}{t}
			\lesssim \norm{f}_{L^2(\R^n)}^2
\end{align*}
and the analogous estimate for $L^{\ast}$ hold due to quadratic estimates, see \eqref{square-functions}.\\
Let us reformulate the assumptions \eqref{P} and \eqref{Pd} also for another case. Let $X$ be a complete Riemannian manifold, with the Riemannian measure $\mu$ on $X$ satisfying the doubling property \eqref{doublingProperty2}, and let $\nabla$ denote the Riemannian gradient. To obtain (i) of \eqref{P}, it is sufficient to assume that a $2$-Poincar\'{e} inequality of the form \eqref{classical-Poincare} holds (with the Lebesgue measure replaced by $\mu$). One can then again choose the pairs $(e^{-t^{2m}L}f, \abbs{\nabla e^{-t^{2m}L}f})$  and $(e^{-t^{2m}L^{\ast}}f, \abbs{\nabla e^{-t^{2m}L^{\ast}}f})$. 
This is due to a certain self-improving property of Poincar\'{e} inequalities on Riemannian manifolds, stating that the interval of all $p$ that satisfy a $p$-Poincar\'{e} inequality, is open. We refer to \cite{KeithZhong} for details.
Sufficient for (ii) is, as for the Euclidean space,  that the mappings $f \mapsto \abbs{\nabla L^{-1/2m}f}$ and $f \mapsto \abbs{\nabla (L^{\ast})^{-1/2m}f}$ are bounded on $L^2(X)$.
\end{Remark}

 The following theorem is a  simplified version of \cite{HajlaszKoskela2}, Theorem 3.2.

\begin{Theorem} \label{Poincare-Thm}
Assume that the pair $(u,g)$ satisfies a $p$-Poincar\'{e} inequality \eqref{Poincare-pair} for some $p>0$. Then there exists some constant $C>0$ such that 
\begin{align*}
		\abs{u(x)-u(y)} \leq C d(x,y) \left(\calM_p g(x) + \calM_p g(y)\right)
\end{align*}
for almost every $x,y \in X$. 
\end{Theorem}

%
 \subsection{Main theorem}
 
We are now ready to state our main theorem.

 \begin{Theorem}  \label{T1-Theorem}
Let $L$ be an operator satisfying the assumptions \eqref{H1}, \eqref{H2} and \eqref{H3}. Additionally, let the assumptions \eqref{P} and \eqref{Pd} be satisfied.
 Let $T: \calD(L) \cap \calR(L) \to L^2_{\loc}(X)$ be a linear operator with $T^{\ast}: \calD(L^{\ast}) \cap \calR(L^{\ast}) \to L^2_{\loc}(X)$, which satisfies the assumptions \emph{\ref{OD1}} and \emph{\ref{OD2}} for some $\gamma>\frac{n+D+2}{2m}$ and let $T(1) \in BMO_L(X), \, T^{\ast}(1) \in BMO_{L^{\ast}}(X)$. \\
 Then $T$ is bounded on $L^2(X)$, i.e. there exists a constant $C>0$ 
 such that for all $f \in L^2(X)$ 
 \[
  	\norm{Tf}_{L^2(X)} \leq C \norm{f}_{L^2(X)}.
 \]
 \end{Theorem}

Let us sketch the main idea of the proof. \\
First, we approximate $T$ by operators associated to $L$, namely, we write with the help of a Calder\'{o}n reproducing formula for $f,g \in L^2(X)$
\begin{align} \label{Cald-repr}
		 \skp{Tf,g}
 		& =  \int_0^{\infty} \int_0^{\infty} \skp{\psi_2(t^{2m}L) T \psi_1(s^{2m}L) \tilde{\psi}_1(s^{2m}L) f, \tilde{\psi}_2(t^{2m}L^{\ast})g} \, \frac{dt}{t} \frac{ds}{s}.
\end{align}
We then decompose the operator  $T$ for each $t>0$, at least formally, in the following way:
 \begin{align} \label{T-decomposition} \nonumber
  	T	&= T(I-e^{-t^{2m}L}) + Te^{-t^{2m}L}  \\ 
 		&= T(I-e^{-t^{2m}L}) + [Te^{-t^{2m}L} - T(1) \cdot A_te^{-t^{2m}L}] + T(1) \cdot A_t e^{-t^{2m}L}.
 \end{align}
 This can be understood as a splitting of the operator into the ``main term'' or ``principal part'' $Te^{-t^{2m}L}$ and the ``error term'' $T(I-e^{-t^{2m}L})$. The main term is then further decomposed  into the term in the squared brackets, which is handled via Poincar\'{e} inequalities and the term $T(1) \cdot A_t e^{-t^{2m}L}$, which can be estimated by application of the theory of paraproducts and use of the assumption $T(1) \in BMO_{L}(X)$.
The idea of such a decomposition is taken from \cite{AxelssonKeithMcIntosh} and \cite{HytoenenMcIntoshPortal}.
In the case $t<s$ in \eqref{Cald-repr}, the error term can easily dealt with via almost orthogonality arguments and quadratic estimates. For $t>s$, one argues via duality and uses the assumption $T^{\ast}(1) \in BMO_{L^{\ast}}(X)$.\\

The boundedness of the occuring paraproduct operator on $L^2(X)$ has been shown in \cite{Frey}, Theorem 4.2.

 \begin{Theorem} \label{paraproduct}
  Assume that $L$ satisfies \eqref{H1}, \eqref{H2} and \eqref{Lp-L2-estimate} of \eqref{H3}. Let $\alpha>0$, $\beta>\frac{n}{4m}$ and let $\psi \in \Psi_{\beta,\alpha}(\Sigma_{\sigma}^0)\setminus \{0\}$,  $\tilde{\psi} \in \Psi(\Sigma_{\sigma}^0)\setminus \{0\}$.
Then the operator $\Pi_{b,L}$, defined for every $f \in L^2(X)$ and every  $b \in BMO_L(X)$ by
 \begin{equation} \label{Paraproduct}
 	\Pi_{b,L}(f):= \int_0^{\infty} \tilde{\psi}(t^{2m}L) [\psi(t^{2m}L)b \cdot A_t(e^{-t^{2m}L}f)] \, \frac{dt}{t},
 \end{equation}
 where $A_t$ is the averaging operator defined in \eqref{avOp},
 is bounded on $L^2(X)$. I.e. there exists some constant $C>0$ such that for every $f \in L^2(X)$ and every $b \in BMO_L(X)$
 \[
  	\norm{\Pi_{b,L}(f)}_{L^2(X)} \leq C \norm{b}_{BMO_L(X)} \norm{f}_{L^2(X)}.
 \]
 Analogously, if $L$ satisfies \eqref{H1}, \eqref{H2} and \eqref{L2-Lq-estimate} of \eqref{H3}, then for every $b \in BMO_{L^{\ast}}(X)$ the operator $\Pi_{b,L^{\ast}}$ is bounded on $L^2(X)$.
 \end{Theorem}

For the treatment of the term in the squared brackets in \eqref{T-decomposition}, we use the following proposition. The idea is taken from \cite{AxelssonKeithMcIntosh}, Proposition 5.5.

\begin{Prop} \label{mainTerm}
Assume that \eqref{P} holds. Let $\{S_t\}_{t>0}$ be a family of linear operators on $L^2(X)$ that satisfies weak off-diagonal estimates of order $\gamma>\frac{n+D+2}{2m}$. Then there exists a constant $C>0$ such that for all $f \in L^2(X)$
\[
  	 \int_0^{\infty} \norm{S_{t^{2m}} e^{-t^{2m}L} f - S_{t^{2m}}(1) \cdot A_t e^{-t^{2m}L} f}_{L^2(X)}^2 \,\frac{dt}{t}
 		\leq C \norm{f}_{L^2(X)}^2.
\]
\end{Prop}

\begin{Proof}
Let $f \in L^2(X)$. The assumption (i) of \eqref{P} yields for every $t>0$ the existence of some function $g_t \in L^2(X)$ such that the pair $(e^{-t^{2m}L}f,g_t)$ satisfies a $p$-Poincar\'{e} inequality for some $p<2$. 
If we can show that there exists some $C>0$, independent of $t$ and $f$, such that
\begin{align} \label{mainTerm-eq1}
		\norm{S_{t^{2m}} e^{-t^{2m}L} f - S_{t^{2m}}(1) \cdot A_t e^{-t^{2m}L} f}_{L^2(X)}^2 
				\leq C t^2 \norm{g_t}_{L^2(X)}^2,
\end{align}
then the assertion of the lemma is a consequence of assumption (ii) of \eqref{P}. \\
Let $t>0$ be fixed and abbreviate $u:=e^{-t^{2m}L}f$. 
 To apply the weak off-diagonal estimates on $S_t$, we decompose $X$ with the help of Lemma \ref{ChristCubes} into ``cubes'' of ``sidelength'' approximately equal to $t$. That is, with the notation of Lemma \ref{ChristCubes}, let $k_0 \in \Z$ be defined by $C_1 \delta^{k_0} \leq t < C_1 \delta^{k_0-1}$ and write $X=\bigcup_{\alpha \in I_{k_0}} Q_{\alpha}^{k_0}$, where the equality holds modulo null sets of $\mu$. By Lemma \ref{ChristCubes} we further know that for every $\alpha \in I_{k_0}$ there exists some $z_{\alpha}^{k_0} \in X$ such that
\begin{align} \label{mainTerm-eq1a}
		B(z_{\alpha}^{k_0},c_1t) \subseteq Q_{\alpha}^{k_0} \subseteq B(z_{\alpha}^{k_0},t)
\end{align}
for some $c_1 \in (0,1)$ independent of $t$ and $\alpha$. 
Moreover, observe that the averaging operator $A_t$ is, by definition, constant on each ``cube'' $Q_{\alpha}^{k_0}$. 
We therefore get
\begin{align} \label{mainTerm-eq2} \nonumber
	  	&  \norm{S_{t^{2m}}u - S_{t^{2m}}(1) \cdot A_t u}_{L^2(X)}^2 
			 =\sum_{\alpha \in I_{k_0}} \norm{S_{t^{2m}}u - S_{t^{2m}}(1) \cdot A_tu}_{L^2(Q_{\alpha}^{k_0})}^2 \\ \nonumber
			& \qquad =\sum_{\alpha \in I_{k_0}} \norm{S_{t^{2m}}(u -\skp{u}_{Q_{\alpha}^{k_0}})}_{L^2(Q_{\alpha}^{k_0})}^2
			 \leq \sum_{\alpha \in I_{k_0}} \left(\sum_{\beta \in I_{k_0}} \norm{S_{t^{2m}}\Eins_{Q_{\beta}^{k_0}} (u -\skp{u}_{Q_{\alpha}^{k_0}})}_{L^2(Q_{\alpha}^{k_0})}\right)^2\\ 
			& \qquad \lesssim \sum_{\alpha \in I_{k_0}} \left(\sum_{\beta \in I_{k_0}} \left(1+\frac{\dist(B_{\alpha},B_{\beta})^{2m}}{t^{2m}}\right)^{-\gamma} \norm{u-\skp{u}_{Q_{\alpha}^{k_0}}}_{L^2(Q_{\beta}^{k_0})} \right)^2.
\end{align}
Observe that due to \eqref{CubeSummation-eq} and  $\gamma>\frac{n}{2m}$
\begin{align} \label{mainTerm-eq1b} 
 	\sup_{\alpha \in I_{k_0}} \sum_{\beta \in I_{k_0}} \left(1+\frac{\dist(B_{\alpha},B_{\beta})^{2m}}{t^{2m}}\right)^{-\gamma}
			\lesssim 1.
\end{align}
The Cauchy-Schwarz inequality then yields that the expression in \eqref{mainTerm-eq2} is bounded by
\begin{align} \label{mainTerm-eq3} \nonumber
		& \sum_{\alpha \in I_{k_0}} \left(\sum_{\beta \in I_{k_0}} \left(1+\frac{\dist(B_{\alpha},B_{\beta})^{2m}}{t^{2m}}\right)^{-\gamma} \right) 
		\left( \sum_{\beta \in I_{k_0}} \left(1+\frac{\dist(B_{\alpha},B_{\beta})^{2m}}{t^{2m}}\right)^{-\gamma} \norm{u-\skp{u}_{Q_{\alpha}^{k_0}}}_{L^2(Q_{\beta}^{k_0})}^2 \right) \\
		 & \qquad \lesssim \sum_{\alpha \in I_{k_0}}  \sum_{\beta \in I_{k_0}} \left(1+\frac{\dist(B_{\alpha},B_{\beta})^{2m}}{t^{2m}}\right)^{-\gamma} \norm{u-\skp{u}_{Q_{\alpha}^{k_0}}}_{L^2(Q_{\beta}^{k_0})}^2.
\end{align}
The term $\norm{u-\skp{u}_{Q_{\alpha}^{k_0}}}_{L^2(Q_{\beta}^{k_0})}^2$ is now handled via the assumed $p$-Poincar\'{e} inequality for the pair $(u,g_t)$. Due to the Cauchy-Schwarz inequality and Theorem \ref{Poincare-Thm} we get
\begin{align} \label{mainTerm-eq4} \nonumber
		& \norm{u-\skp{u}_{Q_{\alpha}^{k_0}}}_{L^2(Q_{\beta}^{k_0})}^2 
		= \int_{Q_{\beta}^{k_0}} \abs{u(x)-\skp{u}_{Q_{\alpha}^{k_0}}}^2 \,d\mu(x) \\ \nonumber
		& \qquad \leq \int_{Q_{\beta}^{k_0}} \left(\frac{1}{V(Q_{\alpha}^{k_0})} \int_{Q_{\alpha}^{k_0}} \abs{u(x)-u(y)} \,d\mu(y)\right)^2 \,d\mu(x) \\ \nonumber
		& \qquad \leq \frac{1}{V(Q_{\alpha}^{k_0})} \int_{Q_{\beta}^{k_0}} \int_{Q_{\alpha}^{k_0}}  \abs{u(x)-u(y)}^2 \,d\mu(y) \,d\mu(x) \\
		& \qquad \lesssim \frac{1}{V(Q_{\alpha}^{k_0})}  \int_{Q_{\beta}^{k_0}} \int_{Q_{\alpha}^{k_0}}
					d(x,y)^2 [\calM_pg_t(x) + \calM_pg_t(y)]^2 \,d\mu(y) \,d\mu(x).
\end{align}
Note that for $x \in Q_{\beta}^{k_0}$ and $y \in Q_{\alpha}^{k_0}$ there holds $d(x,y) \lesssim t\left(1 + \dist(B_{\alpha},B_{\beta})/t\right)$ due to \eqref{mainTerm-eq1a}. Moreover, the doubling property \eqref{DoublingProperty3} and \eqref{mainTerm-eq1a} yield that
\begin{align*}
		\frac{V(Q_{\beta}^{k_0})}{V(Q_{\alpha}^{k_0})} \lesssim \left(1+\frac{\dist(B_{\alpha},B_{\beta})}{t}\right)^D.
\end{align*}
Taking these considerations into account and plugging \eqref{mainTerm-eq4} into \eqref{mainTerm-eq3}, we end up with
\begin{align} \label{mainTerm-eq5} \nonumber
		& \norm{S_{t^{2m}}u - S_{t^{2m}}(1) \cdot A_t u}_{L^2(X)}^2 
		  \lesssim \sum_{\alpha \in I_{k_0}}  \sum_{\beta \in I_{k_0}}
					\left(1+\frac{\dist(B_{\alpha},B_{\beta})^{2m}}{t^{2m}}\right)^{-\gamma}
					\norm{u-\skp{u}_{Q_{\alpha}^{k_0}}}_{L^2(Q_{\beta}^{k_0})}^2 \\ \nonumber
		& \qquad \lesssim  t^2 \sum_{\alpha \in I_{k_0}}  \sum_{\beta \in I_{k_0}}
					\left(1+\frac{\dist(B_{\alpha},B_{\beta})}{t}\right)^{-2m\gamma+2} \\ \nonumber
		& \qquad \qquad	 \qquad \times \left[\int_{Q_{\beta}^{k_0}} [\calM_pg_t(x)]^2 \,d\mu(x)
						 		+ \frac{V(Q_{\beta}^{k_0})}{V(Q_{\alpha}^{k_0})} \int_{Q_{\alpha}^{k_0}} [\calM_p g_t(y)]^2 \,d\mu(y)\right] \\ \nonumber
		& \qquad \lesssim 	t^2 \sum_{\alpha \in I_{k_0}}  \sum_{\beta \in I_{k_0}}
					\left(1+\frac{\dist(B_{\alpha},B_{\beta})}{t}\right)^{-2m\gamma+D+2} \\ \nonumber
			& \qquad \qquad	 \qquad \times \left[\int_{Q_{\beta}^{k_0}} [\calM_pg_t(x)]^2 \,d\mu(x)
						 		+ \int_{Q_{\alpha}^{k_0}} [\calM_p g_t(y)]^2 \,d\mu(y)\right] \\	\nonumber	
		& \qquad \lesssim t^2 \left[ \sum_{\beta \in I_{k_0}} \int_{Q_{\beta}^{k_0}} [\calM_pg_t(x)]^2 \,d\mu(x)
							+ \sum_{\alpha \in I_{k_0}}  \int_{Q_{\alpha}^{k_0}} [\calM_p g_t(y)]^2 \,d\mu(y) \right] \\
		& \qquad \lesssim t^2 \norm{\calM_p g_t}_{L^2(X)} 
					\lesssim t^2 \norm{g_t}_{L^2(X)}^2, 
\end{align}
where we used \eqref{mainTerm-eq1b} with the assumption $\gamma>\frac{n+D+2}{2m}$, the disjointness of the ``cubes'' and the boundedness of $\calM_p$ on $L^2(X)$ for $p<2$ in the last three inequalities.
This shows \eqref{mainTerm-eq1}, which in turn finishes the proof by assumption (ii) of \eqref{P}.
\end{Proof}

The next lemma gives a certain kind of almost orthogonality for operators constructed via $H^{\infty}$-functional calculus.
The first part is a corollary of Lemma \ref{off-diag-st}, the second part is a simple modification of the arguments given there (observe that the roles of $s$ and $t$ are interchanged).

\begin{Lemma} \label{st-estimates}
Let $\alpha,\beta>0$ and let $\psi \in \Psi_{\beta,\alpha}(\Sigma_{\sigma}^0)$. There exists a constant $C>0$ such that for every $s,t>0$
\[
 	\norm{(I-e^{-tL}) \psi(sL) f}_{L^2(X) \to L^2(X)} \leq C \left(\frac{t}{s}\right)^{\min(\alpha,1)} 
\]
and
\[
 	\norm{e^{-tL} \psi(sL)f}_{L^2(X) \to L^2(X)} \leq C \left(\frac{s}{t}\right)^{\beta}.
\]
\end{Lemma}

 Now, we are ready to prove our main theorem.
 
 \begin{Proof}[of Theorem \ref{T1-Theorem}]
 Let $f,g \in L^2(X)$. Let $\alpha\geq 1, \ \beta> \frac{n}{4m} +[\frac{n}{4m}]+1$ and let $\psi_1, \psi_2 \in \Psi_{\beta,\alpha}(\Sigma_{\sigma}^0)\setminus \{0\}$ as given in the assumption.
Corresponding to the functions $\psi_1,\psi_2$, we choose functions $\tilde{\psi}_1, \tilde{\psi}_2 \in \Psi(\Sigma_{\sigma}^0)$ such that $\int_0^{\infty} \psi_1(t) \tilde{\psi}_1(t) \, \frac{dt}{t} =1$ and $\int_0^{\infty} \psi_2(t) \tilde{\psi}_2(t) \, \frac{dt}{t} =1$ and decompose both $f$ and $g$ with the help of the Calder\'{o}n reproducing formula. That is, we write 
 \begin{align*}
  	 \skp{Tf,g}
 		& =  \int_0^{\infty} \int_0^{\infty} \skp{\psi_2(t^{2m}L) T \psi_1(s^{2m}L) \tilde{\psi}_1(s^{2m}L) f, \tilde{\psi}_2(t^{2m}L^{\ast})g} \, \frac{dt}{t} \frac{ds}{s} 
 \end{align*}
and show that the right-hand side is bounded by a constant times $\norm{f}_{L^2(X)} \norm{g}_{L^2(X)}$. In this way, $T$ extends to a bounded operator on $L^2(X)$.\\
For the proof, we split the inner integral into two parts, one over $\{t \in (0,\infty):  0 < t < s\}$, called $J_1$, and the other one over $\{t \in (0,\infty): s \leq t < \infty\}$, called $J_2$. We observe that for the second part $J_2$,  Fubini's theorem yields
 \begin{align*}
  	J_2 &= \int_0^{\infty} \int_s^{\infty} \skp{\psi_2(t^{2m}L) T \psi_1(s^{2m}L) \tilde{\psi}_1(s^{2m}L) f, \tilde{\psi}_2(t^{2m}L^{\ast}) g} \, \frac{dt}{t} \frac{ds}{s} \\
 	&  = \int_0^{\infty} \int_0^s \skp{\tilde{\psi}_1(t^{2m}L) f, \psi_1(t^{2m}L^{\ast}) T^{\ast} \psi_2(s^{2m}L^{\ast}) \tilde{\psi}_2(s^{2m}L^{\ast}) g} \, \frac{dt}{t} \frac{ds}{s}.
 \end{align*}
The last line equals $J_1$ with $T$ replaced by $T^{\ast}$, $L$ by $L^{\ast}$ and the roles of $\psi_1,\tilde{\psi}_1$ and $\psi_2,\tilde{\psi}_2$ interchanged. Note that all our assumptions are symmetric with respect to $T,T^{\ast}$ and $L,L^{\ast}$. Moreover, instead of the weak off-diagonal estimates for $\{T\psi_1(t^{2m}L)\}_t$, assumed in \ref{OD1}, we can take into account the analogous estimates for $\{T^{\ast}\psi_2(t^{2m}L^{\ast})\}_t$, assumed in \ref{OD2}. Thus, it will be sufficient to study only $J_1$. Once we have proven this part, the estimate for $J_2$ will follow by duality.
In the following estimate for $J_1$, we will always assume $0<t<s$.\\
As described in \eqref{T-decomposition}, we decompose  $T$ into the two parts $Te^{-t^{2m}L}$ and $T(I-e^{-t^{2m}L})$ for every $t>0$, which leads to
 \begin{align} \label{T1-eq1b} \nonumber
  J_1 	&=  \int_0^{\infty} \int_0^s \skp{\psi_2(t^{2m}L) T \psi_1(s^{2m}L) \tilde{\psi}_1(s^{2m}L) f, \tilde{\psi}_2(t^{2m}L^{\ast}) g} \, \frac{dt}{t} \frac{ds}{s} \\ \nonumber
 	&  = 	\int_0^{\infty} \int_0^s \skp{\psi_2(t^{2m}L) T e^{-t^{2m}L} \psi_1(s^{2m}L) \tilde{\psi}_1(s^{2m}L) f, \tilde{\psi}_2(t^{2m}L^{\ast}) g} \, \frac{dt}{t} \frac{ds}{s} \\ \nonumber
 	&  \qquad + \int_0^{\infty} \int_0^s \skp{\psi_2(t^{2m}L) T(I-e^{-t^{2m}L}) \psi_1(s^{2m}L) \tilde{\psi}_1(s^{2m}L) f, \tilde{\psi}_2(t^{2m}L^{\ast}) g} \, \frac{dt}{t} \frac{ds}{s} \\
 	&  =: J_M + J_E.
 \end{align}
Let us first turn to the estimation of the error term $J_E$, the main term $J_M$ will be treated below. 
Due to assumption \ref{OD2}, with $\gamma>\frac{n}{2m}$, and Lemma \ref{uniformL2-bound}, $\{\psi_2(t^{2m}L) T\}_{t>0}$ is uniformly bounded on $L^2(X)$. Combining this with Lemma \ref{st-estimates} yields
 \begin{align} \label{T1-eq2a}
  	\norm{\psi_2(t^{2m}L)T(I-e^{-t^{2m}L}) \psi_1(s^{2m}L)}_{L^2(X) \to L^2(X)}
 			& \lesssim \left(\frac{t^{2m}}{s^{2m}}\right)^{\min(\alpha,1)}.
 \end{align}
We can therefore estimate $J_E$ with the help of the Cauchy-Schwarz inequality by
 \begin{align} \label{T1-eq3} \nonumber
  	\abs{J_E}	& \leq \int_0^{\infty} \int_0^s \abs{\skp{\psi_2(t^{2m}L) T(I-e^{-t^{2m}L}) \psi_1(s^{2m}L) \tilde{\psi}_1(s^{2m}L) f, \tilde{\psi}_2(t^{2m}L^{\ast}) g}} \, \frac{dt}{t} \frac{ds}{s} \\ \nonumber
 				& \lesssim \int_0^{\infty} \int_0^s \left(\frac{t^{2m}}{s^{2m}}\right)^{\min(\alpha,1)} \norm{\tilde{\psi}_1(s^{2m}L)f}_{L^2(X)} \norm{\tilde{\psi}_2(t^{2m}L^{\ast})g}_{L^2(X)} \, \frac{dt}{t} \frac{ds}{s} \\ \nonumber
 				& \leq \left(\int_0^{\infty} \int_0^s \left(\frac{t^{2m}}{s^{2m}}\right)^{\min(\alpha,1)} \norm{\tilde{\psi}_1(s^{2m}L)f}_{L^2(X)}^2 \, \frac{dt}{t} \frac{ds}{s} \right)^{1/2} \\
 					&	\qquad \qquad \times\left(\int_0^{\infty} \int_t^{\infty}  \left(\frac{t^{2m}}{s^{2m}}\right)^{\min(\alpha,1)} \norm{\tilde{\psi}_2(t^{2m}L^{\ast})g}_{L^2(X)}^2 \, \frac{ds}{s} \frac{dt}{t} \right)^{1/2},
 	\end{align}
where we also used Fubini's theorem in the last step.
By substitution of $u=\frac{t}{s}$, one easily observes that $\dis \int_0^{s} \left(\frac{t}{s}\right)^{\delta} \, \frac{dt}{t} = \delta^{-1}$ for every $\delta>0$. Since the operator family $\{\tilde{\psi}_1(sL)\}_{s>0}$ satisfies quadratic estimates, see \eqref{square-functions}, the first factor in the last line of \eqref{T1-eq3} can therefore be bounded by 
\begin{align} \label{T1-eq3a}
		\left(\int_0^{\infty} \int_0^s \left(\frac{t^{2m}}{s^{2m}}\right)^{\min(\alpha,1)} \norm{\tilde{\psi}_1(s^{2m}L)f}_{L^2(X)}^2 \, \frac{dt}{t} \frac{ds}{s} \right)^{1/2} 
				\lesssim \norm{f}_{L^2(X)}.
\end{align}
Changing the roles of $s$ and $t$ and using that $\dis \int_t^{\infty}\left(\frac{t}{s}\right)^{\delta} \, \frac{ds}{s}=\delta^{-1}$ for every $\delta>0$, we get the analogous estimate for the second factor in \eqref{T1-eq3} and in summary
\begin{align} \label{T1-eq4}
		\abs{J_E} \lesssim \norm{f}_{L^2(X)} \norm{g}_{L^2(X)}.
\end{align}
  
 To estimate the main term $J_M$, we use the extended decomposition in \eqref{T-decomposition} of $Te^{-t^{2m}L}$ into the two parts $[Te^{-t^{2m}L} - T(1) \cdot A_te^{-t^{2m}L}]$ and $T(1) \cdot A_t e^{-t^{2m}L}$. At the same time, we withdraw the decomposition of the function $f$ by the Calder\'{o}n reproducing formula at scale $s$. To do so, we do not consider $J_M$ itself, but the same expression, now called $J_M^0$, with both paths of integration over the whole interval $(0,\infty)$. This leads to
 \begin{align} \label{T1-eq4a} \nonumber
   J_M^0 &=  \int_0^{\infty} \int_0^{\infty} \skp{\psi_2(t^{2m}L) Te^{-t^{2m}L}  \psi_1(s^{2m}L) \tilde{\psi}_1(s^{2m}L)f, \tilde{\psi}_2(t^{2m}L^{\ast})g} \, \frac{ds}{s} \frac{dt}{t}  \\ \nonumber
 	&  =  \int_0^{\infty} \skp{\psi_2(t^{2m}L) Te^{-t^{2m}L} f, \tilde{\psi}_2(t^{2m}L^{\ast})g} \, \frac{dt}{t} \\ \nonumber
 	&    = \int_0^{\infty} \skp{\psi_2(t^{2m}L) Te^{-t^{2m}L}f - \psi_2(t^{2m}L) T(1) \cdot A_t e^{-t^{2m}L}f, \tilde{\psi}_2(t^{2m}L^{\ast})g} \, \frac{dt}{t} \\ \nonumber
 	& \qquad +  \int_0^{\infty} \skp{\psi_2(t^{2m}L)T(1) \cdot A_t e^{-t^{2m}L}f,\tilde{\psi}_2(t^{2m}L^{\ast})g} \,\frac{dt}{t}\\
 	& =: J_M^1 + J_M^2.
 \end{align}
The term $J_M^2$ is exactly the paraproduct defined in \eqref{Paraproduct}, i.e. $J_M^2 = \skp{\Pi_{T(1),L}(f),g}$,
with the functions $\psi,\tilde{\psi}$ replaced by $\psi_2, \tilde{\psi}_2$. 
Recall that we assumed in \ref{OD2} that $\psi_2 \in \Psi_{\beta,\alpha}(\Sigma_{\sigma}^0)$ for some $\alpha>0$ and $\beta>\frac{n}{4m}$, and moreover assumed $T(1)$ to be an element of $BMO_L(X)$. Thus, $\Pi_{T(1),L}$ is bounded on $L^2(X)$ due to Theorem \ref{paraproduct} and we obtain the estimate
\begin{align} \label{T1-eq5}
		\abs{J_M^2} \lesssim \norm{T(1)}_{BMO_L(X)} \norm{f}_{L^2(X)} \norm{g}_{L^2(X)}.
\end{align}
It remains to find a bound for $J_M^1$. But the major part of this estimate was already done in Proposition \ref{mainTerm} by application of the assumed Poincar\'{e} inequalities \eqref{P}. Thus, if we set $S_{t^{2m}}:=\psi_2(t^{2m}L)T$ and take into account the assumption \ref{OD2} with $\gamma>\frac{n+D+2}{2m}$, then Proposition \ref{mainTerm}, in combination with the Cauchy-Schwarz inequality, yields
 \begin{align} \label{T1-eq6} \nonumber
 	\abs{J_M^1} &=  \abs{ \int_0^{\infty} \skp{\psi_2(t^{2m}L) Te^{-t^{2m}L} f - \psi_2(t^{2m}L)T(1) \cdot A_t e^{-t^{2m}L}f,\tilde{\psi}_2(t^{2m}L^{\ast})g} \, \frac{dt}{t}} \\ \nonumber
 				&  \leq 		\left(\int_0^{\infty} \norm{\psi_2(t^{2m}L) T e^{-t^{2m}L} f - \psi_2(t^{2m}L) T(1) \cdot A_t e^{-t^{2m}L} f}_{L^2(X)}^2 \, \frac{dt}{t}\right)^{1/2} \\ \nonumber
 				& \qquad 	\qquad \times \left(\int_0^{\infty} \norm{\tilde{\psi}_2(t^{2m}L^{\ast})g}_{L^2(X)}^2 \, \frac{dt}{t} \right)^{1/2} \\
 				& \lesssim \norm{f}_{L^2(X)} \norm{g}_{L^2(X)},
 \end{align}
 where we also used quadratic estimates for the operator family $\{\tilde{\psi}_2(tL^{\ast})\}_{t>0}$ in the last step.\\ 
Let us finally observe what we did wrong by considering $J_M^0$ instead of $J_M$. The combination of \eqref{T1-eq5} and \eqref{T1-eq6} provides us with the estimate
\begin{align} \label{T1-eq7}
			\abs{J_M^0} \lesssim \left(\norm{T(1)}_{BMO_L(X)}  + 1\right)\norm{f}_{L^2(X)} \norm{g}_{L^2(X)}.
\end{align}
On the other hand, we have $J_M = J_M^0-J_R$, where the remainder term $J_R$ is defined by
\begin{align*}
		J_R := \int_0^{\infty} \int_s^{\infty} \skp{\psi_2(t^{2m}L) Te^{-t^{2m}L}  \psi_1(s^{2m}L) \tilde{\psi}_1(s^{2m}L)f, \tilde{\psi}_2(t^{2m}L^{\ast})g} \, \frac{dt}{t} \frac{ds}{s}.
\end{align*}
This term can be handled in analogy to the treatment of $J_E$, replacing the estimate \eqref{T1-eq2a} by
\begin{align*}
  	\norm{\psi_2(t^{2m}L) Te^{-t^{2m}L}\psi_1(s^{2m}L)}_{L^2(X) \to L^2(X)}
 		\lesssim \left(\frac{s^{2m}}{t^{2m}}\right)^{\beta},
\end{align*}
which again holds uniformly for all $s,t>0$ according to Lemma \ref{st-estimates} and the uniform boundedness of $\{\psi_2(t^{2m}L) T\}_{t>0}$ in $L^2(X)$.
Together with the Cauchy-Schwarz inequality, the above yields
 \begin{align} \label{T1-eq8} \nonumber
  	\abs{J_R} & \leq \int_0^{\infty} \int_s^{\infty} \abs{\skp{\psi_2(t^{2m}L) Te^{-t^{2m}L}  \psi_1(s^{2m}L) \tilde{\psi}_1(s^{2m}L)f, \tilde{\psi}_2(t^{2m}L^{\ast}) g}} \, \frac{dt}{t} \frac{ds}{s} \\ 
 		&  \lesssim \int_0^{\infty} \int_s^{\infty}  \left(\frac{s^{2m}}{t^{2m}}\right)^{\beta} \norm{\tilde{\psi}_1(s^{2m}L)f}_{L^2(X)} \norm{\tilde{\psi}_2(t^{2m}L^{\ast})g}_{L^2(X)} \, \frac{dt}{t} \frac{ds}{s}.  \end{align}
If we now handle the last line of \eqref{T1-eq8} with the same argument as used in \eqref{T1-eq3} and \eqref{T1-eq3a}, we end up with
\begin{align} \label{T1-eq9}
	\abs{J_R} \lesssim \norm{f}_{L^2(X)} \norm{g}_{L^2(X)}.
\end{align}
By combining \eqref{T1-eq4}, \eqref{T1-eq7} and \eqref{T1-eq9}, and repeating the same procedure for $J_2$ and recalling the splitting $\skp{Tf,g} = J_1 + J_2 = J_E+  J_M^0 - J_R + J_2$,
we finally obtain
\begin{align*}
		\abs{\skp{Tf,g}} \lesssim \left(\norm{T(1)}_{BMO_L(X)} + \norm{T^{\ast}(1)}_{BMO_{L^{\ast}}(X)} + 1\right)\norm{f}_{L^2(X)} \norm{g}_{L^2(X)}.
\end{align*}
This proves the theorem.
 \end{Proof}


\subsection{Extension to $H^p_L(X)$ for $p \neq 2$ and necessary conditions}

If $T$ satisfys off-diagonal estimates \ref{OD1} and \ref{OD2} and is bounded on $L^2(X)$, then the extension to Hardy spaces $H^p_L(X)$ for $p\neq 2$ is almost immediate.
Such a property is similar to the behaviour of Calder\'{o}n-Zygmund operators, in respect of the fact that every Calder\'{o}n-Zygmund operator, that is bounded on $L^2(X)$, is automatically also bounded on $L^p(X)$ for all $p \in (1,\infty)$.\\
We start with the following self-improving property of off-diagonal estimates. The proof is postponed to Section 5.

\begin{Lemma} \label{self-improving-est}
Let $\psi \in \Psi(\Sigma_{\sigma}^0) \setminus \{0\}$ and let $T$ be a linear operator on $L^2(X)$ such that $\{T\psi(tL)\}_{t>0}$ satisfies weak off-diagonal estimates of order $\gamma >\frac{n}{2m}$ on $L^2(X)$. Let $\delta>\gamma$ and let $\varphi \in H^{\infty}(\Sigma_{\sigma}^0)$ with $\abs{\varphi(z)} \lesssim \abs{z}^{\delta}$ for $\abs{z} \leq 1$. Moreover, assume that $\{T\varphi(tL)\}_{t>0}$ is uniformly bounded on $L^2(X)$.
Then $\{T\varphi(tL)\}_{t>0}$ satisfies weak off-diagonal estimates of order $\gamma$ on $L^2(X)$.
\end{Lemma}

We then obtain the following.

\begin{Cor} \label{Lp-bddness}
Let $L$ be an operator satisfying the assumptions \eqref{H1} and \eqref{H2}. 
Let $T:L^2(X) \to L^2(X)$ be a bounded linear operator that satisfies \emph{\ref{OD1}} for some $\gamma>\frac{n}{2m}$. Then $T$ extends to a bounded operator
\begin{align*}
			T: H^p_L(X) \to L^p(X), \quad \qquad 1 \leq p \leq 2,
\end{align*}
and $T^{\ast}$ extends to a bounded operator
\begin{align*}
			T^{\ast} &: L^{p}(X) \to H^p_{L^{\ast}}(X) , \quad \qquad 2 \leq p < \infty, \\
			T^{\ast} &: L^{\infty}(X) \to BMO_{L^{\ast}}(X).
\end{align*}
\end{Cor}

One can obviously obtain the corresponding results for  $T^{\ast},L^{\ast}$ in place of $T,L$, if one uses \ref{OD2} instead of \ref{OD1}.
To obtain boundedness results for $T: L^p(X) \to L^p(X)$, we refer the reader to combine Corollary \ref{Lp-bddness} with Proposition \ref{Hp-equiv}. Following the proof of \cite{HofmannMayborodaMcIntosh}, Proposition 5.6, one can moreover show that $T$ extends to a bounded operator $T: H^1_L(X) \to H^1(X)$, whenever $T^{\ast}(1)=0$.

\begin{Proof}
In order to show that $T$ extends to a bounded operator $T: H^1_L(X) \to L^1(X)$, one combines Proposition \ref{H1-bddness} with Lemma \ref{self-improving-est}. Let $M \in \N$ with $M>\gamma$. Observe that the operator families  $\{(I-e^{-tL})^M\}_{t>0}$, $\{(tLe^{-tL})^{M}\}_{t>0}$ and therefore also  $\{T(I-e^{-tL})^M\}_{t>0}$, $\{T(tLe^{-tL})^{M}\}_{t>0}$ are uniformly bounded on $L^2(X)$. Thus, Lemma \ref{self-improving-est} yields that  $\{T(I-e^{-tL})^M\}_{t>0}$ and $\{T(tLe^{-tL})^{M}\}_{t>0}$ satisfy weak off-diagonal estimates of order $\gamma$. \\
One then uses the interpolation scales for the spaces $L^p(X)$ and $H^p_L(X)$, see Proposition \ref{Interpolation}, and obtains the boundedness of $T: H^p_L(X) \to L^p(X)$ for $1 \leq p \leq 2$.
Since Theorem \ref{Duality} yields that $(H^1_L(X))'=BMO_{L^{\ast}}(X)$ and the space $H^p_{L^{\ast}}(X)$ was defined as the dual space of $H^{p'}_L(X)$ for $2 < p < \infty$ and $\frac{1}{p}+\frac{1}{p'}=1$ (see Definition \ref{Def-Hp}), one finally gets the remaining assertions of the corollary via duality.
\end{Proof}

\begin{Remark}
(i) 
The above results also contain the following necessary conditions for a non-integral operator to be bounded on $L^2(X)$:\\
Let $T: L^2(X) \to L^2(X)$ be a bounded linear operator.
If $T$ satisfies assumption \emph{\ref{OD2}} with $\gamma>\frac{n}{2m}$, then $T(1) \in BMO_L(X)$.
Analogously, if $T$ satisfies assumption \emph{\ref{OD1}} with $\gamma>\frac{n}{2m}$, then $T^{\ast}(1) \in BMO_{L^{\ast}}(X)$.\\
(ii)  Using the relation between elements of $BMO_L(X)$ and Carleson measures, as described in Proposition \ref{LemmaCarlesonMeasure}, one can also formulate the assumption "$T(1) \in BMO_L(X)$" in terms of Carleson measures.
\end{Remark}

\begin{Remark}
Let $L$ be a second order elliptic operator in divergence form.
For the Riesz transform $T=\nabla L^{-1/2}$ one can show that \ref{OD1} is satisfied, cf. e.g. \cite{KatoSquare}, \cite{HofmannMartell}, \cite{BlunckKunstmann2}. However, if $n \geq 3$, \ref{OD2} cannot be satisfied in general:\\
Denote by $(p_{-}(L),p_{+}(L))$ the interior of $L^p$ boundedness of $\{e^{-tL}\}_{t>0}$. Then, as shown in \cite{HofmannMayborodaMcIntosh}, Proposition 9.1, $H^p_L(X)=L^p(X)$ for all $p \in (p_{-}(L),p_{+}(L))$.
 If $\nabla L^{-1/2}$ would satisfy \ref{OD2}, then, according to Corollary \ref{Lp-bddness}, $\nabla L^{-1/2}$ would be bounded on $L^p(X)$ for $p \in [2,p_{+}(L))$. However, in \cite{Auscher} it was proved that $p_{+}(L)>\frac{2n}{n-2}$ and in \cite{AuscherTchamitchian}, Theorem 4.7 (due to Kenig) it was shown that for every $p>2$ there exists some second order elliptic operator in divergence form $L$ with $\nabla L^{-1/2}$ not bounded on $L^p(X)$. 
\end{Remark}

%
 \subsection{A second version under weaker off-diagonal estimates}
\label{sect-weakT1Theorem}

We give in this section a second version of Theorem \ref{T1-Theorem} under weaker estimates than \ref{OD1} and \ref{OD2}. However, we assume in addition that the conservation properties $e^{-tL}(1)=1$ and $e^{-tL^{\ast}}(1)=1$ in $L^2_{\loc}(X)$ are valid and that $T$ is a weakly continuous operator mapping from $L^2(X)$ to $L^2(X)$. The last assumption is stronger than the one in Theorem \ref{T1-Theorem}, but one thinks of an application to some kind of ``truncations'' $T_{\eps}$ of $T$ with uniform $L^2$ bound. See e.g.  Theorem \ref{thm-new-para} below for an example.  The basic idea of the construction is taken from \cite{Bernicot}. The proof, however, is completely different from \cite{Bernicot}, as we cannot use pointwise kernel bounds. The result gives a positive answer to a question raised at the end of \cite{Bernicot}.

\begin{Theorem} \label{weak-T1-Theorem}
Let $L$ be an operator satisfying the assumptions \eqref{H1}, \eqref{H2} and \eqref{H3}. Additionally, let the assumptions \eqref{P} and \eqref{Pd} be satisfied.
 Let $\alpha>0,\beta> \frac{n}{4m}$ and $\psi,\tilde{\psi} \in \Psi_{\beta,\alpha}(\Sigma_{\sigma}^0)$ with $\int_0^{\infty} \psi(t) \tilde{\psi}(t) \,\frac{dt}{t}=1$ and define $\phi \in H^{\infty}(\Sigma_{\sigma}^0)$ via
\[
 	\phi(z):= \int_{\gamma_z} \psi(\zeta) \tilde{\psi}(\zeta) \,\frac{d\zeta}{\zeta}, \qquad z \in \Sigma_{\sigma}^0,
\]
where $\gamma_z(t):=te^{i\arg z}, \; t \in (\abs{z},\infty)$.
Assume that the operator family $\{\phi(tL)\}_{t>0}$ satisfies off-diagonal estimates of order $\gamma>\frac{n+D+2}{2m}$ and moreover, assume that 
\begin{align} \label{conservation-law}
			\phi(tL)(1)=	\phi(tL^{\ast})(1)=1 \qquad \text{in} \ L^2_{\loc}(X)
\end{align}
 for every $t>0$.\\
Let $T: L^2(X) \to L^2(X)$ be a linear, weakly continuous operator such that
$\{\psi(tL) T \phi(tL)\}_{t>0}$ and $\{\psi(tL^{\ast})T^{\ast}\phi(tL^{\ast})\}_{t>0}$ satisfy weak off-diagonal estimates of order $\gamma>\frac{n+D+2}{2m}$ and let $T(1) \in BMO_L(X)$ and $T^{\ast}(1) \in BMO_{L^{\ast}}(X)$.\\
Then $T: L^2(X) \to L^2(X)$ is bounded with a constant independent of the weak continuity parameters of $T$. 
\end{Theorem}

\begin{Remark}
Note that one can get off-diagonal estimates for  $\{\phi(tL)\}_{t>0}$ in the following way. By splitting $\phi(z)=(\phi(z)-e^{-z})+e^{-z}$ for $z \in \Sigma_{\sigma}^0$, one can on the one hand take into account Davies-Gaffney estimates for the semigroup $\{e^{-tL}\}_{t>0}$. On the other hand, it is clear by definition that $\phi(z)-e^{-z} \to 0$ for $\abs{z} \to 0$ and for $\abs{z} \to \infty$. 
Proposition \ref{H-inf-offdiag} then yields the existence of off-diagonal estimates for $\{\phi(tL)-e^{-tL}\}_{t>0}$.\\
With a similar reasoning, one can show that the assumption \eqref{conservation-law} is a consequence of the property $e^{-tL}(1)=e^{-tL^{\ast}}(1)=1$ in $L^2_{\loc}(X)$. This is due to the fact that the latter implies $\psi(tL)(1)=\psi(tL^{\ast})(1)=0$ in $L^2_{\loc}(X)$ for every $\psi \in \Psi_{\beta,\alpha}(\Sigma_{\sigma}^0)$ with $\beta>\frac{n}{4m}$ and $\alpha>0$, see Lemma \ref{psi-remark}.
\end{Remark}

The proof of Theorem \ref{weak-T1-Theorem} is almost equal to the one of Theorem \ref{T1-Theorem}. The only difference is the replacement of the Calder\'{o}n reproducing formula by the representation formula \eqref{T1-weak-eq1}, which is a generalization of a construction in \cite{Bernicot}.

\begin{Proof}
We first observe that by definition of $\phi$ there holds $\lim_{t \to 0} \phi(t)=1$ and $\lim_{t \to \infty} \phi(t)=0$. Since $T$ is weakly continuous, we thus get by functional calculus
\begin{align*}
 	Tf	 = \lim_{t \to 0} \phi^2(tL) T \phi^2(tL)f, \qquad
	0 	 = \lim_{t \to \infty}  \phi^2(tL) T \phi^2(tL)f,
\end{align*}
where the limit is interpreted in the weak sense in $L^2(X)$. 
Again by functional calculus, we obtain from the above as a special form of a Calder\'{o}n reproducing formula that $\skp{Tf,g}$ can be represented as
\begin{align} \label{T1-weak-eq1}
 	\skp{Tf,g} = \skp{\int_0^{\infty} \left(\left[t \frac{d}{dt}\phi^2(tL)\right]T \phi^2(tL)(f)
				+ \phi^2(tL)T\left[t\frac{d}{dt}\phi^2(tL)\right]f\right) \,\frac{dt}{t},g}.
\end{align}
Once having handled the first summand in \eqref{T1-weak-eq1}, in the following called $J$, the second one will work in the same way simply by duality. So let us have a more detailed look at the first part.\\
By definition of $\phi$ there holds $z\phi'(z)=\psi(z)\tilde{\psi}(z)$ for $z \in \Sigma_{\sigma}^0$.  This yields due to functional calculus,
\begin{align} \label{T1-weak-eq2} \nonumber 
 	 \int_0^{\infty} \left[t\frac{d}{dt}\phi^2(tL)\right] T\phi^2(tL)(f) \,\frac{dt}{t} 
		 & = 2 \int_0^{\infty} (tL)\phi'(tL)\phi(tL)T\phi^2(tL)(f) \,\frac{dt}{t} \\
		& = 2 \int_0^{\infty} \psi(tL)\psi_1(tL) T\phi^2(tL)(f) \,\frac{dt}{t},
\end{align}
where we set $\psi_1(z):=\tilde{\psi}(z)\phi(z)$. 
We further decompose $f$ with the help of another Calder\'{o}n reproducing formula as
\begin{align} \label{T1-weak-eq3}
 	f = \int_0^{\infty} \psi(sL) \tilde{\psi}(sL) f \,\frac{ds}{s},
\end{align}
taking into account the assumption $\int_0^{\infty} \psi(t) \tilde{\psi}(t) \,\frac{dt}{t}=1$.
The combination of the two equations \eqref{T1-weak-eq2} and \eqref{T1-weak-eq3} then leads to
\begin{align} \label{T1-weak-eq4}
		J &= 2 \int_0^{\infty} \int_0^{\infty} \skp{\psi(t^{2m}L) T \phi^2(t^{2m}L) \psi(s^{2m}L) \tilde{\psi}(s^{2m}L) f, \psi_1(t^{2m}L^{\ast})g} \,\frac{dt}{t}\frac{ds}{s}.
\end{align}
Similar to the proof of Theorem \ref{T1-Theorem}, we split the inner integral into two parts, one over the interval $\{t \in (0,\infty):  0 < t < s\}$, called $J_1$, and the other one over $\{t \in (0,\infty): s \leq t < \infty\}$, called $J_2$.
In contrast to the proof of Theorem \ref{T1-Theorem}, for lack of symmetry in \eqref{T1-weak-eq4} we cannot handle $J_2$ simply by duality, but it can be dealt with similar to the remainder term $J_R$ in Theorem \ref{T1-Theorem}. \\ 
Thus, let us first turn to $J_2$. 
Lemma \ref{uniformL2-bound} and the assumed weak off-diagonal estimates yield that  $\{\psi(tL)T\phi(tL)\}_{t>0}$ is uniformly bounded on $L^2(X)$.
Moreover, observe that by assumption $\frac{\psi(\zeta)\tilde{\psi}(\zeta)}{\zeta} = \calO(\abs{\zeta}^{-2\alpha-1})$ for $\abs{\zeta} \to \infty$ and consequently, $\phi(z)=\calO(\abs{z}^{-2\alpha})$ for $\abs{z} \to \infty$. Replacing $e^{-z}$ by $\phi(z)$ in Lemma \ref{st-estimates}, it is therefore easy to check that there exists some $\delta>0$ such that
\begin{align*}
		\norm{\phi(t^{2m}L)\psi(s^{2m}L)h}_{L^2(X) \to L^2(X)} \lesssim \left(\frac{s^{2m}}{t^{2m}}\right)^{\delta},
\end{align*}
uniformly in $s,t>0$. With exactly the same arguments as in \eqref{T1-eq8}, we end up with
\begin{align*}
 	\abs{J_2} &\leq \int_0^{\infty} \int_s^{\infty} \abs{\skp{\psi(t^{2m}L)T\phi^2(t^{2m}L) \psi(s^{2m}L) \tilde{\psi}(s^{2m}L)f,\psi_1(t^{2m}L^{\ast})g}} \,\frac{dt}{t}\frac{ds}{s}
 			 \lesssim \norm{f}_{L^2(X)} \norm{g}_{L^2(X)}.
\end{align*}
To handle $J_1$, we apply for every $t>0$ the splitting
\begin{align*}
		T \phi^2(t^{2m}L) = T \phi^2(t^{2m}L)  e^{-t^{2m}L} + T \phi^2(t^{2m}L) (I-e^{-t^{2m}L}),
\end{align*}
representing the splitting of $J_1$ into the main term $J_M$ and the error term $J_E$ just as in \eqref{T1-eq1b}, i.e.
\begin{align*}
 J_1	& = \int_0^{\infty} \int_0^s \skp{\psi(t^{2m}L) T \phi^2(t^{2m}L) e^{-t^{2m}L} \psi(s^{2m}L) \tilde{\psi}(s^{2m}L) f, \psi_1(t^{2m}L^{\ast})g} \,\frac{dt}{t}\frac{ds}{s} \\
	& \qquad + \int_0^{\infty} \int_0^s \skp{\psi(t^{2m}L) T \phi^2(t^{2m}L) (I-e^{-t^{2m}L}) \psi(s^{2m}L) \tilde{\psi}(s^{2m}L) f, \psi_1(t^{2m}L^{\ast})g} \,\frac{dt}{t}\frac{ds}{s} \\
	& =: J_M + J_E.
\end{align*}
 The treatment of $J_E$ works analogously to \eqref{T1-eq3}, using the weak off-diagonal estimates for $\{\psi(tL)T\phi(tL)\}_{t>0}$ instead of assumption \ref{OD2} and the uniform boundedness of $\{\phi(tL)\}_{t>0}$ in $L^2(X)$.
To estimate the main term $J_M$, we also aim to apply a paraproduct estimate and therefore write $J_M = J_M^0 + J_R$ with a remainder $J_R$ that can be handled with the same arguments as in \eqref{T1-eq8}, and
\begin{align*}
		J_M^0 &= \int_0^{\infty} \int_0^{\infty} \skp{\psi(t^{2m}L) T \phi^2(t^{2m}L) e^{-t^{2m}L} \psi(s^{2m}L) \tilde{\psi}(s^{2m}L) f, \psi_1(t^{2m}L^{\ast})g} \,\frac{dt}{t}\frac{ds}{s} \\
			 	&    = \int_0^{\infty} \skp{\psi(t^{2m}L) T  \phi^2(t^{2m}L) e^{-t^{2m}L}f - \psi(t^{2m}L) T \phi^2(t^{2m}L)(1) \cdot A_t e^{-t^{2m}L}f, \psi_1(t^{2m}L^{\ast})g} \, \frac{dt}{t} \\ \nonumber
 	& \qquad +  \int_0^{\infty} \skp{\psi(t^{2m}L)T \phi^2(t^{2m}L)(1) \cdot A_t e^{-t^{2m}L}f,\psi_1(t^{2m}L^{\ast})g} \,\frac{dt}{t}\\
 	& =: J_M^1 + J_M^2,		
\end{align*}
in analogy to \eqref{T1-eq4a}.
Observe that the operator family $\{\psi(tL)T\phi^2(tL)\}_{t>0}$ satisfies weak off-diagonal estimates of order $\gamma>\frac{n+D+2}{2m}$ due to the assumptions and Proposition \ref{weak-composition}.
By taking assumption \eqref{P} into account, we can thus apply Proposition \ref{mainTerm} with $S_{t}:=\psi(tL)T\phi^2(tL)$, which yields the desired estimate for $J_M^1$ just as in \eqref{T1-eq6}.
We finally note that assumption \eqref{conservation-law} yields
\begin{align*}
		J_M^2 = \int_0^{\infty} \skp{\psi(t^{2m}L)T (1) \cdot A_t e^{-t^{2m}L}f,\psi_1(t^{2m}L^{\ast})g} \,\frac{dt}{t} 
					= \skp{\Pi_{T(1),L}f,g},
\end{align*}
and $J_M^2$ can therefore be treated by Theorem \ref{paraproduct} and the assumption $T(1) \in BMO_L(X)$.
\end{Proof}

 \subsection{Application to paraproducts}

We present an application of Theorem \ref{weak-T1-Theorem} to a special type of paraproduct operator (cf. also \cite{Bernicot}, Section 4).
We do this under more restrictive assumptions on $L$.
Let again $L$ be an operator satisfying \eqref{H1}, \eqref{H2} and \eqref{H3}. Additionally, let us assume the following. 
\begin{enumerate}[label= \textbf{(H\arabic*)},ref=H\theenumi]
\setcounter{enumi}{3}
\item \label{H4}  The operator 	$e^{-tL} : L^{\infty}(X) \to L^{\infty}(X)$ is bounded uniformly in $t>0$.
\item \label{H5}  For every $t>0$ there holds $e^{-tL}(1)=1$ in $L^{\infty}(X)$ and $e^{-tL^{\ast}}(1)=1$ in $L^2_{\loc}(X)$.
\end{enumerate}

The assumption \eqref{H5} in particular implies that  $\psi(tL^{\ast})(1)=0$ in $L^2_{\loc}(X)$ for every $t>0$ and every $\psi \in \Psi_{\beta,\alpha}(\Sigma_{\sigma}^0)$, where $\beta>\frac{n}{4m}$ and $\alpha>0$, see Lemma \ref{psi-remark}.

\begin{Def}
	Let $\alpha_1,\beta_1,\alpha_2,\beta_2>0$. Assume that $\psi_1 \in \Psi_{\beta_1,\alpha_1}(\Sigma_{\sigma}^0)\setminus\{0\}$ and $\psi_2 \in \Psi_{\beta_2,\alpha_2}(\Sigma_{\sigma}^0)\setminus\{0\}$ and abbreviate $\tilde{\psi}:=\psi_1 \cdot \psi_2$. For every $f \in L^{\infty}(X)$ and every $g \in L^2(X)$ we define the paraproduct
\begin{align} \label{Def-Para-new}
		\tilde{\Pi}_f(g):= \int_0^{\infty} \tilde{\psi}(t^{2m}L) [e^{-t^{2m}L}g \cdot e^{-t^{2m}L}f ] \,\frac{dt}{t}.
\end{align}
\end{Def}

If one would replace $e^{-t^{2m}L}g$ in \eqref{Def-Para-new} by $\psi(t^{2m}L)g$ for some $\psi \in \Psi(\Sigma_{\sigma}^0)$, then the boundedness of  $\tilde{\Pi}_f$ on $L^2(X)$ would be an immediate consequence of quadratic estimates due to bounded $H^{\infty}$-functional calculus for $L$. 
In our case, we obtain the following.

\begin{Theorem} \label{thm-new-para}
Let $L$ satisfy \eqref{H1}-\eqref{H5} and let the assumptions \eqref{P} and \eqref{Pd} be satisfied.	
For every $f \in L^{\infty}(X)$ let $\tilde{\Pi}_f$ be the operator defined in \eqref{Def-Para-new} with $\min(\alpha_1,\beta_1,\alpha_2,\beta_2)>\frac{n+D+2}{2m}$. 
Then there exists some constant $C>0$ such that for every $f \in L^{\infty}(X)$ and every $g \in L^2(X)$ 
\begin{align*}
		\norm{\tilde{\Pi}_f(g)}_{L^2(X)} \leq C \norm{f}_{L^{\infty}(X)} \norm{g}_{L^2(X)}.
\end{align*} 
\end{Theorem}

For the proof, let us first define suitable approximations of $\tilde{\Pi}_f$. 
Let $f \in L^{\infty}(X)$ be fixed.
We define for every $R>0$ the operator $T_R: L^2(X) \to L^2(X)$ by
\begin{align} \label{approx-op}
	T_R(g):= \int_{1/R}^R \psi_1(t^{2m}L) \Eins_{B(0,R)} \psi_2(t^{2m}L)[e^{-t^{2m}L}g \cdot e^{-t^{2m}L}f] \,\frac{dt}{t}
\end{align}
for every $g \in L^2(X)$.\\

For convenience, let us set $\delta:=\min(\alpha_1,\beta_1,\alpha_2,\beta_2)$. Let $\psi \in \Psi_{\delta,\delta}(\Sigma_{\sigma}^0)$ and choose $\phi \in H^{\infty}(\Sigma_{\sigma}^0)$ according to the assumptions of Theorem \ref{weak-T1-Theorem}, such that $\{\phi(tL)\}_{t>0}$ satisfies off-diagonal estimates of order $\delta$. One can show the following off-diagonal estimates.

\begin{Lemma} \label{off-diag-para1}
Let $f \in L^{\infty}(X)$ and let $R>0$. 
For every $\gamma$ with $0<\gamma<\delta$ there exists some constant $C>0$, independent of $R>0$, such that for arbitrary open sets $E,F$ in $X$, all $g \in L^2(X)$ with $\supp g \subseteq E$ and all $f \in L^{\infty}(X)$
\begin{align*}
 	\norm{\psi(tL)T_Rg}_{L^2(F)} + \norm{\phi(tL)T_R\psi(tL)g}_{L^2(F)} \leq C \left(1+\frac{\dist(E,F)^{2m}}{t}\right)^{-\gamma} \norm{f}_{L^{\infty}(X)} \norm{g}_{L^2(E)}.
\end{align*}
\end{Lemma}

We postpone the proof of the lemma to Section \ref{sect-Proofs}.

\begin{Proof}[of Theorem \ref{thm-new-para}]
We apply Theorem \ref{weak-T1-Theorem} to the approximation operators $T_R$ defined in \eqref{approx-op}. 
First observe that due to the uniform boundedness of the operator families $\{\psi_1(tL)\}_{t>0}, \; \{\psi_2(tL)\}_{t>0}, \; \{e^{-tL}\}_{t>0}$ on $L^2(X)$ and of $\{e^{-tL}\}_{t>0}$ on $L^{\infty}(X)$, every operator $T_R$ is bounded on $L^2(X)$ with the operator norm bounded by some constant depending on $R>0$. \\
Using Lemma \ref{DaviesGaffneyComp}, we obtain from Lemma \ref{off-diag-para1} the required off-diagonal estimates for the operator families $\{\psi(tL)T_R\phi(tL)\}_{t>0}$ and $\{\psi(tL^{\ast})T_R^{\ast} \phi(tL^{\ast})\}_{t>0}$ with constants independent of $R>0$. \\
It remains to check that
\begin{align} \label{T1-Bed-Para}
		\sup_{R>0} \norm{T_R(1)}_{BMO_L(X)} < \infty
		\qquad \text{and}  \qquad
		\sup_{R>0} \norm{T^{\ast}_R(1)}_{BMO_{L^{\ast}}(X)} < \infty.
\end{align}
Starting with the first assertion, let us define for every $h \in H^1_{L^{\ast}}(X)$ a function $H$ by 
\begin{align*}
		H(x,t):= \psi_1(t^{2m}L^{\ast}) h(x) , \qquad (x,t) \in X \times (0,\infty).
\end{align*}
Since $\psi_1 \in \Psi_{\beta_1,\alpha_1}(\Sigma_{\sigma}^0)$ with $\alpha_1>\frac{n}{4m}$, Theorem \ref{Charact-Hp} yields that $H \in T^1(X)$ with $\norm{H}_{T^1(X)} \approx \norm{h}_{H^1_{L^{\ast}}(X)}$. 
Using that $\psi_2 \in \Psi_{\beta_2,\alpha_2}(\Sigma_{\sigma}^0)$ with $\beta_2>\frac{n}{4m}$, there holds on the other hand that the function $F$, defined by
\begin{align*}
		F(x,t):=\psi_2(t^{2m}L) e^{-t^{2m}L} f(x), \qquad (x,t) \in X \times (0,\infty),
\end{align*}
is according to Proposition \ref{LemmaCarlesonMeasure} an element of $T^{\infty}(X)$ with $\norm{F}_{T^{\infty}(X)} \lesssim \norm{f}_{BMO_L(X)}$. Due to the assumption $e^{-tL}(1)=1$ and \cite{DuongYan2}, Proposition 2.5 (generalized to our setting), there actually holds $L^{\infty}(X) \subseteq BMO(X) \subseteq BMO_L(X)$ and therefore $\norm{F}_{T^{\infty}(X)} \lesssim \norm{f}_{L^{\infty}(X)}$.
Again taking into account the assumption $e^{-tL}(1)=1$ in $L^{\infty}(X)$, we thus obtain
\begin{align*}
		\skp{T_R(1),h} &= \int_{1/R}^R \skp{\Eins_{B(0,R)} \psi_2(t^{2m}L)[e^{-t^{2m}L}1 \cdot e^{-t^{2m}L}f], \psi_1(t^{2m}L^{\ast}) h} \,\frac{dt}{t} \\
				& = \int_{1/R}^R \skp{\Eins_{B(0,R)} \psi_2(t^{2m}L) e^{-t^{2m}L}f, \psi_1(t^{2m}L^{\ast}) h} \,\frac{dt}{t} 
				 = \int_{1/R}^R \skp{\Eins_{B(0,R)} F(x,t), H(x,t)} \,\frac{dt}{t}.
\end{align*}
The duality of $T^1(X)$ and $T^{\infty}(X)$, cf. \cite{CoifmanMeyerStein}, Theorem 1, then yields that
\begin{align*}
		\abs{\skp{T_R(1),h}} \lesssim \norm{F}_{T^{\infty}(X)} \norm{H}_{T^1(X)}
				\lesssim \norm{f}_{L^{\infty}(X)} \norm{h}_{H^1_{L^{\ast}}(X)},
\end{align*}
where the implicit constants are independent of $R>0$. 
Due to the duality of $H^1_{L^{\ast}}(X)$ and $BMO_L(X)$, see Theorem \ref{Duality}, we finally obtain that $T_R(1) \in BMO_L(X)$ with
\begin{align*}
		\sup_{R>0} \norm{T_R(1)}_{BMO_L(X)} \lesssim \norm{f}_{L^{\infty}(X)}.
\end{align*}
Coming to the second assertion in \eqref{T1-Bed-Para}, observe that $\psi_1(tL^{\ast})(1)=0$ in $L^2_{\loc}(X)$ due to the assumption $e^{-tL^{\ast}}(1)=1$ in $L^2_{\loc}(X)$ and Lemma \ref{psi-remark}. 
Thus, 
\begin{align*}
		T_R^{\ast}(1)  = \int_{1/R}^R e^{-t^{2m}L^{\ast}} [\psi_2(t^{2m}L^{\ast})\Eins_{B(0,r)}\psi_1(t^{2m}L^{\ast})(1)
				 \cdot \overline{e^{-t^{2m}L}f}]\,\frac{dt}{t} 
				 		 =0
\end{align*}
in $L^2(X)$ and therefore also in $BMO_{L^{\ast}}(X)$.
\end{Proof}

 \subsection{Towards a $T(b)$-Theorem}
In this section, we give a criterion, under which a $T(b)$-Theorem in our setting holds.

\begin{Def} \label{Def-accretive}
  A function $b \in L^{\infty}(X)$ is said to be \emph{accretive} if there exists a constant $c_0>0$ such that $\Re b(x) \geq c_0$ for almost all $x \in X$. 
\end{Def}

We first state two auxiliary results, Lemmata \ref{Tb-Lemma} and \ref{Carleson-Tb}. They represent the major changes in the proof of Theorem \ref{Tb-Theorem} below in comparision to the proof of the $T(1)$-Theorem, Theorem \ref{T1-Theorem}. Their proofs are postponed to Section 5.\\ 
For every $b \in L^{\infty}(X)$, we denote by $M_b$ the multiplication operator defined by $M_b f := b \cdot f$ for all measurable functions $f: X \to \C$.

\begin{Lemma} \label{Tb-Lemma}
Let $\alpha, \beta\geq 1$, $\psi \in \Psi_{\beta,\alpha}(\Sigma_{\sigma}^0)$ and let $T: \calD(L) \cap \calR(L) \to L^2_{\loc}(X)$ be a linear operator such that the operator family $\{T\psi(tL)\}_{t>0}$ satisfies weak off-diagonal estimates of order $\gamma>\frac{n}{2m}$.
 Moreover, let $b \in L^{\infty}(X)$ and assume that there exist $\delta>0$ and $\tilde{\psi} \in \Psi_{\beta_1,\alpha_1}(\Sigma_{\sigma}^0)$ for some $\alpha_1 \geq \alpha$ and $\beta_1 \geq \beta$ such that  $\int_0^{\infty} \psi(t) \tilde{\psi}(t) \, \frac{dt}{t} =1$ and such that there exists some constant $C>0$ with
\begin{equation} \label{schur-est}
 	\norm{\tilde{\psi}(sL) M_b \psi(tL)f}_{L^2(X)} 
		\leq C \min\left( \frac{s}{t}, \frac{t}{s} \right)^{\delta} \norm{f}_{L^2(X)} \norm{b}_{L^{\infty}(X)}
\end{equation}
for all $s,t>0$ and all $f \in L^2(X)$. 
Additionally, assume that there exists some $\eps_0 \in (0,1)$ such that $\eps_0 \beta>\gamma$ and   $(1-\eps_0)\delta>\frac{n}{2m} + \gamma$.\\
Then the operator family $\{TM_b\psi(tL)\}_{t>0}$, originally defined by \eqref{tb-eq1}, satisfies weak off-diagonal estimates of order $\gamma$.
\end{Lemma}

\begin{Remark}
 If one replaces the \emph{weak} off-diagonal estimates by off-diagonal estimates in Lemma \ref{Tb-Lemma}, one no longer needs the assumption $\gamma>\frac{n}{2m}$. Also the assumption $(1-\eps_0)\delta>\frac{n}{2m} + \gamma$ reduces to $(1-\eps_0)\delta>\gamma$. 
The proof in this case follows the one of Lemma \ref{Tb-Lemma} (cf. Section 5), replacing the splitting of $X$ into balls of radius $t$ by a splitting into two complementary sets, as it is done in the proof of Lemma \ref{DaviesGaffneyComp} and Lemma \ref{Carleson-Tb} below.
\end{Remark}

\begin{Lemma} \label{Carleson-Tb}
Let $\alpha>0$, $\beta>\frac{n}{4m}+[\frac{n}{4m}]+1$ and $\psi \in \Psi_{\beta,\alpha}(\Sigma_{\sigma}^0)$. Let $b \in L^{\infty}(X)$ and  assume that there exist $\delta>\beta$ and $\tilde{\psi} \in \Psi_{\beta_1,\alpha_1}(\Sigma_{\sigma}^0)$ for some $\alpha_1 \geq \alpha$ and $\beta_1 \geq \beta$ such that  $\int_0^{\infty} \psi(t) \tilde{\psi}(t) \,\frac{dt}{t}=1$ and such that \eqref{schur-est} is satisfied with $b$ replaced by $\overline{b}$.
Additionally, assume that there exists some $\eps_0 \in (0,1)$ with $\eps_0\beta>\frac{n}{4m}$ and  $(1-\eps_0)\delta>\eps_0\beta+[\frac{n}{4m}]+1$. \\
Then for every $f \in BMO_L(X)$ is
\begin{align*}
		\nu_{\psi,f}:=\abs{\psi(t^{2m}L)M_bf(y)}^2 \frac{d\mu(y)dt}{t}
\end{align*}
a Carleson measure and there exists a constant $C_{\psi}>0$ such that for all $f \in BMO_{L}(X)$
\begin{align*}
		\norm{\nu_{\psi,f}}_{\calC} \leq C_{\psi} \norm{b}_{L^{\infty}(X)}^2 \norm{f}^2_{BMO_L(X)}.
\end{align*}
\end{Lemma}

For $b=1$, one can show, in analogy to Lemma \ref{off-diag-st}, that condition \eqref{schur-est} is satisfied for $\delta=\min(\alpha,\beta,\alpha_1,\beta_1)$.
For arbitrary $b \in L^{\infty}(X)$, such an estimate is no longer obvious. We give in Proposition \ref{schur-suff-cond} below sufficient conditions for \eqref{schur-est} with $\delta=\frac{1}{2m}$. Unfortunately, one needs $\delta>\frac{1}{2m}$ in \eqref{schur-est} for the proof of Theorem \ref{Tb-Theorem} below.\\
 We first require the following assumption, which is a slight modification of assumption \eqref{P}.

\paragraph{Assumption (P1)} Let $\psi \in \Psi(\Sigma_{\sigma}^0)$ be given.  Assume that for every $f \in L^2(X)$ there exists a measurable function $g: X \times (0,\infty) \to \C$ such that
for all $s>0$ there holds $g_s:=g(\,.\,,s) \geq 0$, and the pair $(\psi(s^{2m}L)f,g_s)$ satisfies a $p$-Poincar\'{e} inequality of the form \eqref{Poincare-pair} for some $p<2$ and with constants  $\lambda \geq 1, \; C_P>0$ independent of $s$ and $f$. Moreover, assume that for all $s>0$ there holds $g_s \in L^2(X)$ and $\norm{sg_s}_{L^2(X)} \leq C \norm{f}_{L^2(X)}$ with $C>0$ independent of $f$ and  $s$.\\

If $X$ is the Euclidean space $\R^n$, then the pair $(\psi(s^{2m}L),g_s)$ corresponds to $(\psi(s^{2m}L),\nabla \psi(s^{2m}L))$ and (P1) is satisfied whenever $\{s\nabla\psi(s^{2m}L)\}_{s>0}$ is uniformly bounded on $L^2(X)$.

\begin{Prop} \label{schur-suff-cond}
Let $b \in L^{\infty}(X)$, let $e^{-tL}(b)=b$ and $e^{-tL^{\ast}}(\bar{b})=\bar{b}$ in $L^2_{\loc}(X)$ for all $t>0$.
Let $\alpha>0, \beta>\frac{n+D+2}{2m}$. Let $\psi, \tilde{\psi} \in \Psi_{\beta,\alpha}(\Sigma_{\sigma}^0)\setminus \{0\}$.
Let (P1) be satisfied for $\psi$, $L$ and for $\tilde{\psi}$, $L^{\ast}$. 
Then there exists some $C>0$ such that
\[
 	\norm{\tilde{\psi}(tL)M_b\psi(sL)}_{L^2(X) \to L^2(X)} \leq C \min\left(\frac{s}{t},\frac{t}{s}\right)^{\frac{1}{2m}} \norm{b}_{L^{\infty}(X)}.
\]
 \end{Prop}

\begin{Proof}
We follow the proof of Proposition \ref{mainTerm}. 
Let $f \in L^2(X)$ and $s,t>0$. Due to the assumptions  $\beta>\frac{n}{4m}$ and $e^{-tL}(b)=b$ in $L^2_{\loc}(X)$ for all $t>0$, Lemma \ref{psi-remark} yields $\tilde{\psi}(tL)(b)=0$ in $L^2_{\loc}(X)$ for all $t>0$. Thus, for every $R>0$,
\[
 	 \norm{\tilde{\psi}(t^{2m}L)M_b\psi(s^{2m}L)f}_{L^2(B(0,R))}^2 
		 = \norm{\tilde{\psi}(t^{2m}L)M_b\psi(s^{2m}L)f - \tilde{\psi}(t^{2m}L)b \cdot A_t\psi(s^{2m}L)f}_{L^2(B(0,R))}^2. 
\]
Due to Proposition \ref{H-inf-offdiag}, $\{\tilde{\psi}(tL)\}_{t>0}$ satisfies off-diagonal estimates of order $\beta>\frac{n+D+2}{2m}$.
Observe that the calculations in \eqref{mainTerm-eq2}-\eqref{mainTerm-eq4} do not depend on the special form of $(u,g_t)$ but only on the fact that they satisfy a $p$-Poincar\'e inequality. Hence, in analogy to \eqref{mainTerm-eq5}, we obtain
\begin{align*}
	& \norm{\tilde{\psi}(t^{2m}L)M_b\psi(s^{2m}L)f}_{L^2(X)}^2 \\
		 & \qquad \lesssim \sum_{\alpha \in I_{k_0}}  \sum_{\beta \in I_{k_0}}
					\left(1+\frac{\dist(B_{\alpha},B_{\beta})^{2m}}{t^{2m}}\right)^{-\gamma}
					\norm{\psi(s^{2m}L)f-\skp{\psi(s^{2m}L)f}_{Q_{\alpha}^{k_0}}}_{L^2(Q_{\beta}^{k_0})}^2 \\
		& \qquad \lesssim t^2 \norm{\calM_p g_s}_{L^2(X)} 
					\lesssim \left(\frac{t}{s}\right)^2 \norm{sg_s}_{L^2(X)}^2 \lesssim \left(\frac{t}{s}\right)^2 \norm{f}_{L^2(X)}^2, 
\end{align*}
where the last step is a consequence of (P1). The corresponding estimate against $\frac{s}{t}$ follows by duality. 
\end{Proof}

We are now ready to state the $T(b)$-Theorem. 

 \begin{Theorem}  \label{Tb-Theorem}
Let $L$ be an operator satisfying the assumptions \eqref{H1}, \eqref{H2} and \eqref{H3}. Additionally, let the assumptions \eqref{P} and \eqref{Pd} be satisfied.\\
 Let $T: \calD(L) \cap \calR(L) \to L^2_{\loc}(X)$ be a linear operator with $T^{\ast}: \calD(L^{\ast}) \cap \calR(L^{\ast}) \to L^2_{\loc}(X)$
 such that the assumptions \emph{\ref{OD1}} and \emph{\ref{OD2}} are satisfied for some $\gamma>\frac{n+D+2}{2m}$. \\
Let $b_1,b_2 \in L^{\infty}(X)$ be two accretive functions such that the assumptions of Lemma \ref{Tb-Lemma} are satisfied for the operator families  $\{T\psi_1(tL)\}_{t>0}$ with $b_1$ and for $\{T^{\ast}\psi_2(tL^{\ast})\}_{t>0}$ with $\overline{b}_2$ and such that the assumptions of Lemma \ref{Carleson-Tb} are satisfied for the triples $\psi_1,\overline{b}_1,L^{\ast}$ and $\psi_2,b_2,L$. \\
Moreover, let $T(b_1) \in BMO_L(X)$ and $T^{\ast}(\overline{b}_2) \in BMO_{L^{\ast}}(X)$. \\
 Then $T$ is bounded on $L^2(X)$, i.e. there exists a constant $C>0$ 
 such that for all $f \in L^2(X)$ 
 \[
  	\norm{Tf}_{L^2(X)} \leq C \norm{f}_{L^2(X)}.
 \]
 \end{Theorem}

\begin{Proof}[of Theorem \ref{Tb-Theorem}]
The proof works analogously to the one of Theorem \ref{T1-Theorem}. We will not give the proof in all details, but only state the differences to the one of Theorem \ref{T1-Theorem}.\\
Let $f,g \in L^2(X)$. Let $b_1,b_2 \in L^{\infty}(X)$ be the two accretive functions given in the assumption with constants $c_1$ and $c_2$, respectively. Moreover, let $\alpha\geq 1, \ \beta> \frac{n}{4m} +[\frac{n}{4m}]+1$ and let $\psi_1,\psi_2 \in \Psi_{\beta,\alpha}(\Sigma_{\sigma}^0)\setminus \{0\}$ as given in the assumption. Denote by $\tilde{\psi}_1,\tilde{\psi}_2 \in \Psi(\Sigma_{\sigma}^0)$ the functions given in the assumptions of Lemma \ref{Tb-Lemma} and Lemma \ref{Carleson-Tb} that satisfy   $\int_0^{\infty} \psi_1(t) \tilde{\psi}_1(t) \, \frac{dt}{t} =1$ and $\int_0^{\infty} \psi_2(t) \tilde{\psi}_2(t) \, \frac{dt}{t} =1$.\\
Since $b_1,b_2$ are accretive functions, it will be sufficient to estimate $M_{b_2}TM_{b_1}$ instead of $T$.
In analogy to the proof of Theorem \ref{T1-Theorem}, we first decompose both $f$ and $g$ with the help of the Calder\'{o}n reproducing formula, which yields
\begin{align} \label{Tb-eq1}
  	 \skp{M_{b_2}TM_{b_1}f,g}
 		& =  \int_0^{\infty} \int_0^{\infty} \skp{\psi_2(t^{2m}L) M_{b_2}TM_{b_1} \psi_1(s^{2m}L) \tilde{\psi}_1(s^{2m}L) f, \tilde{\psi}_2(t^{2m}L^{\ast})g} \, \frac{dt}{t} \frac{ds}{s}. 
 \end{align}
The two main differences will be the following.
Observe that due to Lemma \ref{Tb-Lemma} and the assumption $\gamma>\frac{n+D+2}{2m}$, the operator families
\begin{align} \label{Tb-eq2}
			 \{TM_{b_1}\psi_1(tL)\}_{t>0} \qquad \text{and} \qquad \{T^{\ast}M_{\overline{b}_2}\psi_2(tL^{\ast})\}_{t>0}
\end{align}
satisfy weak off-diagonal estimates of order $\gamma$. 
Moreover, together with the assumptions $T(b_1) \in BMO_L(X)$ and $T^{\ast}(\overline{b}_2) \in BMO_{L^{\ast}}(X)$, Lemma \ref{Carleson-Tb} yields that
\begin{align} \label{Tb-eq3}
		 \abs{\psi_2(t^{2m}L)M_{b_2}T(b_1)(y)}^2 \frac{d\mu(y)dt}{t} \qquad \text{and} \qquad
		 \abs{\psi_1(t^{2m}L^{\ast})M_{\overline{b}_1}T^{\ast}(\overline{b}_2)(y)}^2 \frac{d\mu(y)dt}{t}
\end{align}
are Carleson measures.\\
As in the proof of Theorem \ref{T1-Theorem}, it is enough to consider the part $J_1$, where in the inner integral of \eqref{Tb-eq1} one only integrates over the interval $\{t \in (0,\infty) \,:\, 0<t<s\}$.
Then, one also uses the first line of the decomposition \eqref{T-decomposition}, but now applied for the operator $M_{b_2}TM_{b_1}$ instead of $T$.
The error term $J_E$ is then equal to
\begin{align*}
		J_E = \int_0^{\infty} \int_0^s \skp{\psi_2(t^{2m}L) M_{b_2}TM_{b_1} (I-e^{-t^{2m}L}) \psi_1(s^{2m}L) \tilde{\psi}_1(s^{2m}L) f, \tilde{\psi}_2(t^{2m}L^{\ast}) g} \, \frac{dt}{t} \frac{ds}{s}.
\end{align*}
Due to the weak off-diagonal estimates for the operator family $\{T^{\ast}M_{\overline{b}_2}\psi_2(tL^{\ast})\}_{t>0}$ and the fact that
\begin{align*}
		\norm{M_{b_1}(I-e^{-t^{2m}L})\psi_1(s^{2m}L) h}_{L^2(X)} \lesssim \norm{(I-e^{-t^{2m}L})\psi_1(s^{2m}L) h}_{L^2(X)},
\end{align*}
we can simply copy the estimates in \eqref{T1-eq2a}, \eqref{T1-eq3} and \eqref{T1-eq3a} and obtain
\begin{align*}
			\abs{J_E} \lesssim \norm{f}_{L^2(X)} \norm{g}_{L^2(X)}.
\end{align*}
To handle the main term $J_M$, we now split $M_{b_2}TM_{b_1}e^{-t^{2m}L}$ into 
\begin{align*}
		[M_{b_2}TM_{b_1}e^{-t^{2m}L} - M_{b_2}T(b_1) \cdot A_te^{-t^{2m}L}] +M_{b_2}T(b_1) \cdot A_t e^{-t^{2m}L}.
\end{align*}
Then, following the same procedure as in \eqref{T1-eq4a}, we get $J_M^0 = J_M^1 + J_M^2$ with 
\begin{align*}
			J_M^1 = \int_0^{\infty} \skp{\psi_2(t^{2m}L) M_{b_2}TM_{b_1}e^{-t^{2m}L}f - \psi_2(t^{2m}L) M_{b_2}T(b_1) \cdot A_te^{-t^{2m}L}f, \tilde{\psi}_2(t^{2m}L^{\ast})g} \,\frac{dt}{t}
\end{align*}
and 
\begin{align*}
			J_M^2 = \int_0^{\infty} \skp{\psi_2(t^{2m}L) M_{b_2}T(b_1) \cdot A_te^{-t^{2m}L}f,\tilde{\psi}_2(t^{2m}L^{\ast})g} \,\frac{dt}{t}.
\end{align*}
The term $J_M^1$ can again be estimated by application of Proposition \ref{mainTerm}, with a slight modification. We set $S_{t^{2m}}:=\psi_2(t^{2m}L) M_{b_2}T$ and observe that this operator satisfies  weak off-diagonal estimates of order $\gamma>\frac{n+D+2}{2m}$ via \eqref{Tb-eq2}.  It remains to check that the constant function $1$ in Proposition \ref{mainTerm} can be replaced by some arbitrary function $b_1 \in L^{\infty}(X)$, i.e. that one can also obtain the estimate
\begin{align*}
	 \int_0^{\infty} \norm{S_{t^{2m}} M_{b_1} e^{-t^{2m}L} f - S_{t^{2m}}(b_1) \cdot A_t e^{-t^{2m}L} f}_{L^2(X)}^2 \,\frac{dt}{t}
 		\leq C \norm{b_1}^2_{L^{\infty}(X)} \norm{f}_{L^2(X)}^2.
\end{align*}
This can easily be seen in the calculations of \eqref{mainTerm-eq2}, where one can pull the function $b_1$ out of the $L^2$-norm in the last step.\\
The term $J_M^2$  is, up to the multiplication operator $M_{b_2}$, a paraproduct. To handle this term, let us have a short look at the proof of  Theorem \ref{paraproduct} (cf. \cite{Frey}, Theorem 4.2), which states the boundedness of paraproducts on $L^2(X)$. There, one only exploits the fact that $\abs{\psi(t^{2m}L)b(y)}^2 \frac{d\mu(y)dt}{t}$ is a Carleson measure whenever $b \in BMO_L(X)$ and does not explicitly use that $b \in BMO_L(X)$. Since we have by assumptions that $\abs{\psi_2(t^{2m}L)M_{b_2}T(b_1)(y)}^2 \frac{d\mu(y)dt}{t}$ is a Carleson measure, see \eqref{Tb-eq3}, we also get the desired estimate for $J_M^2$.\\
The proof of the remainder term $J_R$, defined by 
\begin{align*}
	J_R :=	\int_0^{\infty} \int_s^{\infty} \skp{\psi_2(t^{2m}L) M_{b_2}TM_{b_1} e^{-t^{2m}L}  \psi_1(s^{2m}L) \tilde{\psi}_1(s^{2m}L)f, \tilde{\psi}_2(t^{2m}L^{\ast})g} \, \frac{dt}{t} \frac{ds}{s},
\end{align*}
is again handled as in \eqref{T1-eq8}, with the same changes as those for the treatment of $J_E$.
\end{Proof}

 \section{Proofs of auxiliary results}
 \label{sect-Proofs}

In this section, we prove Lemmata \ref{self-improving-est}, \ref{off-diag-para1}, \ref{Tb-Lemma} and \ref{Carleson-Tb}. In all cases, one has to establish off-diagonal for certain types of operators.
Except for Lemma \ref{off-diag-para1}, one always transfers  off-diagonal estimates from one operator to another with the help of a Calder\'{o}n reproducing formula.

\begin{Proof}[of Lemma \ref{self-improving-est}]
 Let $t>0$ and let $B_1,B_2$ be arbitrary balls with radius $t$. Let $f,g \in L^2(X)$ with $\supp f \subseteq B_1$ and $\supp g \subseteq B_2$.
Given $\psi \in \Psi(\Sigma_{\sigma}^0) \setminus \{0\}$ from the assumptions, we choose some function $\tilde{\psi} \in \Psi_{\alpha,\beta}(\Sigma_{\sigma}^0)$ with $\alpha, \beta>\gamma$ and $\int_0^{\infty} \psi(s) \tilde{\psi}(s) \,\frac{ds}{s}=1$.
The Calder\'{o}n reproducing formula then yields
\begin{align} \label{int-eq}
 	\skp{T\varphi(t^{2m}L)f,g}
			= \int_0^{\infty} \skp{T\varphi(t^{2m}L) \psi(s^{2m}L) \tilde{\psi}(s^{2m}L) f,g} \,\frac{ds}{s}.
\end{align}
Since $\{T\varphi(tL)\}_{t>0}$ is uniformly bounded on $L^2(X)$, we can without restriction assume that $\dist(B_1,B_2)\leq t$.
We break the integral in \eqref{int-eq} into two parts, one over $(0,t)$, which is called $J_1$, and one over $(t,\infty)$, which is called $J_2$.\\
We first turn to $J_1$. On the one hand, $\{T\psi(sL)\}_{s>0}$ satisfies weak off-diagonal estimates of order $\gamma>\frac{n}{2m}$. Proposition \ref{H-inf-offdiag}, on the other hand, yields that $\{\tilde{\psi}(sL)\varphi(tL)\}_{s,t>0}$ satisfies off-diagonal estimates in $s$ of order $\alpha$, since $\sup_{t>0} \norm{\varphi(t \,\cdot\,)}_{L^{\infty}(\Sigma_{\sigma}^0)} = \norm{\varphi}_{L^{\infty}(\Sigma_{\sigma}^0)} < \infty$.
Hence, the composition of the two operators $\{T\psi(sL)\tilde{\psi}(sL)\varphi(tL)\}_{s,t>0}$ satisfies weak off-diagonal estimates in $s$ of order $\min(\gamma,\alpha)=\gamma>\frac{n}{2m}$ on $L^2(X)$ due to Proposition \ref{weak-composition}.
Using Remark \ref{weak-estimates-st} (which provides us with weak off-diagonal estimates in $s$ for balls of radius $t>s$), with the roles of $s$ and $t$ interchanged, we therefore get
\begin{align*}
 	\abs{J_1}
		& \leq \int_0^{t} \abs{\skp{T\psi(s^{2m}L) \tilde{\psi}(s^{2m}L) \varphi(t^{2m}L) f,g}} \,\frac{ds}{s} \\
		& \lesssim \int_0^{t} \left(\frac{t}{s}\right)^n \left(1+\frac{\dist(B_1,B_2)^{2m}}{s^{2m}}\right)^{-\gamma} \,\frac{ds}{s}
				\norm{f}_{L^2(B_1)} \norm{g}_{L^2(B_2)}.
\end{align*}
Since we assumed $\gamma>\frac{n}{2m}$ and $\dist(B_1,B_2)>t$, there further holds
\begin{align*}
 	& \int_0^{t} \left(\frac{t}{s}\right)^n \left(1+\frac{\dist(B_1,B_2)^{2m}}{s^{2m}}\right)^{-\gamma} \,\frac{ds}{s} 
		 \leq \int_0^{t} \left(\frac{t}{s}\right)^{n-2m\gamma} \left(\frac{\dist(B_1,B_2)^{2m}}{t^{2m}}\right)^{-\gamma} \,\frac{ds}{s} \\
		& \qquad =  \left(\frac{\dist(B_1,B_2)^{2m}}{t^{2m}}\right)^{-\gamma} \int_0^1 u^{2m\gamma-n} \,\frac{du}{u} 
		 \lesssim \left(1+\frac{\dist(B_1,B_2)^{2m}}{t^{2m}}\right)^{-\gamma},
\end{align*}
using the substitution $u=\frac{s}{t}$ in the penultimate step.\\
We turn to $J_2$. We again use that $\{T\psi(sL)\}_{s>0}$ satisfies weak off-diagonal estimates in $s$ of order $\gamma$.
 Lemma \ref{off-diag-st} shows that for every $a>0$ with $a \leq \delta$ and $a<\beta$, there exists a family of operators $\{T_{s,t}\}_{s,t>0}$ such that
\[
 	\tilde{\psi}(sL) \varphi(tL) = \left(\frac{t}{s}\right)^a T_{s,t},
\]
where $\{T_{s,t}\}_{s,t>0}$ satisfies off-diagonal estimates in $s$ of order $\alpha+a$ uniformly in $t>0$.
Proposition \ref{weak-composition} then yields that the family of operators $\{T\psi(sL)T_{s,t}\}_{s,t>0}$ satisfies weak off-diagonal estimates in $s$ of order $\min(\gamma,\alpha+a)=\gamma$. 
Since weak off-diagonal estimates can be applied to smaller balls without any change, we obtain
\begin{align*}
 	\abs{J_2}
		& \leq \int_{t}^{\infty} \abs{\skp{T\psi(s^{2m}L)\tilde{\psi}(s^{2m}L)\varphi(t^{2m}L)f,g}} \,\frac{ds}{s} \\
		& \lesssim \int_{t}^{\infty} \left(\frac{t}{s}\right)^{2ma} \left(1+\frac{\dist(B_1,B_2)^{2m}}{s^{2m}}\right)^{-\gamma} \,\frac{ds}{s} \norm{f}_{L^2(B_1)} \norm{g}_{L^2(B_2)}.
\end{align*}
Since we assumed $\delta>\gamma$ and $\beta>\gamma$, we can fix some $a>\gamma$ with $a\leq \delta$ and $a<\beta$. 
For this choice of $a$ we further get, similar to the treatment of $J_1$, 
\begin{align*}
 	& \int_{t}^{\infty} \left(\frac{t}{s}\right)^{2ma} \left(1+\frac{\dist(B_1,B_2)^{2m}}{s^{2m}}\right)^{-\gamma} \,\frac{ds}{s} 
		 \leq \int_{t}^{\infty} \left(\frac{t}{s}\right)^{2m(a-\gamma)} \left(\frac{\dist(B_1,B_2)^{2m}}{t^{2m}}\right)^{-\gamma} \,\frac{ds}{s} \\
		& \qquad	= \left(\frac{\dist(B_1,B_2)^{2m}}{t^{2m}}\right)^{-\gamma} \int_1^{\infty} u^{-2m(a-\gamma)} \,\frac{du}{u} 
		 \lesssim \left(1+\frac{\dist(B_1,B_2)^{2m}}{t^{2m}}\right)^{-\gamma},
\end{align*}
still assuming that $\dist(B_1,B_2)>t$.
Combining the estimates of $J_1$ and $J_2$ finishes the proof.
\end{Proof}

\begin{Proof}[of Lemma \ref{off-diag-para1}]
Let $E,F$ be two arbitrary open sets in $X$ and let $f \in L^{\infty}(X)$ and $g \in L^2(X)$ with $\supp g \subseteq E$. We begin with the estimation of $\{\psi(tL)T_R\}_{t>0}$.
Let $\delta=\min(\alpha_1,\beta_1,\alpha_2,\beta_2)$ as defined before and fix some $\gamma>0$ with $\gamma<\delta$. Then for every $s,t>0$
\begin{align} \label{schur-cond}
		\norm{\psi(tL) \psi_1(sL)}_{L^2(X) \to L^2(X)}
				 \lesssim \min\left(\frac{s}{t}, \frac{t}{s} \right)^{\delta},
\end{align}
using the same arguments as e.g. in Lemma \ref{off-diag-st}. Hence, due to Minkowski's inequality and the uniform boundedness of the operator families $\{\psi_2(sL)\}_{s>0}$, $\{e^{-sL}\}_{s>0}$ on $L^2(X)$ and $\{e^{-sL}\}_{s>0}$ on $L^{\infty}(X)$ we obtain
\begin{align*}
		\norm{\psi(t^{2m}L)T_R(g)}_{L^2(X)} 
				& \leq \int_0^{\infty} \norm{\psi(t^{2m}L) \psi_1(s^{2m}L) \Eins_{B(0,R)} \psi_2(s^{2m}L)[e^{-s^{2m}L}g \cdot e^{-s^{2m}L} f]}_{L^2(X)} \,\frac{ds}{s} \\
				&  \lesssim \int_0^{\infty} \min\left(\frac{s}{t}, \frac{t}{s} \right)^{2m\delta} \,\frac{ds}{s} \norm{g}_{L^2(X)} \norm{f}_{L^{\infty}(X)}
				\lesssim \norm{g}_{L^2(X)}  \norm{f}_{L^{\infty}(X)}.
\end{align*}
If $\dist(E,F)\leq t$, the above estimate yields the desired conclusion.
Otherwise, let $\rho:=\dist(E,F)>t$, and define $G_1:=\{x \in X\,:\, \dist(x,F)<\frac{\rho}{2}\}$ and $G_2:=\{x \in X\,:\, \dist(x,F)<\frac{\rho}{4}\}$. Then there holds that $G_1,G_2$ are open with $\dist(E,G_1) \geq \frac{\rho}{2}$ and $\dist(F,X\setminus \bar{G}_2) \geq \frac{\rho}{4}$. We split $X$ into $X=\bar{G}_2 \cup X \setminus \bar{G}_2$ and obtain
\begin{align*}
	&	\norm{\psi(t^{2m}L)T_R(g)}_{L^2(F)} \\
				& \quad \leq \int_0^{\infty}  \norm{\psi(t^{2m}L) \psi_1(s^{2m}L) \Eins_{B(0,R)} \psi_2(s^{2m}L)\Eins_{\bar{G}_2}[e^{-s^{2m}L}g \cdot e^{-s^{2m}L} f]}_{L^2(F)} \,\frac{ds}{s} \\
				& \quad \qquad + \int_0^{\infty}  \norm{\psi(t^{2m}L) \psi_1(s^{2m}L) \Eins_{B(0,R)} \psi_2(s^{2m}L)\Eins_{X\setminus \bar{G}_2} [e^{-s^{2m}L}g \cdot e^{-s^{2m}L} f]}_{L^2(F)} \,\frac{ds}{s} \\
				& \quad =:J_{\bar{G}_2} + J_{X\setminus \bar{G}_2}.
\end{align*}
The estimate \eqref{schur-cond}, the uniform boundedness of $\{\psi_2(sL)\}_{s>0}$ on $L^2(X)$ and $\bar{G}_2 \subseteq G_1$ yield
\begin{align*}
		J_{\bar{G}_2}
			&  \lesssim \int_0^{\infty}  \min\left(\frac{s}{t}, \frac{t}{s} \right)^{2m\delta}  
							\norm{\psi_2(s^{2m}L)\Eins_{\bar{G}_2}[e^{-s^{2m}L}g \cdot e^{-s^{2m}L} f]}_{L^2(X)} \,\frac{ds}{s} \\
		 & \lesssim \int_0^{\infty}  \min\left(\frac{s}{t}, \frac{t}{s} \right)^{2m\delta}  
							\norm{e^{-s^{2m}L}g \cdot e^{-s^{2m}L} f}_{L^2(G_1)} \,\frac{ds}{s}.
\end{align*}
Since, on the one hand, $\{e^{-sL}\}_{s>0}$ satisfies Davies-Gaffney estimates and is on the other hand uniformly bounded on $L^{\infty}(X)$, we can estimate the above by a constant times $\norm{f}_{L^{\infty}(X)} \norm{g}_{L^2(E)}$ times
\begin{align*}
		&  \int_0^{\infty}  \min\left(\frac{s}{t}, \frac{t}{s} \right)^{2m\delta}  \left(1+\frac{\dist(E,G_1)^{2m}}{s^{2m}}\right)^{-\gamma} \,\frac{ds}{s}\\
		& \qquad \lesssim \int_0^{t}  \left(\frac{s}{t}\right)^{2m\delta} 
					\left(1+\frac{\dist(E,F)^{2m}}{s^{2m}}\right)^{-\gamma} \,\frac{ds}{s} 
		+ \int_t^{\infty} \left(\frac{t}{s}\right)^{2m\delta} \left(1+\frac{\dist(E,F)^{2m}}{s^{2m}}\right)^{-\gamma} \,\frac{ds}{s}\\
		& \qquad \lesssim \left(1+\frac{\dist(E,F)^{2m}}{t^{2m}}\right)^{-\gamma} \left[ \int_0^{t} \left(\frac{s}{t}\right)^{2m\delta}  \,\frac{ds}{s} + \int_t^{\infty} \left(\frac{t}{s}\right)^{2m\delta}  \left(\frac{t}{s}\right)^{-2m\gamma}\,\frac{ds}{s} \right] \\
		& \qquad \lesssim \left(1+\frac{\dist(E,F)^{2m}}{t^{2m}}\right)^{-\gamma},
\end{align*}
using that $\gamma< \delta$. This gives the desired estimate for $J_{\bar{G}_2}$.\\
For the analogous estimate of $J_{X \setminus \bar{G}_2}$, we instead use the off-diagonal estimates of the operator family $\{\psi(tL)\psi_1(sL)\}_{s,t>0}$. We split $J_{X \setminus \bar{G}_2}$ into two parts $J^1_{X \setminus \bar{G}_2}$ and $J^2_{X \setminus \bar{G}_2}$, representing the integration over $(0,t)$ and $(t,\infty)$, respectively. Considering $J^1_{X \setminus \bar{G}_2}$, we take into account that Lemma \ref{off-diag-st} yields off-diagonal estimates in $t$ of order $\gamma$ for the operator family $\{\psi(tL)\psi_1(sL)\}_{s,t>0}$ with an extra term $\left(\frac{s}{t}\right)^{\gamma}$. In addition, $\{\psi_2(sL)\}_{s>0}$ satisfies off-diagonal estimates in $s$ of order $\gamma$ due to Proposition \ref{H-inf-offdiag}. Lemma \ref{DaviesGaffneyComp} then yields that 
\begin{align} \label{appl-para-eq3} \nonumber
 	J^1_{X \setminus \bar{G}_2}& \leq \int_0^t \norm{\psi(t^{2m}L) \psi_1(s^{2m}L) \Eins_{B(0,R)} \psi_2(s^{2m}L)\Eins_{X\setminus \bar{G}_2} [e^{-s^{2m}L}g \cdot e^{-s^{2m}L} f]}_{L^2(F)} \,\frac{ds}{s} \\ \nonumber
	&  \lesssim \int_0^t \left(\frac{s}{t}\right)^{2m\gamma} \left(1+\frac{\dist(F,X \setminus \bar{G}_2)^{2m}}{t^{2m}}\right)^{-\gamma} \norm{e^{-s^{2m}L}g \cdot e^{-s^{2m}L} f}_{L^2(X)} \,\frac{ds}{s} \\
	&  \lesssim \left(1+\frac{\dist(E,F)^{2m}}{t^{2m}}\right)^{-\gamma} \norm{f}_{L^{\infty}(X)} \norm{g}_{L^2(E)}.
\end{align}
For the part $J^2_{X \setminus \bar{G}_2}$, we in turn use that $\{\psi(tL)\psi_1(sL)\}_{s,t>0}$ satisfies off-diagonal estimates in $s$ of order $\gamma_1$ with an extra term $\left(\frac{t}{s}\right)^{\gamma_1}$, where $\gamma<\gamma_1<\delta$. With similar arguments as in \eqref{appl-para-eq3}, we then obtain
\begin{align} \label{appl-para-eq4} \nonumber
 	J^2_{X \setminus \bar{G}_2}
 		& \leq \int_t^{\infty} \norm{\psi(t^{2m}L) \psi_1(s^{2m}L) \Eins_{B(0,R)} \psi_2(s^{2m}L)\Eins_{X\setminus \bar{G}_2} [e^{-s^{2m}L}g \cdot e^{-s^{2m}L} f]}_{L^2(F)} \,\frac{ds}{s} \\ \nonumber
	&  \lesssim \int_t^{\infty} \left(\frac{t}{s}\right)^{2m\gamma_1}  \left(1+\frac{\dist(F,X \setminus \bar{G}_2)^{2m}}{s^{2m}}\right)^{-\gamma_1} \norm{e^{-s^{2m}L}g \cdot e^{-s^{2m}L} f}_{L^2(X)} \,\frac{ds}{s} \\ \nonumber
	&  \lesssim \left(1+\frac{\dist(E,F)^{2m}}{t^{2m}}\right)^{-\gamma} \norm{f}_{L^{\infty}(X)} \norm{g}_{L^2(E)} 
			\int_t^{\infty}  \left(\frac{t}{s}\right)^{2m\gamma_1}  \left(\frac{t}{s}\right)^{-2m\gamma} \,\frac{ds}{s} \\
	&  \lesssim \left(1+\frac{\dist(E,F)^{2m}}{t^{2m}}\right)^{-\gamma} \norm{f}_{L^{\infty}(X)} \norm{g}_{L^2(E)}.
\end{align}
The combination of the above estimates finally yields the desired conclusion.
Observe that all implicit constants in the inequalities are independent of $R>0$. \\
We continue with the estimation of $\{\phi(tL)T_R\psi(tL)\}_{t>0}$.
By definition of $\phi$ there holds $\abs{\phi(z)} = \calO(\abs{z}^{\delta})$ for $\abs{z} \to \infty$. Hence, using similar arguments as in Lemma \ref{st-estimates},
\begin{align} \label{appl-para-eq5}
		\norm{\phi(tL)\psi_1(sL)}_{L^2(X) \to L^2(X)} \lesssim \left(\frac{s}{t}\right)^{\delta},
	\qquad	\norm{e^{-sL}\psi(tL)}_{L^2(X) \to L^2(X)} \lesssim \left(\frac{t}{s}\right)^{\delta}.
\end{align}
We therefore obtain, again using the uniform boundedness of the occuring operator families on $L^2(X)$ and $L^{\infty}(X)$, respectively, 
\begin{align*}
 	& \norm{\phi(t^{2m}L) T_R \psi(t^{2m}L)g}_{L^2(X)} \\
		& \quad \leq \int_0^{\infty} \norm{\phi(t^{2m}L) \psi_1(s^{2m}L) \Eins_{B(0,R)} \psi_2(s^{2m}L) [e^{-s^{2m}L} \psi(t^{2m}L) g \cdot e^{-s^{2m}L} f]}_{L^2(X)} \,\frac{ds}{s} \\
		& \quad \lesssim \int_0^t \left(\frac{s}{t}\right)^{2m\delta} \norm{\psi_2(s^{2m}L) [e^{-s^{2m}L} \psi(t^{2m}L) g \cdot e^{-s^{2m}L} f]}_{L^2(X)} \,\frac{ds}{s} \\
		& \quad \quad + \int_t^{\infty}  \norm{e^{-s^{2m}L} \psi(t^{2m}L) g \cdot e^{-s^{2m}L} f}_{L^2(X)} \,\frac{ds}{s} \\
		& \quad \lesssim \norm{f}_{L^{\infty}(X)} \norm{g}_{L^2(X)} \int_0^{\infty} \min\left(\frac{s}{t},\frac{t}{s}\right)^{2m\delta} \,\frac{ds}{s} 
	 \lesssim \norm{f}_{L^{\infty}(X)} \norm{g}_{L^2(X)}.
\end{align*}
If $\dist(E,F)\leq t$, the above estimate yields the desired conclusion.
Otherwise, with the notation as before, we split $X$ into $X=\bar{G}_2 \cup X \setminus \bar{G}_2$. 
Let us moreover split the integrals into two parts over $(0,t)$ and $(t,\infty)$. Taking into account the fact that $\{e^{-sL}\psi(tL)\}_{s,t>0}$ satisfies  off-diagonal estimates in $t$ of order $\gamma$ and using $\bar{G}_2 \subseteq G_1$ and \eqref{appl-para-eq5}, we then obtain
\begin{align} \label{appl-para-eq7} \nonumber
  J^1_{\bar{G}_2} 
		& \leq \int_0^t \norm{\phi(t^{2m}L) \psi_1(s^{2m}L) \Eins_{B(0,R)} \psi_2(s^{2m}L) \Eins_{\bar{G}_2} [e^{-s^{2m}L} \psi(t^{2m}L) g \cdot e^{-s^{2m}L} f]}_{L^2(X)} \,\frac{ds}{s} \\ \nonumber
		& \lesssim \int_0^t \left(\frac{s}{t}\right)^{2m\gamma} \norm{e^{-s^{2m}L} \psi(t^{2m}L) g \cdot e^{-s^{2m}L} f}_{L^2(G_1)} \,\frac{ds}{s}\\ \nonumber
		& \lesssim \left(1+\frac{\dist(E,G_1)^{2m}}{t^{2m}}\right)^{-\gamma} \int_0^t \left(\frac{s}{t}\right)^{2m\gamma} \,\frac{ds}{s}
				\norm{f}_{L^{\infty}(X)} \norm{g}_{L^2(E)} \\
		& \lesssim \left(1+\frac{\dist(E,F)^{2m}}{t^{2m}}\right)^{-\gamma}\norm{f}_{L^{\infty}(X)} \norm{g}_{L^2(E)}.
\end{align}
Moreover, for $\eps=\delta-\gamma>0$, the operator family $\{e^{-sL}\psi(tL)\}_{s,t>0}$ satisfies off-diagonal estimates in $t$ of order $\gamma$ with an extra factor $\left(\frac{t}{s}\right)^{\eps}$, where we use that
\begin{align*}
 	e^{-sL}\psi(tL) = \left(\frac{t}{s}\right)^{\eps} (sL)^{\eps} e^{-sL}(tL)^{-\eps} \psi(tL)
\end{align*}
and Proposition \ref{H-inf-offdiag}. Thus, 
\begin{align} \label{appl-para-eq8} \nonumber
 J^2_{\bar{G}_2} 
		& \leq \int_t^{\infty} \norm{\phi(t^{2m}L) \psi_1(s^{2m}L) \Eins_{B(0,R)} \psi_2(s^{2m}L) \Eins_{\bar{G}_2} [e^{-s^{2m}L} \psi(t^{2m}L) g \cdot e^{-s^{2m}L} f]}_{L^2(X)} \,\frac{ds}{s} \\ \nonumber
		& \lesssim \int_t^{\infty} \norm{e^{-s^{2m}L} \psi(t^{2m}L) g \cdot e^{-s^{2m}L} f}_{L^2(G_1)} \,\frac{ds}{s} \\ \nonumber
		& \lesssim \left(1+\frac{\dist(E,G_1)^{2m}}{t^{2m}}\right)^{-\gamma} \int_t^{\infty} \left(\frac{t}{s}\right)^{2m\eps} \,\frac{ds}{s} \norm{f}_{L^{\infty}(X)} \norm{g}_{L^2(E)} \\
		& \lesssim \left(1+\frac{\dist(E,F)^{2m}}{t^{2m}}\right)^{-\gamma} \norm{f}_{L^{\infty}(X)} \norm{g}_{L^2(E)}.
\end{align}
Let us turn to the calculation of $J_{X \setminus \bar{G}_2}$. 
We now use that for $\eps=\delta-\gamma>0$
\begin{align*}
 	\phi(tL)\psi_1(sL) = \left(\frac{s}{t}\right)^{\eps} (tL)^{\eps}\phi(tL) (sL)^{-\eps} \psi_1(sL),
\end{align*}
therefore $\{\phi(tL)\psi_1(sL)\}_{s,t>0}$ satisfies off-diagonal estimates in  $t$ of order $\gamma$ with an extra factor $\left(\frac{s}{t}\right)^{\eps}$. Hence,
\begin{align} \label{appl-para-eq9} \nonumber
  J^1_{X \setminus \bar{G}_2}
		& \leq \int_0^t \norm{\phi(t^{2m}L) \psi_1(s^{2m}L) \Eins_{B(0,R)} \psi_2(s^{2m}L) \Eins_{X \setminus \bar{G}_2} [e^{-s^{2m}L} \psi(t^{2m}L) g \cdot e^{-s^{2m}L} f]}_{L^2(X)} \,\frac{ds}{s} \\ \nonumber
		& \lesssim \left(1+\frac{\dist(F,X \setminus \bar{G}_2)^{2m}}{t^{2m}}\right)^{-\gamma}  \norm{f}_{L^{\infty}(X)} \norm{g}_{L^2(E)}\int_0^t \left(\frac{s}{t}\right)^{2m\eps} \,\frac{ds}{s} \\
		& \lesssim \left(1+\frac{\dist(E,F)^{2m}}{t^{2m}}\right)^{-\gamma} \norm{f}_{L^{\infty}(X)} \norm{g}_{L^2(E)}.
\end{align}
For the remaining part, we apply \eqref{appl-para-eq5} and off-diagonal estimates of $\{\phi(tL)\psi_1(sL)\}_{s,t>0}$ in $s$ of order $\gamma$, which yields
\begin{align} \label{appl-para-eq10} \nonumber
 	J^2_{X \setminus \bar{G}_2}
		& \leq \int_t^{\infty} \norm{\phi(t^{2m}L) \psi_1(s^{2m}L) \Eins_{B(0,R)} \psi_2(s^{2m}L) \Eins_{X \setminus \bar{G}_2} [e^{-s^{2m}L} \psi(t^{2m}L) g \cdot e^{-s^{2m}L} f]}_{L^2(X)} \,\frac{ds}{s} \\ \nonumber
		& \lesssim \int_t^{\infty} \left(1+\frac{\dist(F,X \setminus \bar{G}_2)^{2m}}{s^{2m}}\right)^{-\gamma} 
				\left(\frac{t}{s}\right)^{2m\delta} \,\frac{ds}{s}  \norm{f}_{L^{\infty}(X)} \norm{g}_{L^2(E)} \\
		& \lesssim  \left(1+\frac{\dist(E,F)^{2m}}{t^{2m}}\right)^{-\gamma} \norm{f}_{L^{\infty}(X)} \norm{g}_{L^2(E)},
\end{align}
since $\delta>\gamma$. Combining \eqref{appl-para-eq7} and \eqref{appl-para-eq8} with \eqref{appl-para-eq9} and \eqref{appl-para-eq10} finishes the proof.
\end{Proof}

\begin{Proof}[of Lemma \ref{Tb-Lemma}]
 Let $b \in L^{\infty}(X)$ and let $t>0$. Further, let $B_1,B_2$ be two arbitrary ball in $X$ with radius $t$ and let $f,g \in L^2(X)$ with $\supp f \subseteq B_1$ and $\supp g \subseteq B_2$.
We decompose the given expression with the help of a Calder\'{o}n reproducing formula as
\begin{align} \label{tb-eq1}
 	\skp{TM_b \psi(t^{2m}L) f,g}
		 = \int_0^{\infty} \skp{\tilde{\psi}(s^{2m}L) M_b \psi(t^{2m}L) f, \psi(s^{2m}L^{\ast}) T^{\ast}g} \,\frac{ds}{s},		
\end{align}
where $\tilde{\psi}$ is the function taken from the assumptions, satisfying $\int_0^{\infty} \psi(t) \tilde{\psi}(t) \, \frac{dt}{t} =1$.
We deduce from Lemma \ref{uniformL2-bound} that due to the weak off-diagonal estimates of order $\gamma>\frac{n}{2m}$, the operator family  $\{T\psi(tL)\}_{t>0}$ is uniformly bounded on $L^2(X)$. Together with assumption \eqref{schur-est} and the Cauchy-Schwarz inequality, this yields 
\begin{align*}
 	\abs{\skp{TM_b \psi(t^{2m}L) f,g}}
			& \leq \int_0^{\infty} \norm{\tilde{\psi}(s^{2m}L) M_b \psi(t^{2m}L) f}_{L^2(X)} \norm{\psi(s^{2m}L^{\ast}) T^{\ast}g}_{L^2(X)} \,\frac{ds}{s} \\
			& \lesssim \int_0^{\infty} \min\left(\frac{s}{t},\frac{t}{s}\right)^{2m\delta} \,\frac{ds}{s} \; \norm{b}_{L^{\infty}(X)} \norm{f}_{L^2(B_1)} \norm{g}_{L^2(B_2)} \\
			& \lesssim \norm{b}_{L^{\infty}(X)} \norm{f}_{L^2(B_1)} \norm{g}_{L^2(B_2)}.
\end{align*}
This shows the desired estimate in the case of $\dist(B_1,B_2)\leq t$.
For $\dist(B_1,B_2)>t$, we split the integral in \eqref{tb-eq1} into two parts, one over $(0,t)$, which is called $J_1$, and one over $(t,\infty)$, which is called $J_2$. \\
To handle $J_1$, we cover $X$ with the help of Lemma \ref{ChristCubes} by balls of radius $t$. That is, we have $X = \bigcup_{\alpha \in I_{k_0}} B_{\alpha}$, where $k_0 \in \Z$ is determined by $C_1 \delta^{k_0} \leq t < C_1 \delta^{k_0-1}$, the balls are defined by $B_{\alpha}:=B(z_{\alpha}^{k_0},t)$ and $I_{k_0}, z_{\alpha}^{k_0}$ are as in Lemma \ref{ChristCubes} and Notation \ref{Cube-Notation}.
Applying this decomposition of $X$ and using the Cauchy-Schwarz inequality, we then get
\begin{align*}
 	\abs{J_1} & \leq \int_0^t \abs{\skp{\tilde{\psi}(s^{2m}L) M_b \psi(t^{2m}L) f, \psi(s^{2m}L^{\ast}) T^{\ast}g}} \,\frac{ds}{s} \\
		& \leq \sum_{\alpha \in I_{k_0}}  \int_0^t  \norm{\tilde{\psi}(s^{2m}L) M_b \psi(t^{2m}L) f}_{L^2(B_{\alpha})}\norm{\psi(s^{2m}L^{\ast}) T^{\ast}g}_{L^2(B_{\alpha})} \,\frac{ds}{s}.
\end{align*}
Due to the weak off-diagonal estimates for $\{\psi(sL^{\ast}) T^{\ast}\}_{s>0}$ and Remark \ref{weak-estimates-st}, we have for all $s<t$ and all $\alpha \in I_{k_0}$
\begin{align} \label{Tb-est1}
 	\norm{\psi(s^{2m}L^{\ast}) T^{\ast}g}_{L^2(B_{\alpha})}
			\lesssim \left(\frac{t}{s}\right)^n \left(1+\frac{\dist(B_{\alpha},B_2)^{2m}}{s^{2m}}\right)^{-\gamma} \norm{g}_{L^2(B_2)}.
\end{align}
On the other hand, as a result of Proposition \ref{H-inf-offdiag}, $\{\tilde{\psi}(tL)\}_{t>0}$ and $\{\psi(tL)\}_{t>0}$ satisfy off-diagonal estimates in $t$ of order $\beta_1$ and $\beta$, respectively. Hence, Lemma \ref{DaviesGaffneyComp} shows that $\{\tilde{\psi}(sL)M_b\psi(tL)\}_{s,t>0}$ satisfies off-diagonal estimates in $\max(s,t)$ of order $\beta=\min(\beta,\beta_1)$. Together with assumption \eqref{schur-est}, this yields for all $s<t$ and all $\alpha \in I_{k_0}$
\begin{align} \label{Tb-est2} \nonumber
 	 & \norm{\tilde{\psi}(s^{2m}L) M_b \psi(t^{2m}L) f}_{L^2(B_{\alpha})} \\ \nonumber
			& \qquad \lesssim \min\left\{\left(1+\frac{\dist(B_1,B_{\alpha})^{2m}}{t^{2m}}\right)^{-\beta},
			 \left(\frac{s}{t}\right)^{2m\delta} \right\} \norm{b}_{L^{\infty}(X)} \norm{f}_{L^2(B_1)} \\
			& \qquad \lesssim \left(1+\frac{\dist(B_1,B_{\alpha})^{2m}}{t^{2m}}\right)^{-\eps\beta}
			 \left(\frac{s}{t}\right)^{(1-\eps)2m\delta}  \norm{b}_{L^{\infty}(X)} \norm{f}_{L^2(B_1)},
\end{align}
for every $\eps \in (0,1)$. 
Recall that we assumed the existence of some $\eps_0 \in (0,1)$ such that $\eps_0 \beta> \frac{n}{2m}$ and $(1-\eps_0)\delta > \frac{n}{2m}+\min(\eps_0\beta,\gamma)$. Since we also assumed $\gamma>\frac{n}{2m}$, we therefore have  $\min(\eps_0\beta,\gamma)>\frac{n}{2m}$. This enables us to apply Lemma \ref{composition-lemma} to get
\begin{align} \label{Tb-est3} 
 	 \sum_{\alpha \in I_{k_0}}  \left(1+\frac{\dist(B_1,B_{\alpha})^{2m}}{t^{2m}}\right)^{-\eps_0\beta} \left(1+\frac{\dist(B_{\alpha},B_2)^{2m}}{s^{2m}}\right)^{-\gamma} 
		\lesssim \left(1+\frac{\dist(B_1,B_2)^{2m}}{t^{2m}}\right)^{-\min(\eps_0\beta,\gamma)},
\end{align}
where we estimated the occuring $s$ in \eqref{Tb-est3} simply by $t$.
 Combining the estimates \eqref{Tb-est1} and \eqref{Tb-est2} with \eqref{Tb-est3}, we therefore obtain
\begin{align*}
 	\abs{J_1}	
		& \leq  \sum_{\alpha \in I_{k_0}}  \int_0^t \norm{\tilde{\psi}(s^{2m}L) M_b \psi(t^{2m}L) f}_{L^2(B_{\alpha})} \norm{\psi(s^{2m}L^{\ast}) T^{\ast}g}_{L^2(B_{\alpha})}  \,\frac{ds}{s} \\
		& \lesssim  \left(1+\frac{\dist(B_1,B_2)^{2m}}{t^{2m}}\right)^{-\min(\eps_0\beta,\gamma)}
			\int_0^t \left(\frac{s}{t}\right)^{(1-\eps_0)2m\delta-n}  \,\frac{ds}{s} \norm{b}_{L^{\infty}(X)} \norm{f}_{L^2(B_1)} \norm{g}_{L^2(B_2)} \\
		& \lesssim \left(1+\frac{\dist(B_1,B_2)^{2m}}{t^{2m}}\right)^{-\min(\eps_0\beta,\gamma)}
			 \norm{b}_{L^{\infty}(X)} \norm{f}_{L^2(B_1)} \norm{g}_{L^2(B_2)},
\end{align*}
where in the last step we used the fact that the integral is bounded by a constant independent of $s$ and $t$ due to the assumption $(1-\eps_0)\delta > \frac{n}{2m}$. Therefore, the last line gives the desired estimate for $J_1$. \\
We now turn to the integral $J_2$. As before, we cover $X$ with balls of radius $t$ and use the Cauchy-Schwarz inequality to get
\begin{align*}
 	\abs{J_2} & \leq \int_t^{\infty} \abs{\skp{\tilde{\psi}(s^{2m}L) M_b \psi(t^{2m}L) f, \psi(s^{2m}L^{\ast}) T^{\ast}g}} \,\frac{ds}{s} \\
		& \leq \sum_{\alpha \in I_{k_0}}  \int_t^{\infty}  \norm{\tilde{\psi}(s^{2m}L) M_b \psi(t^{2m}L) f}_{L^2(B_{\alpha})} \norm{\psi(s^{2m}L^{\ast}) T^{\ast}g}_{L^2(B_{\alpha})}\,\frac{ds}{s}.
\end{align*}
On the one hand, we again use the weak off-diagonal estimates for $\{\psi(sL^{\ast})T^{\ast}\}_{s>0}$, applied to balls of radius $t$ by embedding them into larger balls of radius $s$, and get for $s>t$ and $\alpha \in I_{k_0}$
\begin{align} \label{Tb-est4}
 	 \norm{\psi(s^{2m}L^{\ast}) T^{\ast}g}_{L^2(B_{\alpha})}
			\lesssim \left(1+\frac{\dist(B_{\alpha},B_2)^{2m}}{s^{2m}}\right)^{-\gamma} \norm{g}_{L^2(B_2)}.
\end{align}
In analogy to \eqref{Tb-est2}, on the other hand, we obtain by application of Proposition \ref{H-inf-offdiag}, Lemma \ref{DaviesGaffneyComp} and assumption \eqref{schur-est} for every $s>t$ and every $\alpha \in I_{k_0}$
\begin{align} \label{Tb-est5} 
 	  \norm{\tilde{\psi}(s^{2m}L) M_b \psi(t^{2m}L) f}_{L^2(B_{\alpha})} 
			 \lesssim \left(1+\frac{\dist(B_1,B_{\alpha})^{2m}}{s^{2m}}\right)^{-\eps_0\beta}
			 \left(\frac{t}{s}\right)^{(1-\eps_0)2m\delta}  \norm{b}_{L^{\infty}(X)} \norm{f}_{L^2(B_1)}.
\end{align}
Lemma \ref{composition-lemma} in turn yields that for all $s>t$
\begin{align} \label{Tb-est6} \nonumber
 	& \sum_{\alpha \in I_{k_0}} 
		\left(1+\frac{\dist(B_1,B_{\alpha})^{2m}}{s^{2m}}\right)^{-\eps_0\beta} \left(1+\frac{\dist(B_{\alpha},B_2)^{2m}}{s^{2m}}\right)^{-\gamma} \\
		& \qquad \lesssim \left(\frac{s}{t}\right)^{(n+\eps)} \left(1+\frac{\dist(B_1,B_2)^{2m}}{s^{2m}}\right)^{-\min(\eps_0\beta,\gamma)}
\end{align}
for arbitrary $\eps>0$. 
As above, the combination of \eqref{Tb-est4}, \eqref{Tb-est5} and \eqref{Tb-est6} provides us with
\begin{align} \label{Tb-est7} \nonumber 
 	 \abs{J_2} 
			&  \leq   \sum_{\alpha \in I_{k_0}}  \int_t^{\infty}  \norm{\tilde{\psi}(s^{2m}L) M_b \psi(t^{2m}L) f}_{L^2(B_{\alpha})}\norm{\psi(s^{2m}L^{\ast}) T^{\ast}g}_{L^2(B_{\alpha})} \,\frac{ds}{s} \\ \nonumber
			&  \lesssim \int_t^{\infty} \left(1+\frac{\dist(B_1,B_2)^{2m}}{s^{2m}}\right)^{-\min(\eps_0\beta,\gamma)}
								\left(\frac{t}{s}\right)^{(1-\eps_0)2m\delta-(n+\eps)} \,\frac{ds}{s} \\
			& \qquad \qquad \times	\norm{b}_{L^{\infty}(X)} \norm{f}_{L^2(B_1)} \norm{g}_{L^2(B_2)}.
\end{align}
Finally observe that the integral in \eqref{Tb-est7} can in view of the assumption $\dist(B_1,B_2)>t$ be bounded by a constant times
\begin{align*}
	&\left(\frac{\dist(B_1,B_2)^{2m}}{t^{2m}}\right)^{-\min(\eps_0\beta,\gamma)}
			\int_t^{\infty} \left(\frac{t}{s}\right)^{-2m\min(\eps_0\beta,\gamma)} \left(\frac{t}{s}\right)^{(1-\eps_0)2m\delta-(n+\eps)} \,\frac{ds}{s} \\
	& \qquad \lesssim \left(1+\frac{\dist(B_1,B_2)^{2m}}{t^{2m}}\right)^{-\min(\eps_0\beta,\gamma)},
\end{align*}
since we postulated  $(1-\eps_0)\delta-\min(\eps_0\beta,\gamma)>\frac{n+\eps}{2m}$ for sufficiently small $\eps>0$.
\end{Proof}

\begin{Proof}[of Lemma \ref{Carleson-Tb}]
We set $M:=[\frac{n}{4m}]+1$. Then $BMO_L(X)=BMO_{L,M}(X)$ according to Definition \ref{DefBMO-uniform}.\\
We follow the proof of \cite{HofmannMayboroda}, Lemma 8.3, replacing the operator family $\{(tL)^Me^{-tL}\}_{t>0}$ by the operator family $\{\psi(tL)M_b\}_{t>0}$. The corresponding term $I_1$ can be handled with just the same methods, once one has checked that $\{\psi(tL)M_b\}_{t>0}$ satisfies quadratic estimates and off-diagonal estimates of order $\beta>\frac{n}{4m}$ and that this are the only properties of $\{(tL)^Me^{-tL}\}_{t>0}$ that are used in  \cite{HofmannMayboroda}, Lemma 8.3, for $I_1$.\\ 
For the term $I_2$, it needs a more careful treatment. What is essential for this part is the fact that the operator family $\{(tL)^Me^{-tL}(tL)^{-k}\}_{t>0}$, now replaced by $\{\psi(tL)M_b (tL)^{-k}\}_{t>0}$, satisfies off-diagonal estimates of order $\beta-k>\frac{n}{4m}$ for every $1 \leq k \leq M$. If one can establish these estimates, the proof for the second part $I_2$ can be copied from the one of \cite{HofmannMayboroda}, Lemma 8.3.\\
Thus, let us show, in analogy to Lemma \ref{Tb-Lemma}, that $\{\psi(tL)M_b (tL)^{-k}\}_{t>0}$ satisfies off-diagonal estimates of some order larger than $\frac{n}{4m}$. 
Let $E,F$ be two open sets in $X$ and let $g \in \calD(L^{-k})$ with $\supp g \subseteq E$, $h \in L^2(X)$ with $\supp h \subseteq F$.  Via the Calder\'{o}n reproducing formula, we write
\begin{align*}
		  \skp{\psi(t^{2m}L)M_b (t^{2m}L)^{-k} g,h} 
				 = \int_0^{\infty} \skp{\psi(t^{2m}L)M_b \psi(s^{2m}L) \tilde{\psi}(s^{2m}L) (t^{2m}L)^{-k} g,h}  \,\frac{ds}{s}.
\end{align*}
Due to the Cauchy-Schwarz inequality, the uniform boundedness of $\{\psi(sL)(sL)^{-k}\}_{s>0}$ and assumption \eqref{schur-est} we then obtain 
\begin{align*}
	& \abs{\skp{\psi(t^{2m}L)M_b (t^{2m}L)^{-k} g,h} } \\
				& \qquad \leq \int_0^{\infty} \left(\frac{s}{t}\right)^{2mk} 
								\norm{\psi(s^{2m}L)(s^{2m}L)^{-k} g}_{L^2(X)} \norm{\tilde{\psi}(s^{2m}L^{\ast}) M_{\bar{b}} \psi(t^{2m}L^{\ast}) h}_{L^2(X)}\,\frac{ds}{s} \\
				& \qquad \lesssim \int_0^{\infty} \min\left(\frac{s}{t},\frac{t}{s}\right)^{2m\delta} \left(\frac{s}{t}\right)^{2mk} \,\frac{ds}{s} \norm{b}_{L^{\infty}(X)} \norm{g}_{L^2(X)} \norm{h}_{L^2(X)} \\
				& \qquad \lesssim \norm{b}_{L^{\infty}(X)} \norm{g}_{L^2(E)} \norm{h}_{L^2(F)},
\end{align*}
where for the case $s>t$ we take into account that $\delta>M$ and therefore $\delta>k$ for all $1 \leq k \leq M$.
This yields the desired estimate for $\dist(E,F)\leq t$.\\
For the case $\rho:=\dist(E,F)>t$, we define the sets $G_1:=\{x \in X\,:\, \dist(x,F) < \frac{\rho}{2}\}$ and $G_2:=\{x \in X\,:\, \dist(x,F) < \frac{\rho}{4}\}$ and then split $X$ into $X = \bar{G}_2 \cup X \setminus \bar{G}_2$. By construction $G_1,G_2$ are open with $\dist(E, G_1) \geq \frac{\rho}{2}$ and $\dist(F, X \setminus \bar{G}_2) \geq \frac{\rho}{4}$.
Using that $\bar{G}_2 \subseteq G_1$, this leads to 
\begin{align*}
		& \abs{\skp{\psi(t^{2m}L)M_b (t^{2m}L)^{-k} g,h} } \\
				& \quad \leq \int_0^{\infty} \left(\frac{s}{t}\right)^{2mk} 
								\norm{\psi(s^{2m}L)(s^{2m}L)^{-k} g}_{L^2(G_1)} \norm{\tilde{\psi}(s^{2m}L^{\ast}) M_{\bar{b}} \psi(t^{2m}L^{\ast}) h}_{L^2(G_1)}\,\frac{ds}{s}
\end{align*}
\begin{align*}
				& \qquad + \int_0^{\infty} \left(\frac{s}{t}\right)^{2mk} 
								\norm{\psi(s^{2m}L)(s^{2m}L)^{-k} g}_{L^2(X \setminus \bar{G}_2)} \norm{\tilde{\psi}(s^{2m}L^{\ast}) M_{\bar{b}} \psi(t^{2m}L^{\ast}) h}_{L^2(X\setminus \bar{G}_2)}\,\frac{ds}{s} \\
				& \quad =: J_{G_1} + J_{X\setminus \bar{G}_2}.
\end{align*}
For the term $J_{X \setminus \bar{G}_2}$ we get via Lemma \ref{DaviesGaffneyComp}, applied to $\{\tilde{\psi}(sL^{\ast})M_{\bar{b}} \psi(tL^{\ast})\}_{s,t>0}$, assumption \eqref{schur-est} and the uniform boundedness of  $\{\psi(sL)(sL)^{-k}\}_{s>0}$ 
\begin{align} \label{Carleson-Tb-eq4} \nonumber
	J_{X \setminus \bar{G}_2}
			& \lesssim \int_0^{\infty} \min\left(\left(1+\frac{\dist(F,X \setminus \bar{G}_2)^{2m}}{\max(s,t)^{2m}}\right)^{-\beta}, \min\left( \frac{s}{t}, \frac{t}{s} \right)^{2m \delta}\right) \left(\frac{s}{t}\right)^{2mk}  \,\frac{ds}{s} \\
				& \qquad  \times \norm{b}_{L^{\infty}(X)} \norm{g}_{L^2(E)} \norm{h}_{L^2(F)}.
\end{align}
Since by construction $\dist(F, X \setminus \bar{G}_2) \gtrsim \dist(E,F)>t$, we can bound the integral in \eqref{Carleson-Tb-eq4} in a similar way as in the proof of Lemma \ref{Tb-Lemma} by a constant times
\begin{align*}
		& \int_0^t \min\left(\left(1+\frac{\dist(E,F)^{2m}}{t^{2m}}\right)^{-\beta}, \left(\frac{s}{t}\right)^{2m \delta}\right) \left(\frac{s}{t}\right)^{2mk}  \,\frac{ds}{s} \\
		&  \qquad + \int_t^{\infty} \min\left(\left(1+\frac{\dist(E,F)^{2m}}{s^{2m}}\right)^{-\beta},\left( \frac{t}{s} \right)^{2m \delta}\right) \left(\frac{s}{t}\right)^{2mk}  \,\frac{ds}{s} \\
		&  \lesssim \left(1+\frac{\dist(E,F)^{2m}}{t^{2m}}\right)^{-\beta} \int_0^t\left(\frac{s}{t}\right)^{2mk}  \,\frac{ds}{s} \\
		&  \qquad + \left(\frac{\dist(E,F)^{2m}}{t^{2m}}\right)^{-\eps_0\beta} \int_t^{\infty} \left(\frac{t}{s} \right)^{-2m\eps_0\beta} \left( \frac{t}{s} \right)^{2m(1-\eps_0) \delta- 2mk} \frac{ds}{s} \\
		& \lesssim \left(1+\frac{\dist(E,F)^{2m}}{t^{2m}}\right)^{-\eps_0\beta}
\end{align*}
for $\eps_0 \in (0,1)$ as given in the assumptions with $(1-\eps_0)\delta>\eps_0 \beta+k$ for all $1 \leq k \leq M$. \\
It remains to estimate $J_{G_1}$. Observe that $\{\psi(sL)(sL)^{-k}\}_{s>0}$ satisfies off-diagonal estimates of order $\beta-k$ due to Proposition \ref{H-inf-offdiag}. With the help of assumption \eqref{schur-est}, we therefore obtain
\begin{align} \label{Carleson-Tb-eq5} \nonumber
	J_{G_1} & \lesssim \int_0^{\infty} \left(1+\frac{\dist(E,G_1)^{2m}}{s^{2m}}\right)^{-(\beta-k)} \min\left( \frac{s}{t}, \frac{t}{s} \right)^{2m \delta} \left(\frac{s}{t}\right)^{2mk} \,\frac{ds}{s}\\
				& \qquad  \times \norm{b}_{L^{\infty}(X)} \norm{g}_{L^2(E)} \norm{h}_{L^2(F)}.
\end{align}
Using the fact that $\dist(E,G_1) \gtrsim \dist(E,F)>t$ and the assumption $\delta>\beta$, we can show that the integral in \eqref{Carleson-Tb-eq5} is bounded by a constant times
\begin{align*}
		&  \left(1+\frac{\dist(E,F)^{2m}}{t^{2m}}\right)^{-(\beta-k)} \int_0^t  \left( \frac{s}{t} \right)^{2m \delta} \left(\frac{s}{t}\right)^{2mk} \,\frac{ds}{s}\\
		&  \qquad + \left(\frac{\dist(E,F)^{2m}}{t^{2m}}\right)^{-(\beta-k)} \int_t^{\infty} \left(\frac{t}{s}\right)^{-2m(\beta-k)} \left(\frac{t}{s} \right)^{2m \delta}\left(\frac{s}{t}\right)^{2mk} \,\frac{ds}{s}\\
		&  \lesssim  \left(1+\frac{\dist(E,F)^{2m}}{t^{2m}}\right)^{-(\beta-k)}.
\end{align*}
In summary, the above estimates yield that the operator family $\{\psi(tL)M_b (tL)^{-k}\}_{t>0}$ satisfies off-diagonal estimates of order $\min(\beta-k,\eps_0\beta)>\frac{n}{4m}$ for every $1 \leq k \leq M$.
\end{Proof}

\addcontentsline{toc}{section}{References}
\small{

}

\small{\textsc{Institute for Analysis, Karlsruhe Institute of Technology (KIT), Kaiserstr. 89, D-76128 Karlsruhe, Germany} \\ \textit{E-mail address:} \texttt{dorothee.frey@kit.edu, peer.kunstmann@kit.edu}}

\end{document}